\documentclass[11pt]{article}

\usepackage{pdfpages}
\usepackage{multirow,multicol,soul}
\usepackage[normalem]{ulem} 
\usepackage{tabularx}
\usepackage{longtable,eurosym}
\usepackage{booktabs}
\usepackage{amsmath,bbm,amsfonts,latexsym,}
\usepackage{amssymb,epsfig}
\usepackage{float,graphicx,subfigure}
\usepackage[utf8]{inputenc}
\usepackage[counterclockwise]{rotating}
\usepackage{lscape,setspace}
\usepackage{url}
\usepackage{color}
\usepackage[T1]{fontenc}
\usepackage[square,sort,comma,numbers]{natbib}
\setcitestyle{authoryear,open={(},close={)}}
\usepackage[scale={0.75,0.9},centering,includeheadfoot]{geometry}
\usepackage{graphicx,float}
\usepackage[colorlinks=true,linkcolor=blue,citecolor=blue,urlcolor=blue]{hyperref}
\usepackage{enumerate}
\usepackage{comment}
\usepackage{tikz}
\usepackage{longtable}

\DeclareMathOperator*{\argmin}{arg\,min}

\newcommand{\calB}{{\cal B}}

\newcommand{\calH}{{\cal H}}

\newcommand{\calS}{{\cal S}}

\newcommand{\eps}{{\epsilon}}

\newcommand{\la}{\lambda}

\numberwithin{equation}{section}

\def\E{{\mathbb{E}}}

\def\eps{\epsilon}
\newcommand{\half}{\frac{1}{2}\:}

\newcommand{\one}{\mathbbm{1}}
\newcommand{\real}{{\mathbb{R}}}
\newcommand{\pr}{{\mathbb{P}}}
\newcommand{\qed}{$\Box$}

\newtheorem{theorem}{Theorem}[section]

\newtheorem{lemma}[theorem]{Lemma}

\newtheorem{proposition}[theorem]{Proposition}

\newtheorem{definition}[theorem]{Definition}

\newtheorem{rem}[theorem]{Remark}

\newtheorem{Remark}[theorem]{Remark}

\newtheorem{properties}[theorem]{Properties}
\newtheorem{Example}[theorem]{Example}

\graphicspath{{plots/}}

\begin{document}

\title{
Extremal behaviour and convergence rates for sample--based geometric quantiles and half space depths \\[1ex]
}


\author{Sibsankar SINGHA\textsuperscript{(a)} \\[.5ex]
Supervisors: Marie KRATZ\textsuperscript{(b)} and Sreekar VADLAMANI\textsuperscript{(c)}
\\[1ex]
\small
\textsuperscript{(a)} TIFR-CAM, Bangalore, India \& ESSEC CREAR, France; Email: sibsankar@tifrbng.res.in \\
\small
\textsuperscript{(b)} ESSEC Business School, CREAR, Cergy-Pontoise, France; Email: kratz@essec.edu \\
\small
\textsuperscript{(c)} TIFR-CAM, Bangalore, India; Email: sreekar@tifrbng.res.in
}

\maketitle

\begin{abstract}
   We consider the empirical versions of geometric quantile and halfspace depth, and study their extremal behaviour as a function of the sample size. The objective of this study is to establish connection between the rates of convergence and tail behaviour of the corresponding underlying distributions. The intricate interplay between the sample size and the parameter driving the extremal behaviour forms the main result of this analysis. In the process, we also fill certain gaps in the understanding of population versions of geometric quantile and halfspace depth.
\end{abstract}
       
 
\newpage
    
\tableofcontents
\newpage

\section{Introduction}



Due to recent advances in high-dimensional statistics, there is a renewed interest in developing tools to better understand the geometric structure of datasets. Numerous multivariate quantiles and statistical depth functions have been proposed to establish ranks and identify outliers in multivariate data. While quantiles are defined analytically using the inverse of the cumulative distribution function, depth functions take a more geometric approach, employing halfspaces, paraboloids, and projections to measure centrality from a global perspective. This results in an ordering of observations from the center outward. Consequently, the philosophies behind quantiles and depth functions appear to be distinct. However, both concepts offer a geometric perspective on ordering in a multivariate setup. In fact, in the case of univariate data, quantiles and depth functions are conceptually related through a functional relationship, making them inseparable.

These geometric tools of depths and quantiles offer nonparametric descriptions of a data set in a multidimensional space, making them quite useful for statistical inference problems (e.g. \cite{Serfling2006}), among which classification and regression (see e.g. \cite{Cuevas2007,Dutta2016,Hubert2010,Hallin2010,Paindaveine2012,Rousseeuw2004,Struyf1999}), for learning theory (see e.g. \cite{Bousquet2002,Koltchinskii2002,Koltchinskii2006,Koltchinskii2011,Panchenko2003,Gine2006}), for bootstrap (e.g. \cite{Cuevas2006}), outliers or anomaly detection (e.g. \cite{Staerman2022}), applications to geometry (e.g. \cite{Kong2012}) and multivariate risk analysis. 

Numerous depth functions have been introduced and studied, starting with Mahalanobis distance depths (\cite{Mahalanobis1936,Liu1993,Zuo2000}), the well-known and used Tukey or halfspace depth (\cite{Tukey1975}), going on, for instance, with simplicial (volume) depths (\cite{Oja1983,Liu1990}), onion depths (\cite{Barnett1976,Eddy1982}), all notions of spatial depths (\cite{Dudley1992,Chaudhury1996,Koltchinskii1997,Vardi2000,Mottonen2005}), the projection depth (\cite{Donoho1992,Zuo2003,Dutta2012,Nagy2020}), the zonoid depth (\cite{Dyckerhoff1996,Koshevoy1997,Koshevoy2003}), 
local depths (\cite{Agostinelli2011,Paindaveine2013}). We refer to \cite{Hallin2010,Mosler2002,Mosler2013,Kuelbs2016,Chernozhukov2017,Nagy2020,Nagy2021,Mosler2021} for theoretical and practical aspects (as well as computational) of depth functions, as well as to \cite{Nagy2022} and references therein for halfspace depths.


Our motivation lies in exploring asymptotics, specifically the examination of the behavior of depth-based multivariate quantiles as they approach extreme regions, both in terms of population measures and empirical data. Furthermore, our objective is to comprehend the connection between the extreme behavior of a probability measure (whether it exhibits a light or heavy tail) and the corresponding geometric measures associated with it.


We focus on two prominent geometric measures: geometric quantiles, introduced by \cite{Chaudhury1996}, and halfspace depth, as described by \cite{Tukey1975}. The selection of these specific geometric measures stems from the core objective of our research. The fundamental question of characterising the tail behaviour of a probability measure using these geometric measures inherently addresses whether they capture essential aspects of the underlying probability distribution. In fact, it was demonstrated by \cite{Koltchinskii1997} and subsequently utilized by \cite{Dhar2014} for proposing a test statistic, that geometric quantiles uniquely identify the underlying probability measure. See also \cite{Konen2023} for details on analytical inversion of geometric quantiles to recover their underlying probability measure. However, the same does not hold true for halfspace depth, as shown by \cite{Nagy2021a}. On the other hand, \cite{Struyf1999} established that halfspace depth uniquely identifies measures with finite support (e.g., empirical measures). This uniqueness property (under constraint for halfspace depths) indicates a direct correspondence between these two geometric measures of our interest and the underlying probability measures. It is therefore natural to look for clearer connection between the extremal behaviours of these geometric measures and their underlying probability measures. This motivates our study.


Numerous results are already available for the population-based analysis of these geometric measures, particularly in terms of their asymptotic (extreme) behaviour. For geometric quantiles, we refer to \cite{Girard2017, Paindaveine2021}. For halfspace depth, assuming multivariate regularly varying distributions, the work by \cite{He2017} is noteworthy. 


Considering practical applications, the questions regarding asymptotics become even more critical when examining sample versions of these two geometric measures. This forms the essence of the paper: We establish convergence rates for the sample versions and investigate the extreme behavior of the geometric measures based on the nature of the underlying distribution

Concerning geometric quantiles, proofs are developed with classical tools of probability, solving the paradoxical issue raised by \cite{Girard2017} in the population setting, and using results by \cite{Chaudhury1996,Paindaveine2021} and Glivenko--Cantelli theorem for the sample version.

In the case of halfspace depth, recall that (see \cite{Donoho1992}), as the sample size increases, the halfspace depth for a sample converges almost surely to the halfspace depth for the underlying distribution. To obtain rates of decay of halfspace depth, our approach builds mainly on the theory of empirical processes (\cite{Shorack2009}) and weighted empirical processes indexed by sets by \cite{Alexander1987}. The latter paper and its powerful results, helped us prove our own results in a rather direct and elegant way. 


In the scenario where the tail of the distribution is light, \cite{Burr2017} employed a novel geometric methodology to derive uniform convergence rates for halfspace depth. Through the reorganization of halfspaces into one-dimensional families, the authors successfully attained improved convergence bounds for the sample version of halfspace depth, surpassing the typical Glivenko-Cantelli bounds. This advancement was particularly evident when considering exponential decay in the underlying distribution.

Our approach allows one to establish convergence results for both light and heavy tails. 

Our main results are illustrated for the population and sample versions. The geometric quantiles have been programmed on python using the algorithm developed by \cite{Dhar2014}, while the computation of halfspace depth uses R-packages developed in \cite{Liu2019,Pokotylo2019} to evaluate Tukey depth and contours. In this latter case, the method relies on an approximation of the true depth proposed by \cite{Dyckerhoff2004} (see also \cite{Dyckerhoff2021}), for which rates of convergence have been obtained by \cite{Nagy2020,Mosler2021}. We note here that the computation of depths is indeed quite a challenge, taken by various research teams to develop relevant softwares (see also e.g. \cite{Genest2019,Liu2015,Mahalanobish2015}). 

{\it Notation.} All the analysis in this paper is on $\real^d$, unless otherwise stated. The centered unit open ball and the unit sphere in $\real^d$ are denoted by $B^d$ and $S^{d-1}$, while $\langle\cdot,\cdot\rangle$ and $\|\cdot\|$, denote the Euclidean inner product and $l^2$-norm, respectively, in $\real^d$.


{\it Structure of the paper.} The paper primarily consists of two main sections,
where we explore the asymptotic properties of the 
geometric quantiles and halfspace depths. 
Section~\ref{sec:multi-quantiles} considers the population versions of the two geometric measures, completing the literature on the topic.
In the third section, core of the paper, we study their empirical counterparts. 
In particular, we investigate the asymptotic behaviour of their sample versions in relation to the sample size. We provide illustrative examples using samples drawn from bivariate distributions with both light and heavy tails. Concluding remarks are stated in Section~\ref{sec:concl}. The proofs for all the results are presented in Section~\ref{sec:proofs}.


\section{Multivariate geometric measures} 
\label{sec:multi-quantiles}

As in the univariate case, it is expected that multivariate quantiles or depth functions encode the tail behaviour of the underlying probability measure. It is this line of thought that we explore in this section, studying the asymptotic behaviour of geometric quantiles and halfspace depth {\it vis-\`{a}-vis} of the tail characteristics of the underlying probability measure, while also discussing the {\it status quo} of the subject to contextualise our results.

\subsection{Geometric quantiles}
\label{ss:geomQ}

The motivation for geometric quantiles can be traced to univariate notion of quantile, which are defined as generalized inverse of cumulative distribution function. However, such definition does not have a natural extension to the multivariate case: firstly, a cumulative distribution does not have a clear interpretability in higher dimensions and, secondly, the generalized inverse of a cumulative distribution function in higher dimensions would be a subset of the domain, and not a point as in the univariate case, again raising the concerns of interpretability.

It is well known that in the univariate case, the $\alpha$-th quantile of $\pr$ is defined as  (see \cite{Ferguson1967})
$$
\argmin_{q \in \real} \left[\int \left(|x-q|-|x|\right)d\pr(x) -(2\alpha -1)q\right].
$$ 
This univariate representation of quantiles forms the basis for the following generalisation.

\subsubsection{State of the art: Definitions and Properties}
\label{sss:knownGeomQ}

\begin{definition}[Geometric quantile \cite{Chaudhury1996}]
Let $X$ be a $\real^d$- valued random vector for $d \geq 1$, such that the induced measure $\pr$ is not supported on any straight line (or any of its subsets), and $u \in \real ^d$ be such that $\|u\|<1$. The geometric quantile $q_X(u)$ is defined as the solution of the following optimisation problem:
    \begin{equation}
    \label{eqn:geometric-quantile-dfn}
	    \argmin_{q \in \real^d}\{\E\big[\|X-q\|-\|X\|\big]-\langle u,q \rangle\}.
    \end{equation}
\end{definition}
In the literature, $u$ is called the {\it index vector} that controls the centrality of the quantile: When $u$ is close to the origin, the corresponding quantiles are close to the center of the distribution and referred to as central quantiles, and, when $u$ is close to the boundary of the unit ball, then the corresponding quantiles are referred to as extreme quantiles. It is noteworthy that the minimiser of the optimisation problem in \eqref{eqn:geometric-quantile-dfn} has a unique solution for $d\ge 2$ whenever the support of $\pr$ is not contained in a single straight line.
Writing 
$$
\psi(u,q)= \E\big[\|X-q\|-\|X\|\big]-\langle u,q \rangle,
$$ 
the existence and uniqueness of the solution of the optimization problem \eqref{eqn:geometric-quantile-dfn} for every fixed $u$ is a consequence of the strict convexity and continuity of $\psi$, together with the fact that $\psi \to \infty$ as $\|q\| \to \infty$. 

Note that the argument for uniqueness breaks down in $d\ge 2$ whenever the distribution of $X$ is supported on a single straight line. However, when $d=1$, the geometric quantile coincides with the usual univariate quantile defined as generalised inverse of a cumulative distribution function. In this case, the optimisation problem \eqref{eqn:geometric-quantile-dfn} may not have a unique minimiser, since the distribution lies on a straight line. We refer the reader to \cite{Chaudhury1996} for all detailed proofs and arguments.

As is the case with most optimisation problems, getting closed form expressions for the solution of \eqref{eqn:geometric-quantile-dfn} is very difficult, even for simple probability distributions. Therefore, one often looks for characterising properties of the solution in order to gain insight into the solution. We shall list below some known properties of geometric quantiles for non--atomic distributions whose support is not contained in any unit dimensional affine subspace. We refer the readers to the original references for the precise statements and proofs. 
%
\begin{properties} \label{property:set-1} 
Some known characteristic properties of geometric quantiles from \cite{Chaudhury1996} and \cite{Koltchinskii1997}: 
\begin{itemize}
    \item From \cite{Chaudhury1996}
    \begin{itemize}
        \item[(i)] The quantile map $q_X:u \mapsto q_X(u)$ defined in \eqref{eqn:geometric-quantile-dfn} is an homeomorphism from $B^d(0)$ to $\real^d$. (See also \cite{Koltchinskii1997}.)     
        \item[(ii)] For any fixed $u\in B^d(0)$, there exists a unique $q\in \real ^d$ such that \eqref{eqn:geometric-quantile-dfn} holds iff the pair $(u,q)\in \real^d \times \real^d$ satisfies 
        \begin{equation} \label{eqn: exp of sign}
            u=-\E\Bigg[\frac{X-q}{\|X-q\|}\Bigg]
        \end{equation}
        This characterisation of geometric quantiles is often more useful for analysis. 
        \item[(iii)] \label{prop:homogenious}
        The quantile map is positive homogeneous, i.e. $q_{c\, X}(u)=c \, q_X(u)$ for any positive constant $c > 0$.
        \item[(iv)] \label{prpn:translation invariance}
        For any $C \in \real^d$, $q_{X+C}(u)=q_X(u)+C$. 
        \item[(v)] The geometric quantile is rotational equivariant, i.e. $q_{AX}=A\,q_{X}$ for any orthogonal matrix $A$.
        \end{itemize}
    \item From \cite{Koltchinskii1997}
    \begin{itemize}
        \item[(vi)] Geometric quantiles uniquely characterise the distribution. (See also \cite{Chaudhury1996} for a discussion.)
        \end{itemize}
    \end{itemize}
\end{properties}

Additionally, more analytical results can be extracted by imposing further restrictions on the underlying distribution. Specifically: 
\begin{properties}[\textnormal{From \cite{Girard2017}}]\label{property:set-2}
Assume that the distribution of $X$ is rotationally symmetric about the origin. Then, 
\begin{itemize}
	\item[(i)] The quantile map $u \mapsto q_X(u)$ commutes with every linear isometry of $\real^d$. The norm of $q_X(u)$ only depends upon the norm of $u$. As a result, the iso--quantile contours are spheres centered at the center of symmetry.
	\item[(ii)] The directions of the quantile $q_X(u)$ and the index vector $u$ are the same if $u\neq 0$. 
	\item[(iii)] The function $\|u\|\to \|q_X(u)\|$ is a continuous strictly increasing function on $[0,1)$.
	\item[(iv)] $\|q_X(u)\| \to \infty$ as $\|u\|\to 1$.
	\item[(v)] $\frac{q_X(v)}{\|q_X(v)\|} \to u$, as $v \to u \in S^{d-1}, \,\, \|v\|< 1$.
\end{itemize}
\end{properties}

\begin{figure}[H]
    \centering
    \includegraphics[width=.5\linewidth]{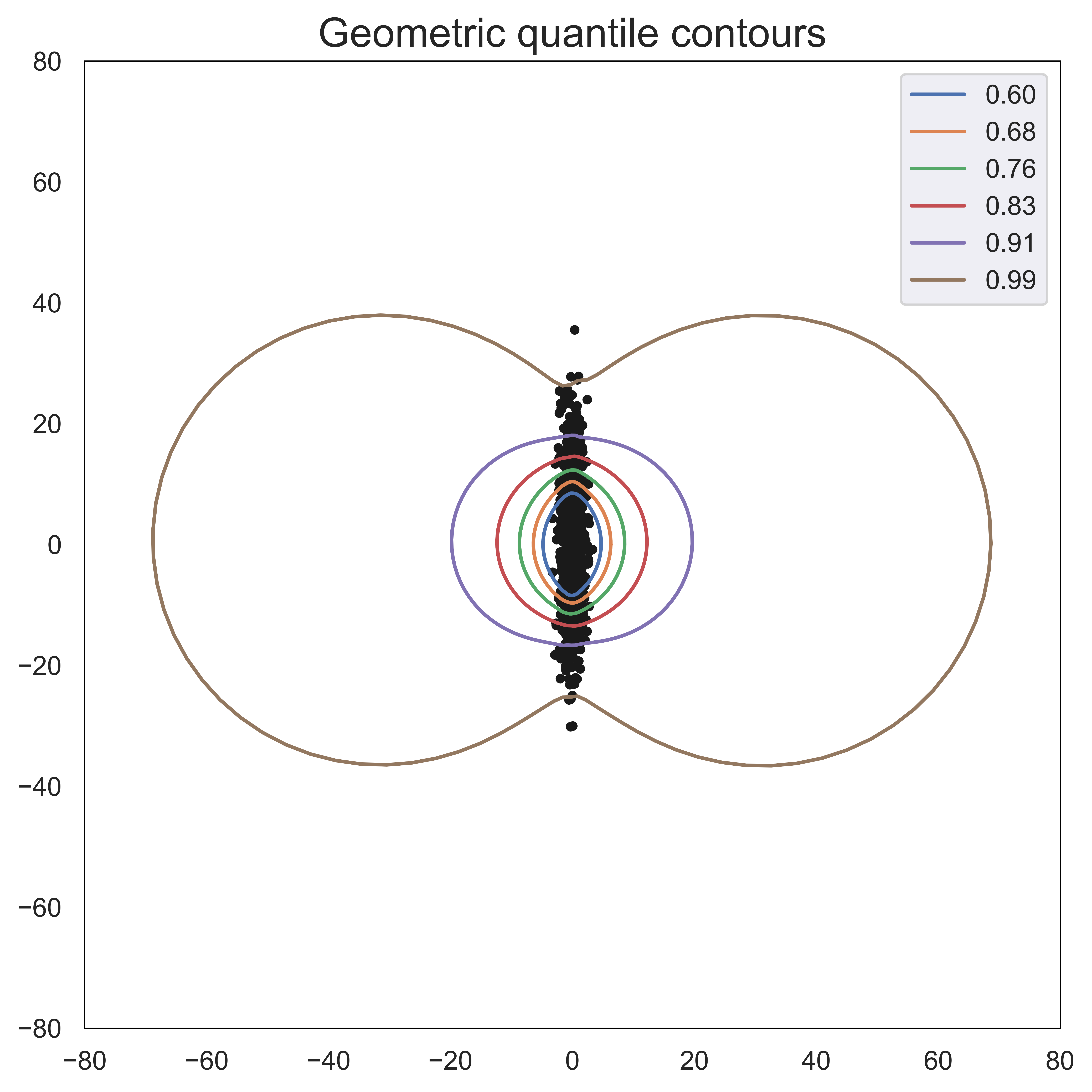}
    \parbox{380pt}{\caption{\sf Samples of size $1000$ are drawn from a Gaussian distribution with mean zero and diagonal covariance matrix $diag(1,100)$. Observe that the convexity present in the underlying sampling distribution is not reflected in the isoquantile contours.}
    \label{fig:GeometricQuantileContours}}
\end{figure}

\begin{rem}
\leavevmode
\begin{itemize}
\item[(i)] Similar results can be obtained when the distribution is assumed to be symmetric around a point other than the origin.
\item[(ii)] Quantile contours, even for elliptically symmetric distributions, are not convex; see Figure~\ref{fig:GeometricQuantileContours}. In particular, quantile contours are not convex for non-spherical distributions. 
\item[(iii)] Some asymptotic properties stated above for spherically symmetric distributions hold true even for non-symmetric distributions. For instance, $\|q(u)\| \to \infty$ as $\|u\|\to 1$ for any distribution -- even atomic; see \cite[Theorem 2.1]{Girard2017}. A simple and interesting consequence of this property is that the quantile contours do not remain contained within the convex hull of the support of the distribution. 
\item[(iv)] Additionally, the asymptotic direction of a geometric quantile matches that of its index vector -- for any distribution including the atomic ones; see \cite[Theorem 2.1]{Girard2017}.
\end{itemize}
\end{rem}
\vspace{-4ex}

\subsubsection{Discriminating tail behaviours}\label{sss:newGeomQ}

We know that all geometric quantiles increase to infinity as the index vector approaches the boundary of unit ball. A natural question is then to investigate the rate of increase of these quantiles for different distributions.  
While in the univariate case, the asymptotic behaviour of quantiles characterises the tail behaviour of the distribution, we may want to question this property for multivariate quantiles.
A first answer has been given in \cite[Theorem 2.2]{Girard2017}, where the authors analyse the asymptotic behaviour of geometric quantiles for non atomic distributions by introducing an appropriate normalisation, obtaining a {\it first order expansion} of both magnitude and direction of $q_X(u)$ as $u$ approaches the boundary of the unit ball.
The authors also observe that this asymptotic behaviour appears counter intuitive as the extreme quantiles appear to be dependent on {\it central parameters} like covariance. In particular, if two distributions have the same covariance matrix, then their {\it first order} asymptotic expansion of quantiles are identical. 

But, since geometric quantiles uniquely characterize the underlying probability measure, we need to address this issue.
Hence, we further investigate the asymptotic behaviour of geometric quantiles. A natural first approach is simply to look for higher order asymptotic expansions that may help distinguish between distributions. 
Here, we focus mainly on second order expansions of both magnitude and direction of $q(\alpha u)$. Similar higher order asymptotic expansions can be worked out with additional moment conditions. It must be noted that, in this way, we obtain the rate of convergence, as $\alpha \to 1$, of  $\|q_X(\alpha u)\|^2(1-\alpha)$ to its limit $\frac12 (\text{tr}\Sigma - u^T\Sigma u)$, $\Sigma$ denoting the covariance matrix of the vector $X$. This involves moments of order three and is elicited in Theorem~\ref{thm:higher-order}.   

As using such expansions rely on the existence of moments, we suggest an alternative approach without moment conditions and  
provide in Theorem~\ref{thm:upperbd-q} an upper bound for $\|q_X(\alpha u)\|$ as $\alpha\to 1$.

\begin{theorem}\label{thm:higher-order} 
Let $u \in S^{d-1}$ and $\{u,w_1,...,w_{d-1}\}$ be an orthonormal basis of $\real^d$.
\begin{enumerate}
	\item[(i)]  If $\E\|X\|^2<\infty$, we have
		\begin{equation}\label{eqn:2nd-order}
			\begin{aligned}
			&\|q_X(\alpha u)\| \Big[ q_X(\alpha u)-\big\{\|q_X(\alpha u)\|u + \E(X - \langle X,u\rangle u)\big\}\Big]\qquad\\ 
			&\; \underset{\alpha\to 1}{\longrightarrow}\, -\half\| \E(X-\langle X,u\rangle u)\|^2 u +\sum_{k=1}^{d-1} \text{cov}(\langle X,u\rangle , \langle X, w_k\rangle)\, w_k.
			\end{aligned}
		\end{equation}
	\item[(ii)] If $\E\|X\|^3 < \infty$, then, as $\alpha\to 1$,  
		\begin{eqnarray}\label{eqn:3rd-order}		
        &&\!\!\!\!\!\!\!\!\!\!\!\!\!\!\!\!\!\!\!\!\!  \|q_X(\alpha u)\|\Big[\|q_X(\alpha u)\|^2 (1 - \alpha ) -  \half (\text{tr } \Sigma - u^T\Sigma u)\Big] \, \longrightarrow \nonumber \\
        && \E\Big(\langle X,u\rangle \big[\|X-\langle X,u\rangle u\|^2 - \langle X,\E\left(X-\langle X,u\rangle u\right)\rangle \big]\Big)-\sum_{k} \text{cov}(\langle X,u\rangle, \langle X,w_k\rangle)\, \E(\langle X,w_k\rangle ) \nonumber\\
        && = \sum_{k=1}^{d-1}\Big[\text{cov}\big(\langle X,u\rangle\langle X,w_k\rangle, \langle X,w_k\rangle\big) -\text{cov}(\langle X,u\rangle, \langle X,w_k\rangle)\, \E(\langle X,w_k\rangle ) \Big]
		\end{eqnarray}
	If $\langle X,u\rangle$ and $\langle X,w_k\rangle$ are independent for every $k=1,\ldots,d-1$, then, the result can be simplified to
	\begin{eqnarray*}
&&\!\!\!\!\!\!\!\!\!\!\!\!\!\!\!\!\!\!\!\!\! \|q_X(\alpha u)\|\Big[\|q_X(\alpha u)\|^2 (1 - \alpha ) -  \half (\text{tr } \Sigma - u^T\Sigma u)\Big] \, \longrightarrow \nonumber \\
&& \E(\langle X,u\rangle) \left( tr \Sigma - u^T\Sigma u\right) -\sum_{k} \text{cov}(\langle X,u\rangle, \langle X,w_k\rangle)\, \E(\langle X,w_k\rangle )
		\end{eqnarray*}
\end{enumerate}
\end{theorem}
Note that this result is quite informative for a very heavy tail distribution, for which there is no second moment, or a moderate heavy one, for which the third moment exists. Notice that a gap remains in Theorem~\ref{thm:higher-order} when there exists a second moment but not a third one; this drawback is due to the chosen approach, using an expansion in our analysis. 
So, we investigate it empirically in the Appendix, considering various examples. In particular we vary the number $\nu$ of degrees of freedom of a multivariate Student distribution to cover the theoretical gap. We observe that the heavier the tail, the slower the convergence in Equation~\eqref{eqn:3rd-order}, 
especially in the absence of moments of order larger or equal to 3 (case $\nu=2.2$).


Moreover, using another approach, we can provide an upper bound for $\|q_X(\alpha u)\|$ with no moments assumption:
\begin{theorem}\label{thm:upperbd-q}
    Let $u \in S^{d-1}$, $0<\alpha < 1$ and $k_\alpha>0$.
    Then, we have
    $$
    \|q_X(\alpha u)\| \le \frac{2\,k_\alpha \, \pr(\|X\|\le k_{\alpha})}{1-2\pr(\|X\|\,>\, k_\alpha)-\alpha}.
    $$
\end{theorem}



\subsection{Halfspace depth}
\label{ss:depth}


\begin{definition}[Depth function]
A depth function corresponding to a distribution $F$ is a non-negative function $D(x,F)$ defined at every point $x \in \real^d$, which provides an outward ordering from the {\it center of the distribution}.
\end{definition}
Since it is desirable that depth functions decrease to zero in every direction from the median/center, the {\it center} of a distribution is defined as the point of maxima of the depth function. In case the maximum is attained at multiple (finitely many) points then the centroid of all such points is called the median (and the center). 


We shall focus on the following notion of halfspace depth. It is widely used, and is also a good representative of the class of depth functions, as it satisfies most of the desirable properties for depth functions (see \cite{Mosler2021}, Table 2). 
\begin{definition}[Halfspace depth; \textnormal{\cite{Tukey1975}}] 
    \label{def:HD}
    For a probability distribution $\pr$ defined on $\real^d$, the halfspace depth is given by: 
    $HD(x,\pr)= \inf \{\pr(H):\,H \in \calH_x\},$ where $\calH_x$ denotes the set of halfspaces in $\real^d$ containing $x\in\real^d$. Specifically, if $\pr$ has a probability density function $f$, then
$\displaystyle
HD(x,\pr)=\inf_{|p|=1}\int_{\{y: \langle y-x,p\rangle \geq 0\}} f(y) \,  dy.
$
\end{definition}
Intuitively, a high depth point indicates that it is more central, while a low depth point denotes a relatively extremal point.

When no possible confusion, we will denote $HD(x,\pr)$ by $HD(x)$.
\vspace{-2ex}
\begin{figure}[H]
    \centering
    \includegraphics[width=.6\linewidth]{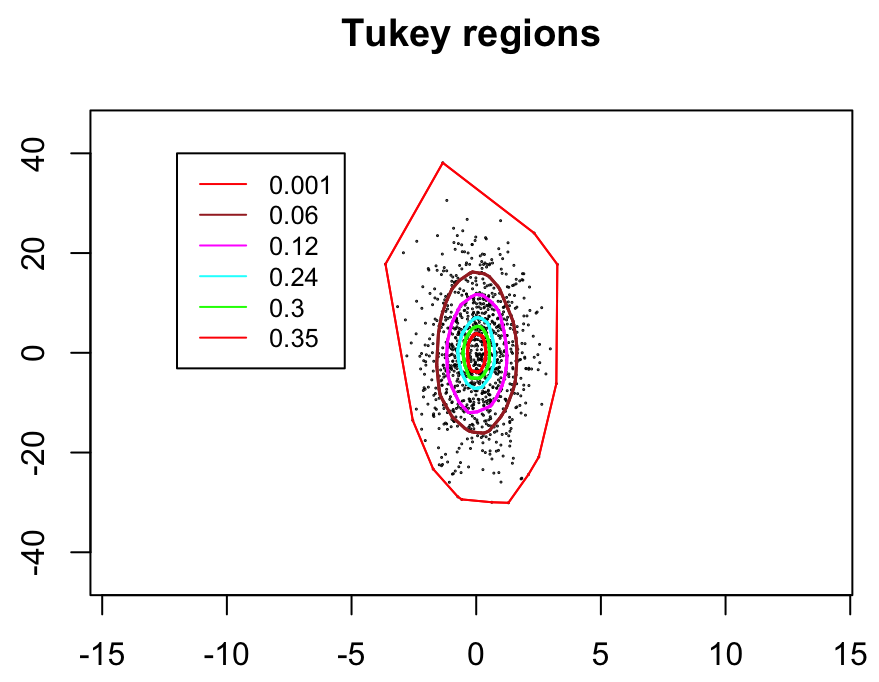}
    \vspace{-2ex}
    \parbox{340pt}{\caption{\sf Representation of the Tukey contours for 6 different depths, considering a sample of $1000$ observations (black points) from a mean zero Gaussian distribution with covariance $\text{diag}(1,100)$.}
    \label{fig:Tukey-regions}}
\end{figure}
This is illustrated in Figure~\ref{fig:Tukey-regions}, when considering a Gaussian sample, using the R-Package from \cite{Barber2022}, based on \cite{Liu2019}. Note that the symmetry present in the underlying sampling distribution may not be seen precisely in the simulated sample. Isodepth contours are drawn for different values of depth. Observe that the asymmetry is more evident in the extremes, than in the bulk region.


\paragraph{Decay rate of halfspace depth -}

Recall that our motivation to study the multivariate quantiles is to understand the connection between the extremal behaviour of a probability measure and properties of multivariate quantile functions. It is, therefore, natural to analyse and understand the asymptotics of depth functions along the same lines as in the case of geometric quantiles. First, we consider the important example of asymptotically elliptically symmetric distributions. Then, we question the halfspace asymptotic decay depending on the light or heaviness of the underlying probability measure. While the heavy tail case has been studied in \cite{He2017}, we complete the picture by providing the corresponding result for the light tail case.

\begin{definition}[Asymptotically elliptically symmetric]
A probability measure defined on $\real^d$ with density $f$ is said to be asymptotically elliptically symmetric with respect to some symmetric, non--singular matrix $\Sigma$, if $$\lim_{t\to\infty} \sup_{\|\Sigma x\|=\|\Sigma y\|} \frac{f(tx)}{f(ty)} = 1.$$
\end{definition}

\begin{theorem}\label{thm:asymp-depth-ell-symm}
If the probability density function $f$ of a distribution $\pr$ is asymptotically elliptically symmetric with respect to $\Sigma$, then the halfspace depth corresponding to $f$ satisfies 
$$\displaystyle \lim_{t\to\infty} \sup_{\|\Sigma x\|=\|\Sigma y\|}\frac{HD(tx,\pr)}{HD(ty,\pr)} =1.$$
\end{theorem}

%

Next, we shall state a generic result comparing distributional asymptotics with that of the induced halfspace depth, when marginals of the joint distribution may have different asymptotic behaviours. For that, we first recall another equivalent way to define the halfspace depth in terms of one dimensional distributions (see e.g. \cite{Donoho1992}), as follows:
\begin{equation}\label{eqn:hdepth-project}
    HD(x,\pr) = \min_{h:|h|=1}(1- F_h(\langle h,x\rangle))
\end{equation}
where $F_h$ is the c.d.f. of the univariate projection of $\pr$ onto the $h$ direction.
As a result, we have the following simple, but interesting, observation:
\begin{proposition}
    \label{thm:hd-decay-survivalMargins}
    Let $\pr$ be a probability measure on $(\real^d,\calB(\real^d))$ with cumulative distribution function $F$ and marginals $F_i$, $1\le i \le d$.
    Then, we have
    $$
    HD(t\,x,\pr)\leq \min_{1\le i\le d} \left(1-F_i(t \, x_i)\right).
    $$
    where $x_i$ denotes the $i$-th coordinate of $x$. 
\end{proposition}

As a consequence, if the $i$-th marginal distribution has a light tail, then the decay of halfspace depth cannot be slower than exponential. If all marginals are heavy-tailed, then the bound will correspond to the one with the smallest tail index. This gives us a simple tool to discriminate between light and heavy tails. Note that one could choose different sets of marginals, and the upper bound would still hold true. Indeed, in practice, this method is useful if one has some prior information of the distribution along certain marginals.

The result given in Proposition~\ref{thm:hd-decay-survivalMargins} also applies to projection depth functions, which are of great use because of their properties (see e.g. \cite{Mosler2021} for definitions and illustrations).

In fact, we can go further and obtain the precise decay of the halfspace depth in the light tail case, namely: 

\begin{theorem} 
\label{thm:hd-decay-light-tail}
Let $\pr$ be a probability measure on $(\real^d,\calB(\real^d))$.
\begin{itemize}
    \item[(i)] If its moment generating function $M_{\pr}(h) = \int_{\real^d}e^{\langle h,y\rangle}\,\pr(dy)$ exists for some $h\neq 0$, then, the halfspace depth is also light tailed:
\begin{equation}\label{eqn:light-tail}
    HD(tx,\pr) = {\textrm O}(e^{-t})\quad \text{as}\quad t\to \infty.
\end{equation}
In particular, if any marginal distribution of $\pr$ has exponential moments, then \eqref{eqn:light-tail} holds true.
    \item[(ii)] Assume there exists a positive function $f$ such that 
    \begin{equation}\label{eqn:unif-direction-density-decay}
        \liminf_{R \to \infty} \inf_{\theta\in S^{d-1}}\dfrac{\pr\big( \langle\theta, X\rangle \geq R\big)}{f(R)} >0.
    \end{equation}
Then, there exists $t_0>0$ such that, 
    \begin{equation*}
        HD(t x, \pr) \geq c\,\, f(t\|x\|), \,\,\,\,\, \forall \,\,t > t_0,
    \end{equation*}
    for some constant $c>0$ that does not depend on $x$.   
\end{itemize}
\end{theorem}

\begin{Remark}\label{rk:boundsP-lightdecay}~
   \begin{enumerate}
   \item Assuming $M_{\pr}(h)$, defined in Theorem~\ref{thm:hd-decay-light-tail}, is finite for some $h\neq 0$, then combining Proposition~\ref{thm:hd-decay-survivalMargins} with \eqref{eqn:light-tail}, gives
   $$
   HD(tx,\pr) = {\textrm O}\left(e^{-t}\wedge \min_{1\le i\le d} \left(1-F_i(t \, x_i)\right)\right)\quad \text{as}\quad t\to \infty,
   $$
   where $F_i$ are the marginal distributions of $\pr$.
   \item Assuming \eqref{eqn:unif-direction-density-decay}, in addition to the finiteness of $M_{\pr}(h)$ for some $h\neq 0$, and combining the upper and lower bounds obtained in the above theorem, we have for sufficiently large $t$
    $$
        c\,f(t\|x\|) \le HD(t x,\pr) \le K\,\left(e^{-t}\wedge \min_{1\le i\le d} \left(1-F_i(t \, x_i)\right)\right),
    $$
    for any $x\ne 0$ with $K,\,c>0$.
    \item An immediate consequence of the lower bound is that, whenever $\pr$ satisfies condition \eqref{eqn:unif-direction-density-decay}, we have
    $$
    \pr(H) \ge c\, f(t\|x\|) , \qquad\forall \,\,H\in\calH_{tx}.
    $$
   \end{enumerate}
\end{Remark}

Theorem~\ref{thm:hd-decay-light-tail} completes the study on the relationship between the asymtotic behavior of the halfspace depth and its underlying probability measure,  since, in the heavy-tailed case, the decay rate has been provided by \cite[Proposition~2]{He2017}   under some classical conditions when considering a regularly varying framework. We recall that a positive measurable function $G$ is said to be {\it regularly varying (RV)} at infinity with index $-\beta$, with $\beta>0$, denoted by $G\in RV_{-\beta}$, if for any $x>0$,  $\displaystyle \frac{G(tx)}{G(t)} \underset{t\to\infty}{\longrightarrow} x^{-\beta}$ (if $\beta=0$, $G$ is said to be slowly varying at infinity).

Let us enunciate Proposition~2 in \cite{He2017} when assuming that the probability measure has a density that is RV.
\begin{proposition}[\textnormal{Proposition~2 in \cite{He2017}, adapted to our notation and case}]
\label{thm:hd-decay-heavy-tail}
    Let $\pr$ be a probability measure on $(\real^d,\calB(\real^d))$ with non--vanishing continuous density $f$ on $\real^d$, such that the map $y\mapsto \|y\|^d f(y)$ is bounded in every compact neighbourhood of the origin, and there
    exist a positive function $\la:\real^d \to \real^{+}$ and a function $V\in RV_{-\alpha}$, with $\alpha>0$, such that
    \begin{equation*}
        \left|\frac{f(ty)}{t^{-d}V(t)} -\la(y)\right| \underset{t\to\infty}{\longrightarrow} 0,\,\,\,\,\,\forall y\neq 0,
        \quad \text{and} \quad
        \sup_{\|y\|= 1}\left| \frac{f(ty)}{t^{-d}V(t)} -\la(y)\right| \underset{t\to\infty}{\longrightarrow} 0.
    \end{equation*}
    Then,  we have
    \begin{equation}\label{eqn:half-sp-decay-reg-varying}
        \lim_{t \to \infty} \frac{HD(tx,\pr)}{V(t)} = \inf_{\substack{|p|=1\\p^Tx>0}} \int_{\{z:p^T(z-x)\geq 0\}} \la(z) dz.
    \end{equation}      
\end{proposition}

In view of the results obtained, we conclude that the asymptotic behaviour of halfspace depth does reflect well the asymptotic behaviour of the underlying probability measure. We shall also question if depth functions mirror the tail behaviour of the underlying distribution, when going from population to sample versions, which shall form the crux of our investigation in Section~\ref{sec:empirical-quantiles}.


\section{Sample versions and applications}
\label{sec:empirical-quantiles}

While our results concerning asymptotic behaviour of multivariate quantiles are of theoretical interest, they also set the stage for the sample versions of multivariate quantiles for which results on rate of convergence are crucial for applications.

There are several statistical problems wherein our results can readily be used. For instance, comparing statistical distributions through their respective samples is an old problem. Clearly, this problem is, mostly, ill--defined given a finite sample. Indeed, several simplifying assumptions are needed to make the problem tractable. 


\subsection{Empirical multivariate geometric quantiles}\label{sec:empirical-geometric}

Early theoretical results on the sample asymptotics of empirical geometric quantiles have been obtained in \cite{Chaudhury1996}, and more recently in \cite{Girard2015,Paindaveine2021}. We further pursue this direction, allowing for simultaneous increase of sample size and the magnitude of the index vector. Moreover, our results can serve as tools for distinguishing various distributions based on their tail behaviour. 

Let us begin with recalling the definition of the empirical or sample version of the geometric quantile.
\begin{definition}[Sample geometric quantile; \textnormal{\cite{Chaudhury1996}}]
Let $\{X_n\}_{n\ge 1}$ be i.i.d sample drawn from $\pr$ defined on $\real^d$. Assume $u \in \real ^d$ be such that $\|u\|<1$. Then the sample geometric quantile, denoted by $\hat{q}_n( u)$, is defined as (dropping the dependence on the distribution for notational simplicity)
\begin{equation}\label{eqn:sample geometric-quantile-dfn}
\displaystyle \hat{q}_n( u)=	
\argmin_{q \in \real^d} \left\{\frac1n\sum_{i=1}^n\big(\|X_i-q\|-\|X_i\|\big)-\langle u,q \rangle\right\}.
\end{equation}
\end{definition}

As can be observed immediately from the definition of sample geometric quantiles, the underlying empirical measure does not satisfy the continuity properties needed to prove the results stated in the previous section. Therefore, we begin with elementary asymptotic estimates in Theorem~\ref{thm:empirical-geom-quant-asymp1}, which can be seen as improvements of Theorems $2$ \&  $3$ of \cite{Paindaveine2021}. 

The broad set of assumptions needed for the underlying measure $\pr$ on $\real^d$ are:\\[-5ex]
\begin{itemize}
\item[\it (A1)] $\pr$ does not have any atom, and its support is not contained in any unit dimensional affine subspace of $\real^d$.
\item[\it (A2)] The density function of $\pr$ is bounded on every compact subset of $\real^d$.
\end{itemize}
\begin{theorem}\label{thm:empirical-geom-quant-asymp1}
Let $\{X_n\}_{n\ge 1}$ be an i.i.d. sample drawn from $\pr$ on $\real^d$ satisfying assumptions (A1) and (A2). Let $\{\alpha_n\}_{n\ge 1}\subset (0,1)$ be such that  $\displaystyle \lim_{n \to \infty} \alpha _n =1$. Then, for every $u\in S^{d-1}$:
\begin{itemize}
    \item[(i)] $\displaystyle\lim_{n \to \infty} \|\hat q_{n}({\alpha_n}u)\|=+ \infty$, a.s.
    \item[(ii)] $\displaystyle\lim_{n \to \infty} \frac{\hat q_{n}({\alpha_n}u)}{\|\hat q_{n}({\alpha_n}u)\|}=u$, a.s.
\end{itemize}
\end{theorem}

Subsequently, we graduate to finer results revealing the asymptotic character of sample geometric quantiles in Theorems \ref{thm:empirical-geom-quant-asymp2} and \ref{thm:numerics}. 
\begin{theorem}\label{thm:empirical-geom-quant-asymp2}
Under the same assumptions as in Theorem~\ref{thm:empirical-geom-quant-asymp1}, suppose $\{\alpha_n\}$ also satisfies
\begin{equation}\label{eq:cdtionAlpha}
    \sum_{n=1}^{\infty}\exp(-n(1-\alpha_n)^2)<\infty. 
\end{equation}
For $u,v\in S^{d-1}$, we have:
\begin{itemize}
\item[(i)] If $\E\|X_1\| <\infty$, then, 
\begin{equation}\label{eqn:empirical-geom-quant-asymp2}
    \| \hat q_n(\alpha_n u)\|\left\langle\alpha_n u- \frac{\hat q_n(\alpha_n u)}{\|\hat q_n(\alpha_n u)\|}, v\right\rangle \overset{a.s.}{\underset{n\to\infty}{\longrightarrow}} -\E\big\langle X-\langle X, u \rangle u,v \big\rangle.
\end{equation}
\item[(ii)] If $\E\|X_1\|^2<\infty$ and $\displaystyle n (1-\alpha_n)^3 \to \infty$ as $n\to\infty$, then,
\begin{equation}\label{eqn:empirical-geom-quant-asymp3}
    \| \hat q_n(\alpha_n u)\|^2\left\langle\alpha_n u-\frac{\hat q_n(\alpha_n u)}{\| \hat q_n(\alpha_n u)\|},\frac{\hat q_n(\alpha_n u)}{\| \hat q_n(\alpha_n u)\|}\right\rangle \overset{a.s.}{\underset{n\to\infty}{\longrightarrow}} -\half \E \|X-\langle X,u\rangle u\|^2.
\end{equation}
\end{itemize}
\end{theorem}

\begin{theorem}\label{thm:numerics}
Under the assumptions of Theorem~\ref{thm:empirical-geom-quant-asymp2},
we have, for $u \in S^{d-1}$:
\begin{enumerate} 
\item[(i)] If $\E\|X_1\| <\infty$, then,
$$\hat{q}_n(\alpha_n u)-\left\{ \|\hat{q}_n(\alpha_n u)\|u+\frac{1}{n}\sum_{i=1}^{n}\left(X_i-\langle X_i,u \rangle u\right)\right\} \overset{a.s.}{\underset{n\to\infty}{\longrightarrow}} 0.$$ 
\item[(ii)] If $\E\|X_1\|^2<\infty$ and $\displaystyle n (1-\alpha_n)^3 \to \infty$ as $n\to\infty$, then 
\begin{equation}
\|\hat{q}_n(\alpha_n u)\|^2(1-\alpha_n) \overset{a.s.}{\underset{n\to\infty}{\longrightarrow}} \half (\text{tr } \Sigma -\langle\Sigma u,u\rangle) \,\, > 0.
\end{equation}
where $\Sigma$ denotes the covariance matrix corresponding to $\pr$.
\end{enumerate}
\end{theorem}
Theorems \ref{thm:empirical-geom-quant-asymp2} and \ref{thm:numerics} can be seen as sample versions of results in \cite{Girard2017}.
We could also provide further decay rates assuming higher moments, as done for the population version. 
%
%


\subsection{Empirical multivariate halfspace depth}
\label{sec:empirical-HD}

As noted earlier in the introduction, our motivation for studying geometric quantiles and halfspace depths is that they are known to capture certain behaviour of the underlying distribution. Contrary to geometric quantiles, the halfspace depths do not uniquely characterise the underlying distribution (see \cite{Nagy2021a}). However, in \cite{Nagy2022}, the author listed eight situations for which the halfspace depths uniquely identify the underlying distribution, empirical measure being one among them, as proven in \cite{Struyf1999}.

While analysing the rate of decay of halfspace depth of empirical measures, it is appealing to compare the decay rate of halfspace depths of empirical measures with those of the measures from which the samples have been generated. This rate of decay is evaluated in the following theorem.
\begin{theorem}
\label{thm:empHDdecay-vs-popHD}
Let $\pr$ be a probability measure defined on $(\real^d,\calB(\real^d))$, and let $g$ be the capacity function corresponding to $\pr$, as defined in \cite{Alexander1987}[p.382]. 
Consider a sequence $\{\gamma_n\}_{n\ge 1}$ satisfying the following conditions
$$
\centering
        (C1) \qquad n^{-1}\log (g(\gamma_n)) = {\textrm o}(\gamma_n)\quad \text{and} \quad n^{-1}\log \log n = {\textrm o}(\gamma_n),
$$  
and a sequence $\{t_n\}_{n\ge 1}$ with $\displaystyle t_n \underset{n\to\infty}{\to}\infty$, such that  
$$
(C2) \quad HD(t_n x,\pr)>\gamma_n \quad \text{for large enough}\;n. 
$$ 
Let $\{X_k\}_{k\ge 1}$ be an i.i.d. sample drawn from $\pr$, and $HD(\cdot,\pr_n)$ be the halfspace depth for the empirical measure $\pr_n = \frac1n \sum_{i=1}^n \delta_{X_i}$. Then, for any $x\neq 0$, we have
\begin{equation}\label{eqn:ratio-emp-pop-HD}
    \left|\dfrac{HD(t_n x, \pr_n)}{HD(t_n x, \pr)}-1\right| \underset{n\to\infty}{\to} 0\quad a.s.
\end{equation}
\end{theorem}

The expression of interest \eqref{eqn:ratio-emp-pop-HD} in Theorem~\ref{thm:empHDdecay-vs-popHD} is difficult to work with, hence we shall first simplify this expression in a straightforward way, with the bound given in Lemma~\ref{prop:hdepth-sup-emp-ratio-gl}. It will need to be adapted to our context (see Proposition~\ref{prop:hdepth-sup-emp-ratio}) when developing the proof.

\begin{Remark}
Condition $(C2)$ implies that in order to apply the above result in any setting, we must have a reasonable way of estimating $HD(t_n x,\pr)$. Indeed this is the case, as we have provided estimates for halfspace depths in Section \ref{ss:depth}, which we are going to use for this purpose. Moreover, as will be seen in the proof, the condition of the halfspace depth is transferred to an appropriate decay condition on the tail probabilities of the measure $\pr$.
\end{Remark}

\begin{lemma}\label{prop:hdepth-sup-emp-ratio-gl}
For any $x\neq 0$, we have the following inequality:
       $$ 
        \left|\dfrac{HD(t_n x, \pr_n)}{HD(t_n x, \pr)}-1\right| \le \sup_{H \in \calH_{t_nx}} \left|\frac{\pr_n(H)}{\pr(H)} - 1 \right|.
       $$  
\end{lemma}

The bound thus obtained has been analysed by many researchers in one or the other form. Specifically, the rate of convergence of $\displaystyle \sup_{H \in \calH_{tx}} \left|\frac{\pr_n(H)}{\pr(H)} - 1 \right|$ has been an object of interest since long with the Glivenko--Cantelli theorem being one of the earliest in this direction. Later, the rate of convergence was obtained by several authors, e.g. \cite{Alexander1987,Gine2006,Shorack2009,Wellner1992,Burr2017} in different scenarios with specific assumptions. It is noteworthy that most of the results in this direction use the specific structure of $\calH_x$ and Dvoretzky--Kiefer--Wolfowitz (DKW) inequality (\cite{Dvoretzky1953}). The specific structure we refer to is called the Vapnik--Chervonenkis (VC) class (see \cite{Vapnik1971,Dudley1984,Alexander1987,Talagrand2003}.The idea of VC class has its roots in statistical learning wherein one is interested in identifying the class of functions to characterise convergence of probability measures. Specifically, a class $\calS$ of sets shatters a finite set $F$ if, given $G\subset F$, $\exists S\in\calS$ for which $G = F \cap S$. A class $\calS$ of sets is called a VC class if for some integer $n$, $\calS$ does not shatter any set of cardinality $n$.
In our analysis, we shall use the approach of \cite{Alexander1987} on VC class (without resorting to the DKW inequality). Specifically, we shall invoke
Theorem~5.1 of the \cite{Alexander1987},
as it is a powerful and crucial result for the proof of Theorem~\ref{thm:empHDdecay-vs-popHD}. 

Notice that, without any assumption on $\pr$, Theorem~\ref{thm:empHDdecay-vs-popHD} presents a general result about the rate of decay of the empirical halfspace depth and allows us to compare the halfspace depth of the parent measure $\pr$ and of the empirical measure $\pr_n$. Indeed, \eqref{eqn:ratio-emp-pop-HD} can be expressed as
$$
1-\varepsilon \le \frac{HD(t_n x, \pr_n)}{HD(t_n x, \pr)} \le 1+\varepsilon, \quad \forall n>n_0,
$$
from which we can deduce that there exist random $c_1, c_2$ with $0< c_1 \le c_2$ such that, for $n$ large enough,
\begin{equation}\label{eqn:emp-pop-hdepth}
    c_1\, HD(t_nx, \pr)\,\le \, HD(t_n x, \pr_n)  \,\le  \,c_2 HD(t_nx, \pr).
\end{equation}
However, as seen in Section \ref{ss:depth}, the rate of decay of halfspace depth is closely related to the tail behaviour of $\pr$, which, in view of \eqref{eqn:emp-pop-hdepth}, implies that the rate of decay of the empirical halfspace depth can be estimated as a function of the tail behaviour of the parent measure $\pr$. The following two theorems establish this connection.

When assuming a multivariate regularly varying framework as in Proposition~\ref{thm:hd-decay-heavy-tail} (for geometric quantiles), we obtain the following rate of convergence for the empirical halfspace depth:
\begin{theorem}\label{thm:emp-hdepth-mrv-decay}
    Let $\pr$ be a probability measure on $(\real^d,\calB(\real^d))$ with non--vanishing continuous density $f$ on $\real^d$, such that the map $y\mapsto \|y\|^d f(y)$ is bounded in every compact neighbourhood of the origin, and there exist a positive function $\la:\real^d \to \real^{+}$ and a function $V\in RV_{-\alpha}$, with $\alpha>0$, such that
    \begin{equation*}
        \left|\frac{f(ty)}{t^{-d}V(t)} -\la(y)\right| \underset{t\to\infty}{\longrightarrow} 0,\,\,\,\,\,\forall y\neq 0,
        \quad \text{and} \quad
        \sup_{\|y\|= 1}\left| \frac{f(ty)}{t^{-d}V(t)} -\la(y)\right| \underset{t\to\infty}{\longrightarrow} 0.
    \end{equation*}
Let $\{X_n\}_{n\ge 1}$ be an i.i.d sample drawn from $\pr$, and $g$ be the capacity function corresponding to $\pr$. Then we have, for any $x\neq 0$,
    \begin{equation}
        \lim_{n \to \infty}\frac{HD(t_n x,\pr_n)}{V(t_n)} = \inf_{\substack{|p|=1\\p^Tx>0}} \int_{\{z:\langle z-x,p\rangle \geq 0\}} \la(z) dz, \,\,\,\,\,\,\,\text{a.s.}
    \end{equation}
whenever $(t_n)\nearrow \infty$ with $n$ is such that $V(t_n)>\, \gamma_n$ (for any large $n$), and $(\gamma_n)$ satisfies Condition $(C1)$ given in Theorem \ref{thm:empHDdecay-vs-popHD}.  
\end{theorem}
\begin{Example} For $\alpha>0$, consider a multivariate regularly varying distribution with index $-\alpha$, then $V(t)=t^{-\alpha}$. Therefore Theorem~\ref{thm:emp-hdepth-mrv-decay} holds if $t_n\nearrow \infty$ and $t_n \leq \gamma_n^{1/\alpha}$. Now, by choosing $\gamma_n=n^{-\beta}$ with $0<\beta<1$, the condition $V(t_n)>\, \gamma_n$ gives a speed of $t_n \leq n^{\beta/\alpha}$.
\end{Example}

Let us turn to the light tail case, for which we can provide, under distinct conditions, a lower and an upper bound for the asymptotics of the halfspace depth. Through examples, we observe that the general bounds can be tight as in the exponential case (see Example~\ref{exple:exp-normal}), with the lower bound of the order of the upper one, showing that this order is optimal in the general case. In the Gaussian case, the gap between the bounds is much larger and can be improved via a direct computation as given in Example~\ref{exple:exp-normal}.
\begin{theorem}\label{thm:emp-hdepth-light-decay}
Let $\{X_n\}_{n\ge 1}$ be an i.i.d sample drawn from a probability measure $\pr$ defined on $(\real^d,\calB(\real^d))$, $g$ be the capacity function corresponding to $\pr$, and $x$ be any unit vector.
\begin{itemize}
    \item[(i)] Assume that the moment generating function $M_{\pr}(h) = \int_{\real^d}e^{\langle h,y\rangle}\,\pr(dy)$ of $\pr$ exists for some $h\neq 0$. Choosing $(\gamma_n)$ satisfying Condition $(C1)$ in Theorem~\ref{thm:empHDdecay-vs-popHD}, we have,
    \begin{center}$\displaystyle
    HD(t_n x,\pr_n) = {\textrm O}_{\pr}(e^{-t_n}) \quad \text{when}\quad t_n\underset{n\to\infty}{\longrightarrow}\infty \quad\text{and}\quad t_n < \log(\gamma_n^{-1}).
    $ \end{center}
    \item[(ii)] Assume there exists a positive function $f$ satisfying 
    \eqref{eqn:unif-direction-density-decay}, then, we have
      \begin{center}$\displaystyle
            HD(t x, \pr_n) 
            \geq c\,\, f(t\|x\|), \,\,\,\,\, \text{ for large}\;t,
       $ \end{center}
        for some random variable $c>0$ that does not depend on $x$ and $n$.  
    \item[(iii)] Combining the conditions of (i) and (ii), and choosing $t_n$ such that $f(t_n)> \gamma_n$, we can write 
    \begin{center}$\displaystyle c\,\gamma_n \le HD(t_n x,\pr_n) \le K\,e^{-t_n}, \quad \text{with}\; K\; \text{a positive random number}\;\text{and}\; c>0, \text{ as in } (ii).
    $ \end{center}
\end{itemize}
\end{theorem}

\begin{Remark}~
    Using Proposition~\ref{thm:hd-decay-survivalMargins} for the empirical measure $\pr_n$ and the relation \eqref{eqn:emp-pop-hdepth} between halfspace depth of $\pr$ and of $\pr_n$, the upper bound given in (i) can be refined as
        \begin{equation}\label{eq:upperBd-improved}
            HD(t_n x,\pr_n) \le K\,\left(e^{-t_n}\wedge \min_{1\le i\le d} \left(1-F_i(t_n \, x_i)\right)\right)\quad \text{as}\quad t_n\to \infty \;\text{with}\;n,
        \end{equation} 
    where $F_i$ denotes the $i$-th marginal distribution of $\pr$, and $x_i$ denotes the $i$-th component of $x$, for $i=1,\cdots,d$.
\end{Remark}

\begin{Example}\label{exple:exp-normal}~     
    \begin{enumerate}
        \item {\bf Exponential case.} 
        Let $Y$ be a random vector with probability density function $h$ given by
    $\displaystyle h(y)=ke^{-\|y\|}$ for $y \in \real^d$.
    Then, we have     
\begin{equation}\label{eq:ExponExple}
    \lim_{R \to \infty} \inf_{\theta\in S^{d-1}}\dfrac{\pr\big( \langle\theta, Y\rangle \geq R\big)}{e^{-R}} >0.
\end{equation}
\textnormal{In order to prove \eqref{eq:ExponExple}, observe that, for any $\theta_1, \theta_2\in S^{d-1}$, 
$\displaystyle
\pr\big( \langle\theta_1, Y\rangle \geq R\big)=\pr\big( \langle\theta_2, Y\rangle \geq R\big)$, 
due to the spherical symmetry of the distribution of $Y$.
Therefore, we can consider any direction and, choosing $\theta=(1,0,...,0)$ and using Minkowski inequality, we can write,
\begin{align*}
        \inf_{\theta\in S^{d-1}}\pr\big( \langle\theta, Y\rangle \geq R\big)
        &=\pr\big( Y_1 \geq R\big)
        =\int_{-\infty}^{\infty}\cdots \int_{-\infty}^{\infty}\int_R^{\infty} ke^{-\|y\|}dy\\
        & \geq \int_{-\infty}^{\infty}\cdot \cdot \cdot \int_R^{\infty} ke^{-(|y_1|+...+|y_d|)}dy =K \int_R^{\infty} e^{- y_1} dy_1
        =Ke^{-R},
\end{align*}
for some constant $K>0$.
We deduce that \;$\displaystyle \frac{\inf_{\theta\in S^{d-1}}\pr\big( \langle\theta, Y\rangle \geq R\big)}{e^{-R}}>0$ for large enough $R$, hence \eqref{eq:ExponExple}.} \hfill \qed
        \item {\bf Gaussian case.} Consider a multivariate standard normal distribution. We choose $f(R)=e^{-\frac{R^2}{2}}$ (based on the Mill's ratio). Now, setting $\gamma_n=n^{-\beta}$ with $0<\beta <1$, the condition $f(t_n)> \gamma_n$ gives $t_n \leq \sqrt{2 \beta \log n}$, and, for such $t_n$, we have
        $\displaystyle 
        HD(t_n x,\pr_n)> c\,n^{-\beta}$. 
        The upper bound as given in (i), of order $e^{-t_n}$, is then
        quite large compared with the lower bound. Nevertheless, it can be improved, simply by considering the Gaussian marginal distributions as in \eqref{eq:upperBd-improved}. In such a case, the upper bound for $\displaystyle t_n = \sqrt{2 \beta \log n}$ becomes
        $$ 
        HD(t_n x,\pr_n)\le C\,n^{-\beta}, \text{ for some random variable $C>0$, independent of $x$ and $n$.} 
        $$ 
    \end{enumerate}
\end{Example}
    

\subsection{Illustration and discussion} 

In view of the various theoretical results obtained on the sample multivariate geometric quantiles and halfspace depths, in terms of rates of growth or decay, respectively, or of tail behavior, we turn to empirical illustrations of those results, considering light and heavy tailed-distributions.
Geometric quantiles are evaluated using the algorithm given in \cite{Chaudhury1996}, while halfspace depth contours are evaluated using the R--package developed in \cite{Pokotylo2019}.

\subsubsection{Geometric quantiles}
\label{sss:plots-geom}

\begin{figure}[H]
\centering
\begin{minipage}{0.49\textwidth}
    \includegraphics[width=.7\linewidth]{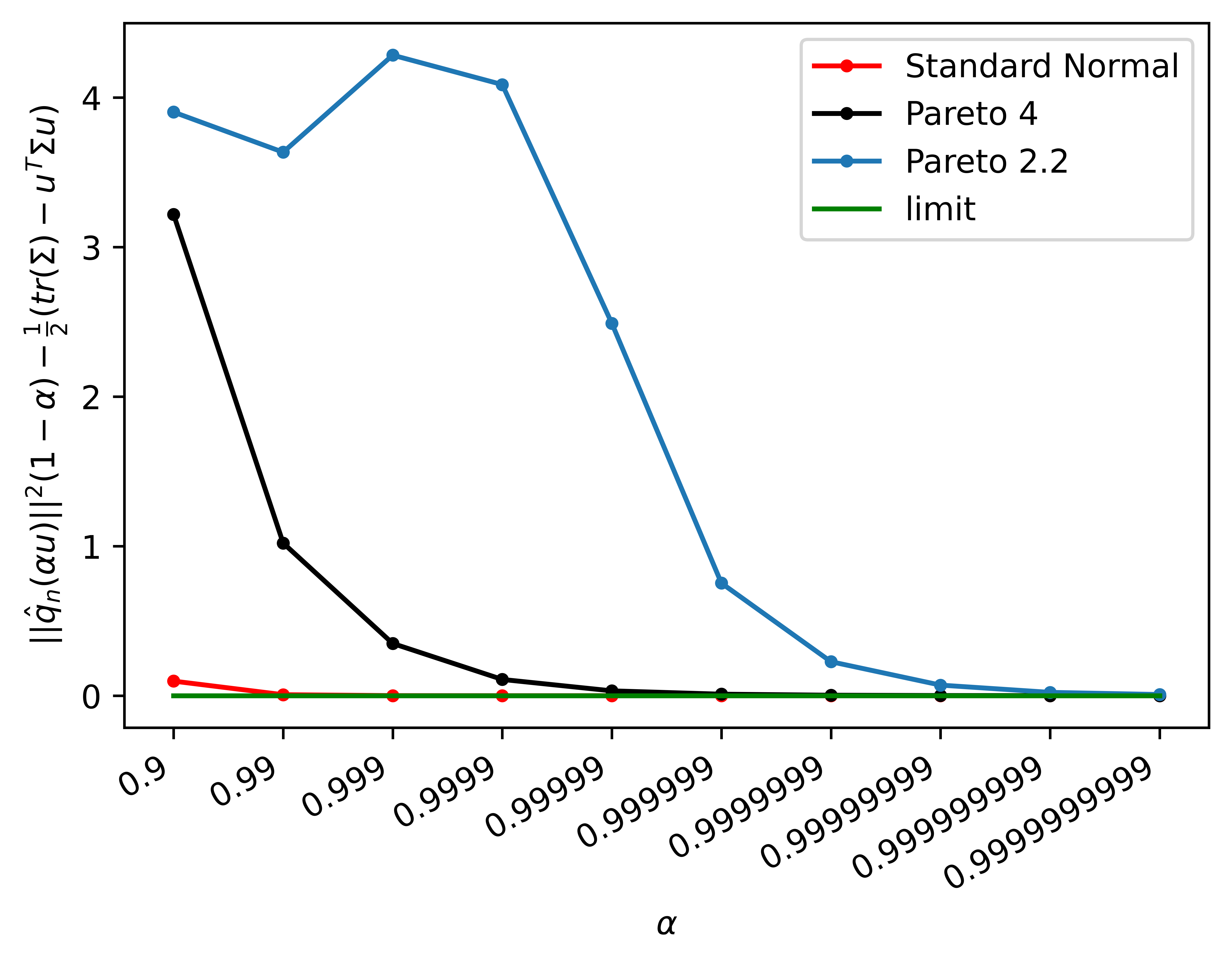}
\end{minipage}
\hfill
\begin{minipage}{0.49\textwidth}
    \includegraphics[width=.7\linewidth]{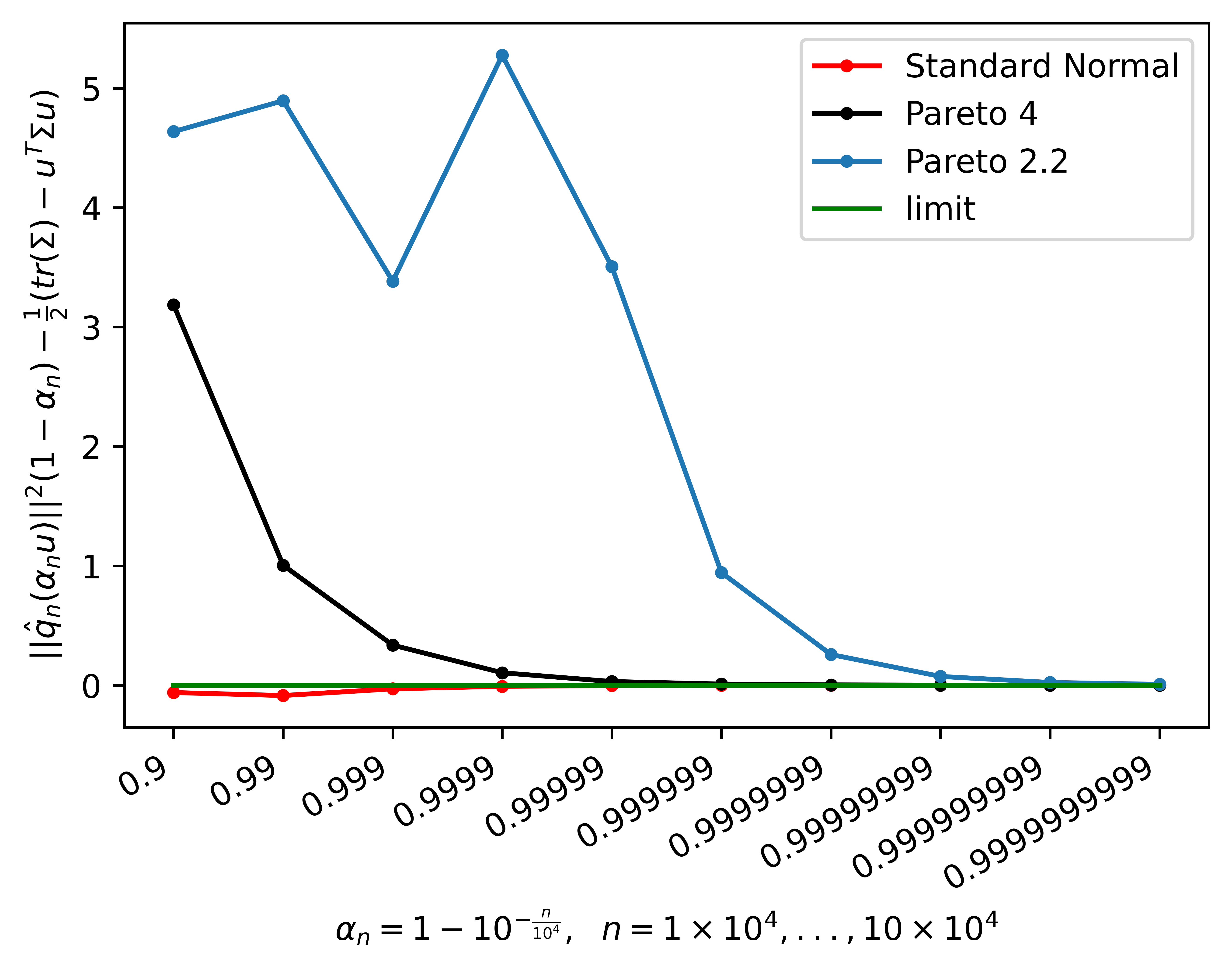}
\end{minipage}    
\parbox{460pt}{\caption{\sf \small Comparing the convergence for different types of tail distributions, using the 2nd order characterisation of geometric quantiles. The $x$-scale is chosen for the $\alpha$-values to look equispaced. Left plot: Number of observations $n=10^5$ is fixed; quantiles are computed corresponding to the index $\alpha\in\{\alpha(k)=1-10^{-k},\,k=1,\cdots,10\}$ and the unit vector $u$ in the direction $(1,1)$. Right plot: Growing sample at each point ({\it i.e.} for each $\alpha_n$); 10\% regularly growing sample, from $10^4$ to $10^5$ simulated values ($\alpha_n=1-10^{-n/10^4}$ for $n=10^4,\cdots,10^5$).}
\label{fig:GeomQuantilesOrder2-gl}}
\end{figure}

\begin{figure}[H] 
    \centering
        \includegraphics[width=.4\linewidth]{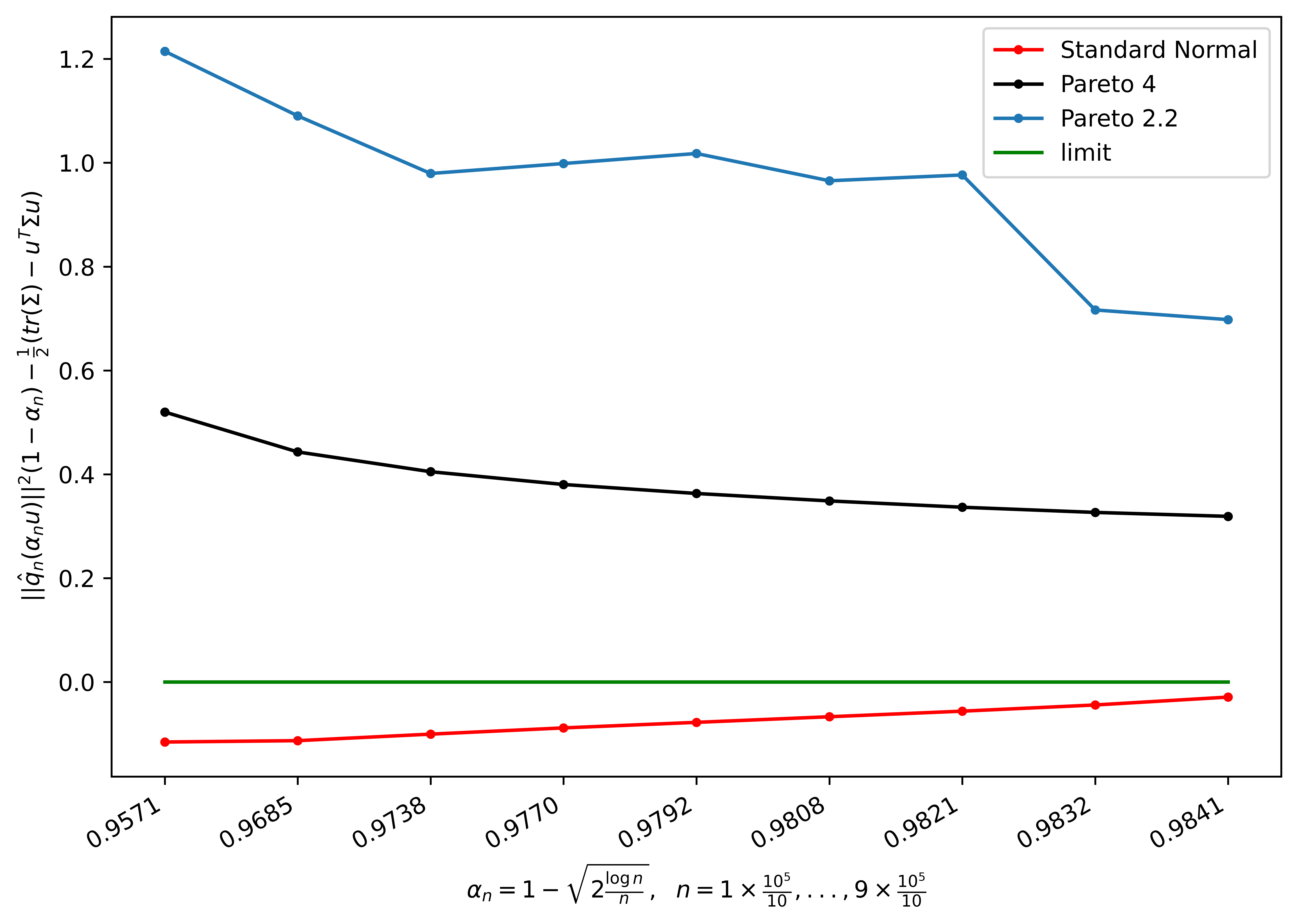}  
    \parbox{440pt}{\caption{\sf \small Slowly growing sample at each point ({\it i.e.} for each $\alpha_n$), up to $10^5$ simulated values, choosing $\alpha_n$ satisfying Condition~\eqref{eq:cdtionAlpha} (e.g. $\alpha_n=1-\sqrt{2\log n /n}$), which gives values of $\alpha_n$ going roughly from $0.90$ to $0.99$ (with a partition of 10, i.e. 10 values of $\alpha_n$). Direction at which quantiles are computed is $(1,1)$.}
    \label{fig:GeomQuantilesOrder2-GammaN}}
\end{figure}

We illustrate the application of Theorem~\ref{thm:numerics} in Figures~\ref{fig:GeomQuantilesOrder2-gl} \& \ref{fig:GeomQuantilesOrder2-GammaN}, using the second order characterisation of geometric quantiles, considering the bivariate ($d=2$) standard Gaussian distribution for a light tail example and Pareto($\delta$) for a heavy tail one, for clear comparison. For both bivariate distributions, we assume for simplicity the components (marginals) to be independent and identically distributed. For Pareto($\delta$)-distribution (with marginal survival distribution given by $x^{-\delta}$, for $x>1$), we vary the parameter $\delta$ to have the existence of the $3$rd moment, or lack thereof, but always a finite $2$nd moment ($\delta>2$). The second moment condition is necessary to illustrate our results. Recall Theorem~\ref{thm:numerics}, where under some conditions, together with the existence of the 2nd moment, we have 
\begin{equation}\label{eq:lim-GeomQ-2ndCharac}
    y(\alpha_n):=\|\hat{q}_n(\alpha_n u)\|^2(1-\alpha_n) - \half (\text{tr } \Sigma -\langle\Sigma u,u\rangle) \underset{n\to\infty}{\longrightarrow} 0 \quad \text{a.s.}
\end{equation}
We evaluate $y(\alpha_n)$ as a function of $\alpha_n$, and plot it in Figures~\ref{fig:GeomQuantilesOrder2-gl} \& \ref{fig:GeomQuantilesOrder2-GammaN}. Besides comparing its convergence towards $0$ with regard to the tail behaviour of the underlying probability measure, we also test the rate of convergence depending on the sample size, fixed or growing, considered to compute the quantiles. 
  
\begin{figure}[H] 
    \centering
    \begin{minipage}{0.48\textwidth}
        \includegraphics[width=.7\linewidth]{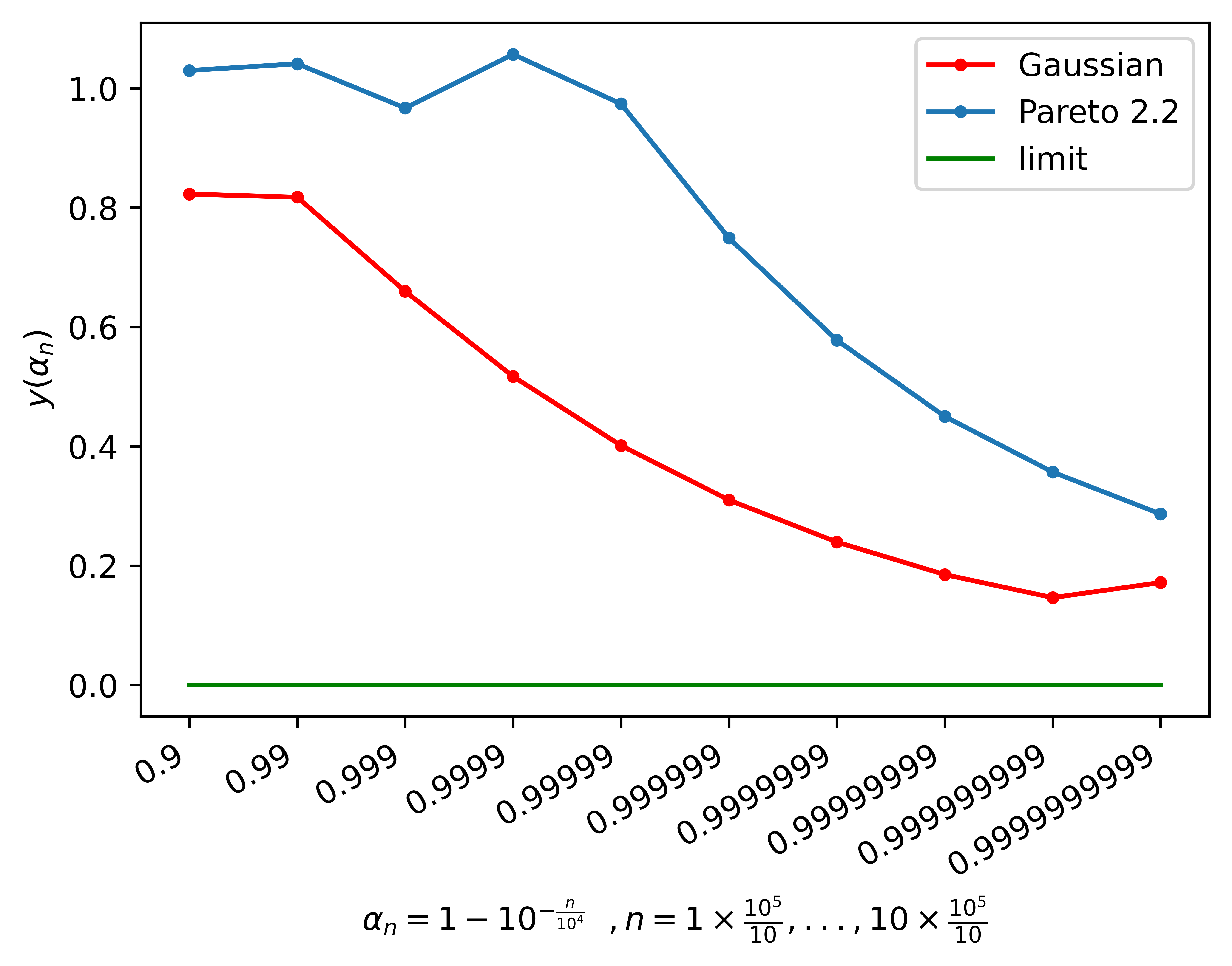}  
    \end{minipage}
    \hfill
    \begin{minipage}{0.48\textwidth}
        \includegraphics[width=.7\linewidth]{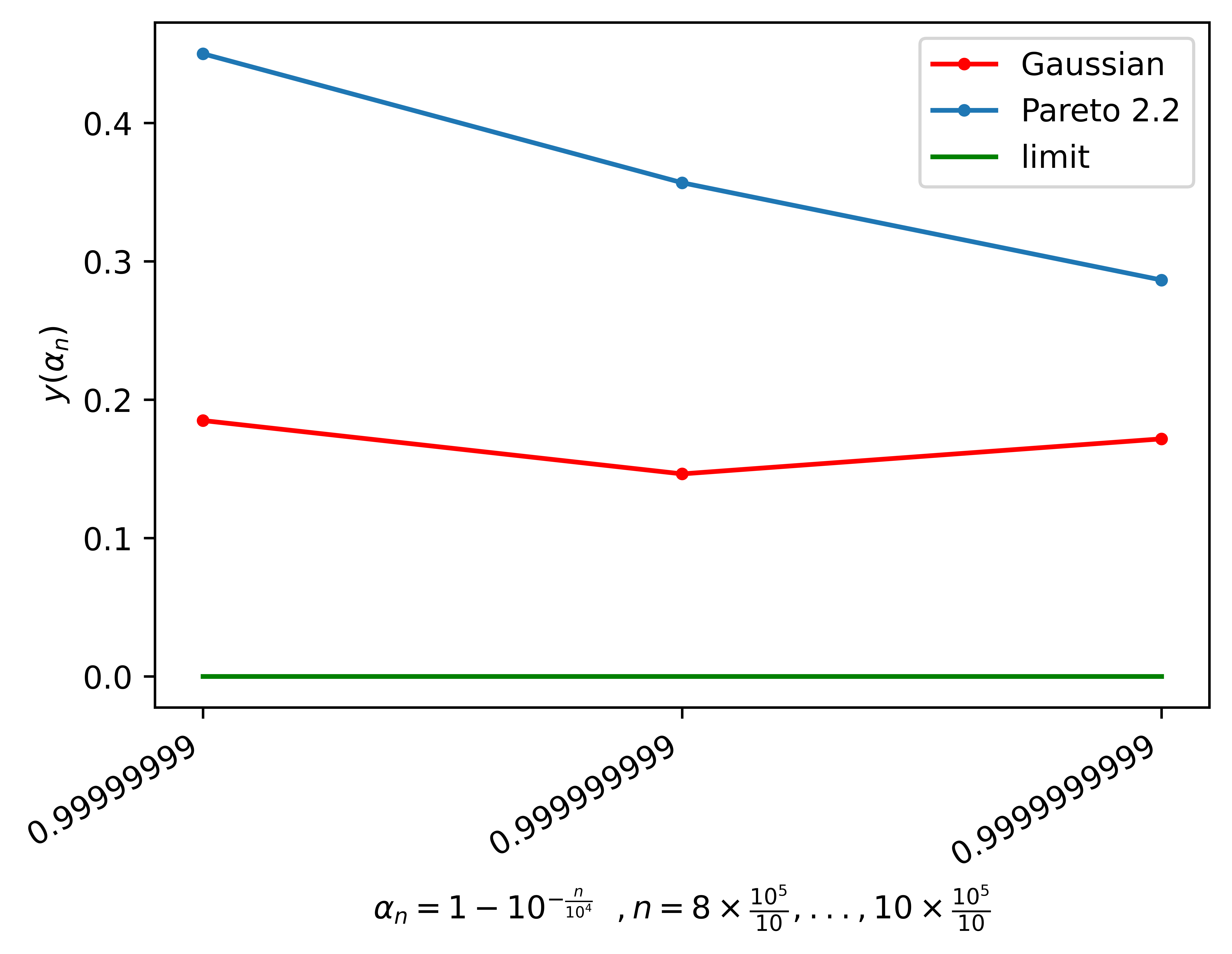}
    \end{minipage}
    \\
    \begin{minipage}{0.48\textwidth}
        \includegraphics[width=.7\linewidth]{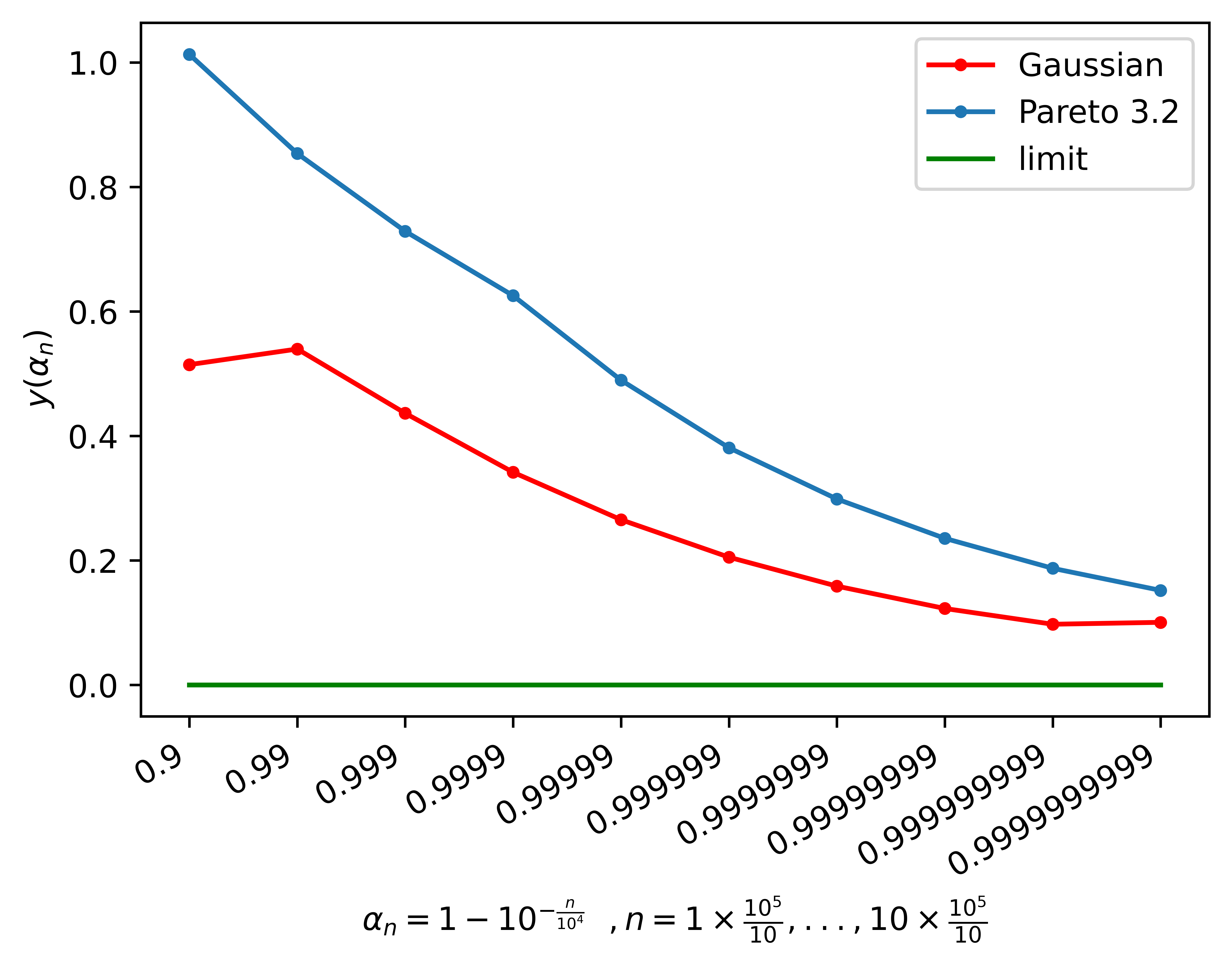}  
    \end{minipage}
    \hfill
    \begin{minipage}{0.48\textwidth}
        \includegraphics[width=.7\linewidth]{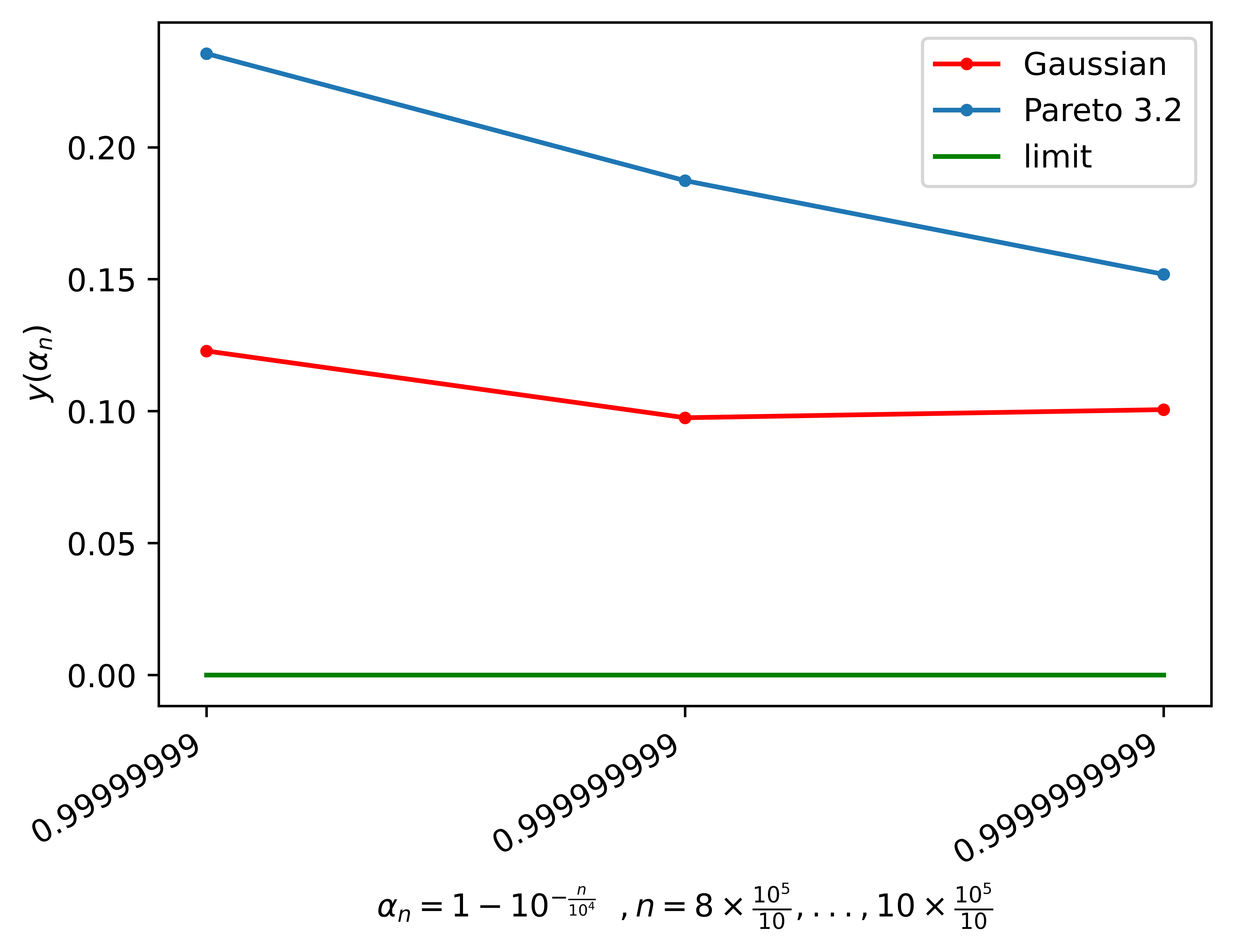}
    \end{minipage}
    \parbox{440pt}{\caption{\sf \small Plot of $y(\alpha_n)$ in terms of $\alpha_n=1-10^{-n\,10^{-4}}$, $n=10^4k$ for $k=1,2,\cdots,n$, for a bivariate centered Gaussian and Pareto($\delta$) distributions sharing the same covariance matrix. First row: $\delta=2.2$, 2nd row: $\delta=3.2$. On each row, the right plot does a zoom on the last 3 values of $\alpha_n$. The $y$-scale is chosen as $y^{1/2}$. }
    \label{fig-samecov}}
\end{figure}

We note here that, although our aim is to study the tail asymptotics of the underlying distribution via the two chosen geometric measures (as seen in Section \ref{ss:geomQ}), the real applications would invariably involve finite samples drawn from the parent distribution, and so do the numerical illustrations.  
Nevertheless, any finite sample, however large, cannot mimic the tail behaviour of the parent distribution. We, therefore, present our illustrations in two scenarios: one with fixed sample size, and the other scenario in which we increase the sample size with increasing $\alpha_n$. In the latter case, we further consider two sub--scenarios: $\alpha_n$'s increase exponentially to $1$, and for each increment in $\alpha_n$, the sample size is increased linearly. In the second sub--scenario, $\alpha_n$ is made to vary as a function of the sample size at an {\it optimal} rate, as identified through Condition~\eqref{eq:cdtionAlpha}. The explicit relationship between tail decay and rate of growth of geometric quantiles is exhibited in Figure~\ref{fig-samecov} by plotting $y(\alpha_n)$ for different distributions: a centered Gaussian and a Pareto($\delta$) family of distributions such that they share the same second moment. We perform this experiment for $\delta=2.2$ and $\delta=3.2$, respectively, to illustrate Theorem~\ref{thm:higher-order} obtained at the population level, where we proved that the third moment can help distinguish the asymptotic behaviour of the geometric quantile whenever the distributions share the same covariance matrix. For simplicity, all the illustrations are shown for the case when $u=(1,1)$. 

It clearly appears in all figures that plotting $y(\alpha_n)$ as a function of $\alpha_n$ gives a visual test to differentiate between distributions based on their tail behaviour.

In Figure~\ref{fig:GeomQuantilesOrder2-gl}, we can clearly identify the difference of behaviour for the various distributions: the lighter the tail, the faster the convergence. 
Note also that increasing the number of observations when increasing $\alpha_n$ does not have much impact on the observed pattern, when comparing the left and right plots of Figure~\ref{fig:GeomQuantilesOrder2-gl}. We surmise that it is due to the large number of (simulated) observations, $10^5$. 
\vspace{-1ex}
\vspace{-1ex}
Let us turn to Figure~\ref{fig:GeomQuantilesOrder2-GammaN}, where the plot $y(\alpha_n)$ is made when taking into account the theoretical bound for $\alpha_n$ given in Theorem~\ref{thm:numerics}, (ii). Note that the plot corresponds  to the average of $y(\alpha_n)$ computed $10$ times, each with different seed.
In Figure~\ref{fig:GeomQuantilesOrder2-GammaN}, we also observe a distinct behaviour of the function $y(\alpha_n)$ according to the tail of the distribution. The convergence towards the limit $0$ is evidently very slow, especially for the Pareto$(2.2)$-distribution (with no third moment). Nevertheless, note that $\alpha_n$ is still `far' from $1$, with the range of $\alpha_n$ corresponding only to the first interspace on the $x$-axis of Figure~\ref{fig:GeomQuantilesOrder2-gl}. 
\vspace{-1ex}
Next, we consider the case where the bivariate Gaussian and Pareto distributions share the same second moment (i.e. the same variance for each component of the two distributions), to match the case of Theorem~\ref{thm:higher-order}. For Pareto, we choose the parameters $2.2$ (no third moment) and $3.3$ (finite third moment), respectively. The two corresponding plots of $y(\alpha_n)$ are given in Figure~\ref{fig-samecov}, with a growing sample as indicated on the plots, and a zoom (right plots) made when $\alpha(n)$ is extremely close to $1$.
In Figure~\ref{fig-samecov}, the curve corresponding to Pareto($\delta$) is closer for $\delta=3.2$ to that of the Gaussian than for $\delta=2.2$, as expected. Larger the $\delta$, faster the convergence. Nevertheless, the function $y(\alpha(n))$ for Pareto($\delta$) remains distinct from that for Gaussian even asymptotically, for $\alpha(n)$ extremely close to $1$, as can be seen on the right plots. While expected for $\delta=3.2$ due to the theoretical result on the rate of convergence under the assumption of a finite third moment, it looks more surprising for $\delta=2.2$. 
%

\subsubsection{Tukey depths}
\label{sss:plots-Tukey}

Now, we illustrate  Theorems~\ref{thm:emp-hdepth-mrv-decay}~\&~\ref{thm:emp-hdepth-light-decay}, considering the asymptotic behaviour of halfspace depth (or Tukey depths) for light and heavy tailed distributions. 
We take the same examples as for geometric quantiles, namely Gaussian and Pareto distributions for comparison, as we have characterised the rate of convergence according to the tail behaviour.
\vspace{-1ex}
\begin{figure}[H]
    \centering
    \begin{minipage}{0.48\textwidth}
        \includegraphics[width=0.8\linewidth]{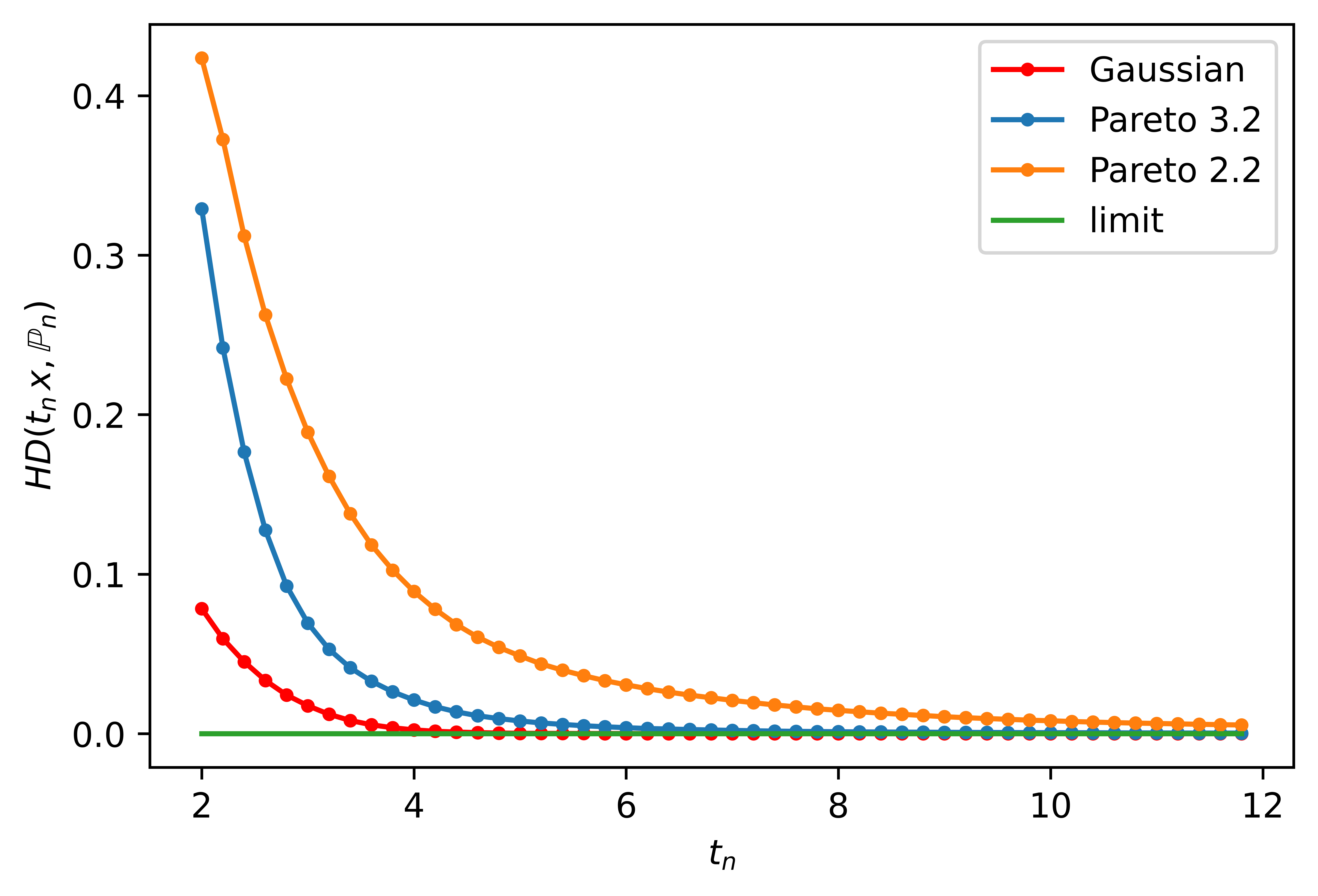}
    \end{minipage}
    \hfill 
    \begin{minipage}{0.48\textwidth}
       \includegraphics[width=0.8\linewidth]{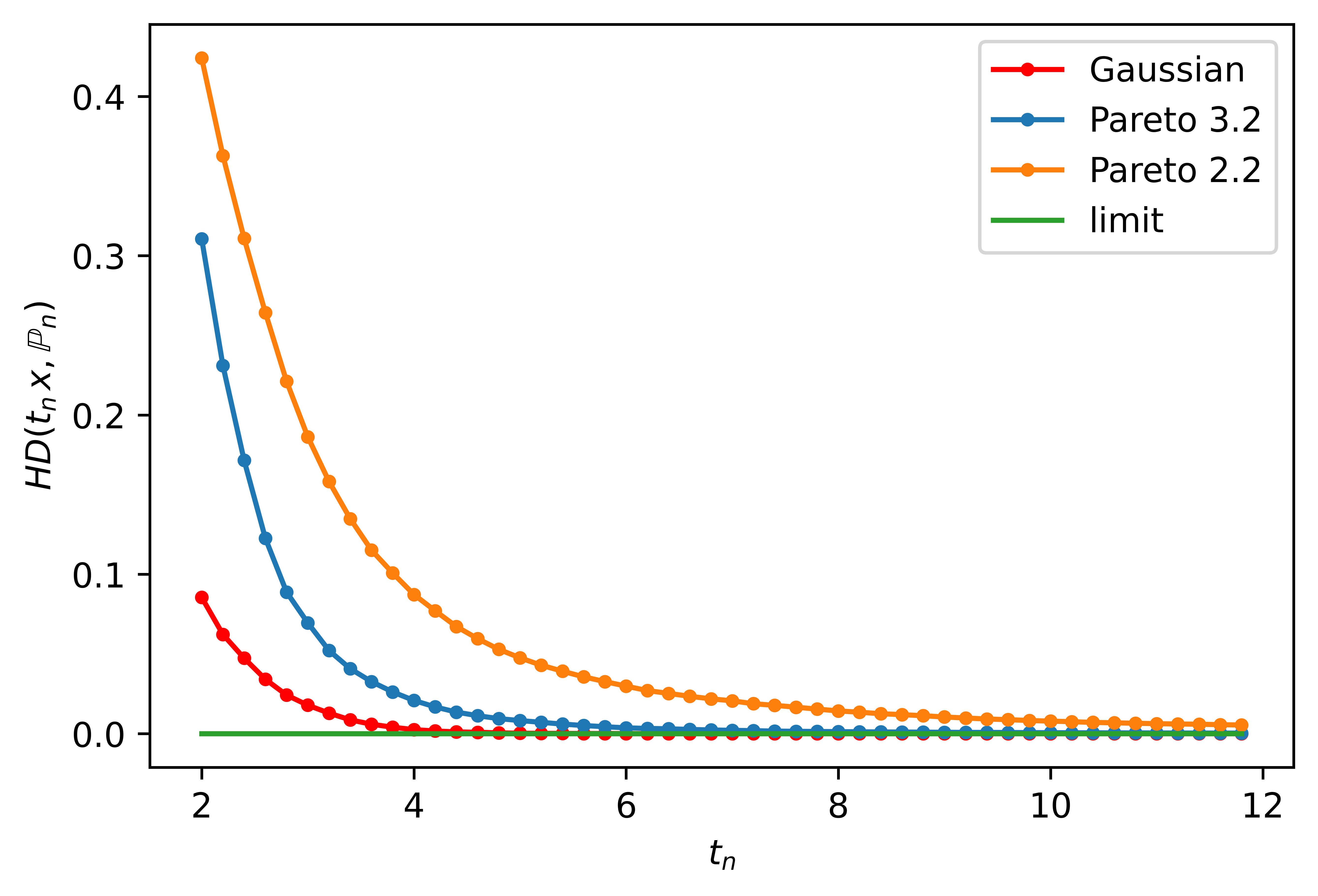}
    \end{minipage}   
\parbox{420pt}{\caption{\sf \small Tukey depths are computed at points in direction $x=(1,1)$ and given in terms of $(t_n)$ growing linearly in $n$ ($t_n=1.8+n.10^{-4}$, with $n=10^5 k/50$, $k=1,2\cdots,50$). Samples are taken from independent bivariate Pareto with parameter $2.2$ and $3.2$, respectively, and Gaussian distribution with diagonal covariance matrix $diag(2,2)$. Number of observations is $10^5$. Left plot: fixed sample. Right plot: growing sample, with a partition of 50}
\label{fig:ratecvdepth-asymp-linearTime}}
\end{figure}
\vspace{-2ex}
In Figure~\ref{fig:ratecvdepth-asymp-linearTime}, we plot $HD(t_n\,x, \pr_n)$ as a function of $t_n$ growing linearly, 
choosing for $x$ the direction $(1,1)$, $\pr_n$ coming from, respectively, bivariate standard Gaussian and Pareto($\delta$) distributions, with independent components, and  $\delta=1.9$, $2.2$ and $3.2$, to span the spectrum from very heavy to moderately heavy tail. 
As previously for the geometric quantiles, we compare $t_n$ when taking a fixed sample (left plot) and a growing one (right plot).
Comparing the different depths according to the type of distributions, from very light to moderate heavy (with no third moment), we clearly observe a different rate of convergence towards $0$. The heavier is the distribution, the slower is the convergence. 

\begin{figure}[H]
    \centering
    \begin{minipage}{0.32\textwidth}
        \includegraphics[width=0.8\linewidth]{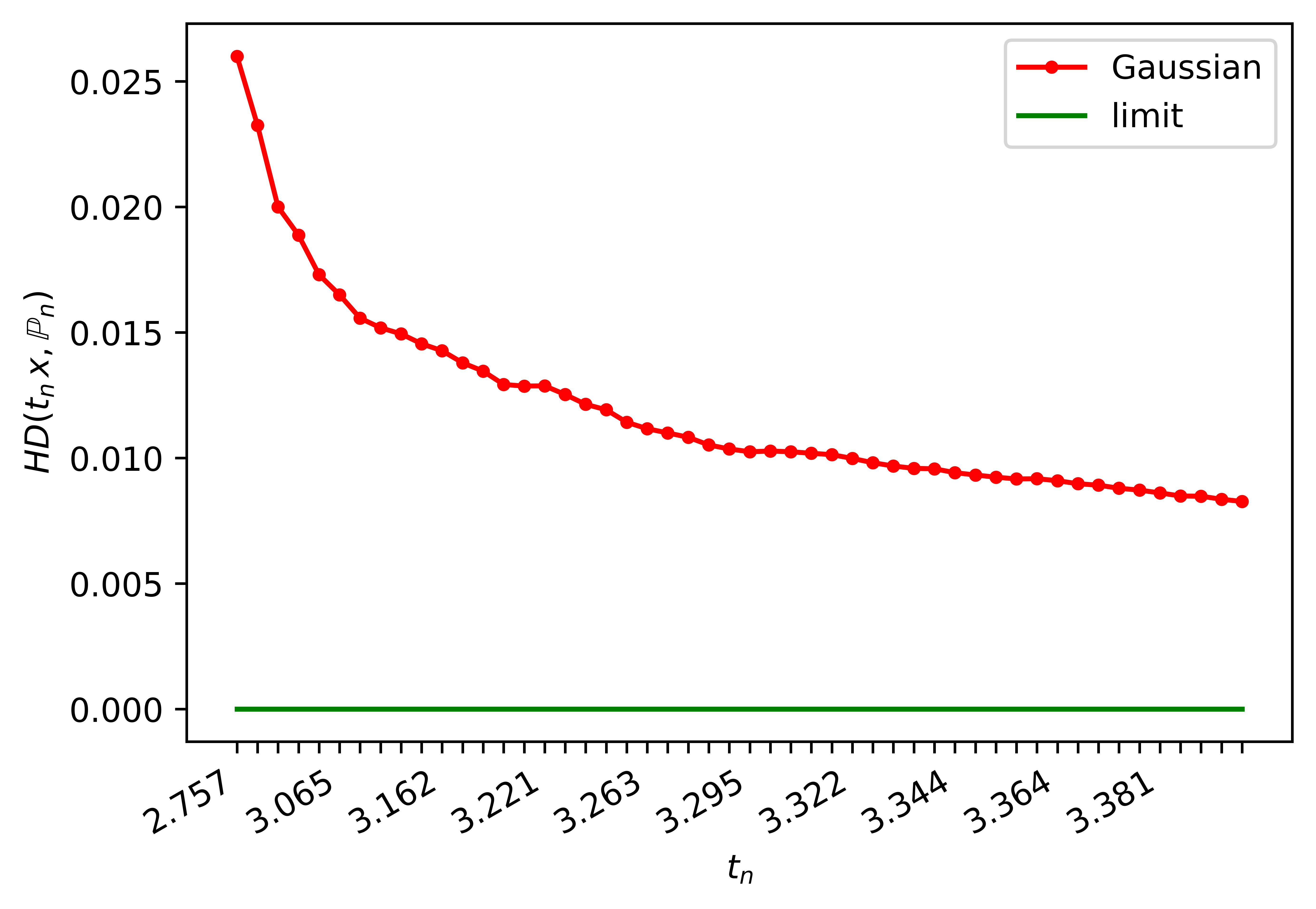}
    \end{minipage}
    \hfill 
    \begin{minipage}{0.32\textwidth}
        \includegraphics[width=0.7\linewidth]{dpth_gs_nba_par}
    \end{minipage}
    \hfill
    \begin{minipage}{0.32\textwidth}
        \includegraphics[width=0.8\linewidth]{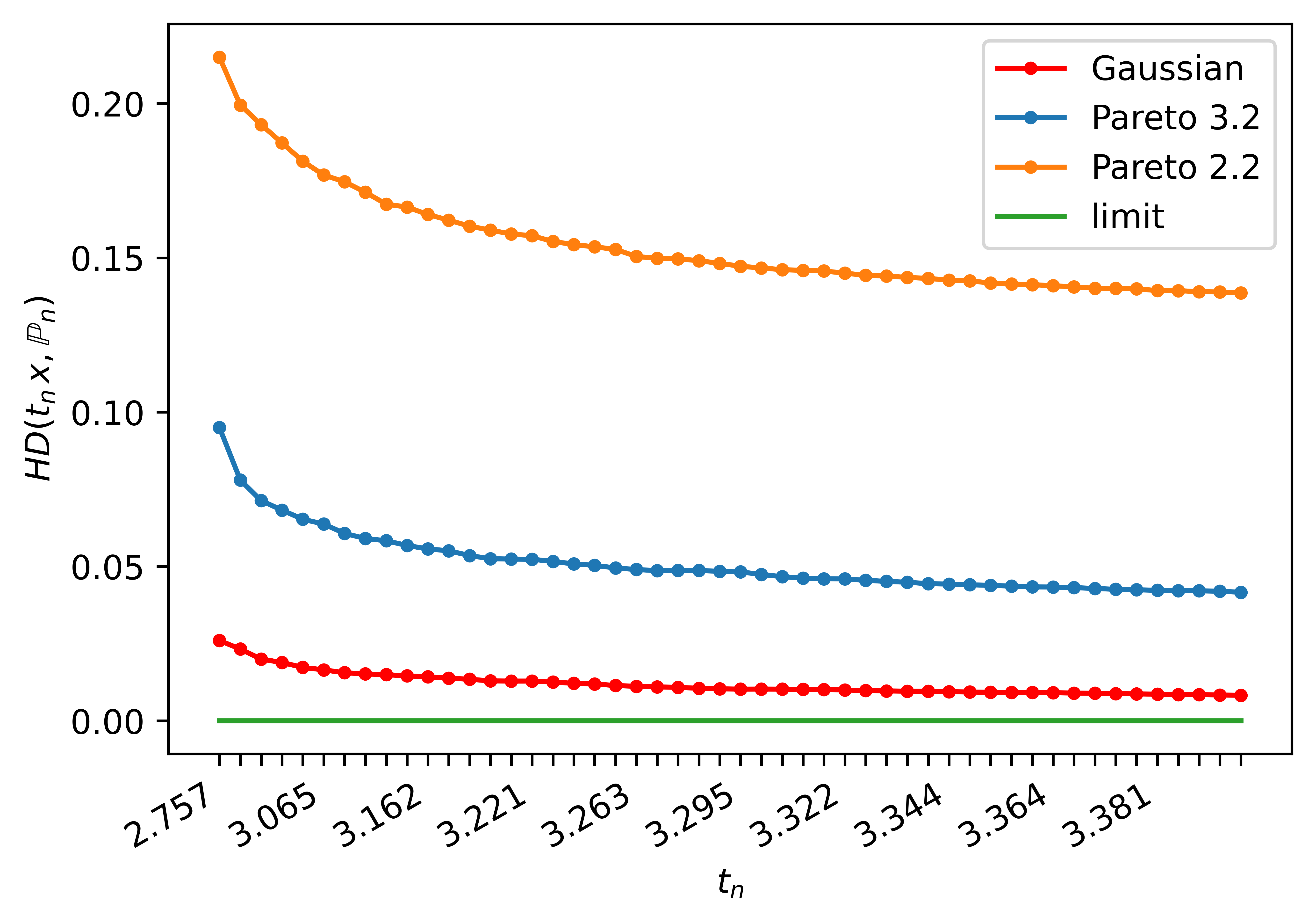}
    \end{minipage}    
\parbox{450pt}{\caption{\sf \small Tukey depths are computed at points in direction $x=(1,1)$ given in terms of $(t_n)$. Samples are taken from independent bivariate Pareto with parameter $1.9$, $2.2$ and $3.2$, respectively, and Gaussian distribution with diagonal covariance matrix $diag(2,2)$. Number of observations is $10^5$. Left plot: Halfspace depth for the Gaussian sample. Middle plot: Comparing the halfspace depth behaviours when considering Pareto($\delta$) samples, varying $\delta$, $t_n$ corresponding to $t_n(\text{Pareto}(3.2))$. Right plot: Halfspace depth behaviours for Gaussian and Pareto($\delta$) samples (choosing $\delta>2$) and for $t_n=t_n(\text{Gaussian})$.}
\label{fig:ratecvdepth-asymp}}
\end{figure}

Next, building on the rate of convergence found in the light tail (see Theorem~\ref{thm:emp-hdepth-light-decay} and Example~\ref{exple:exp-normal} ) and the heavy tail case (see Theorem~\ref{thm:emp-hdepth-mrv-decay}), we plot $HD(t_n x,\pr_n)$ as a function of $t_n$. Given the very different speeds of convergence obtained for the light versus heavy tails, we first give a plot for the Gaussian sample only, then a plot for Pareto($\delta$) samples with varying $\delta$, so that we can appreciate the different behaviour and convergence depending on the heaviness. The sequence $\{t_n\}$ is chosen according to the type of distribution. For the Gaussian case (see Figure~\ref{fig:ratecvdepth-asymp}, left plot), $t_n=\sqrt{\log n}$ (choosing $\beta=1/2$ in Example~\ref{exple:exp-normal}, Gaussian case). For Pareto distributions (see Figure~\ref{fig:ratecvdepth-asymp}, middle plot), we consider 
$t_n=n^{-\frac{\beta}{2\delta}}$, with $0<\beta<1$ also chosen as $1/2$ and the Pareto parameter $\delta$ corresponding to the less heavy, {\it i.e.} $\delta=3.2$ (since the lighter the tail, the faster the convergence towards $0$). Finally, we provide a last plot (see Figure~\ref{fig:ratecvdepth-asymp}, right plot) comparing the Gaussian and Pareto cases, choosing the scaling $t_n$ associated with the Gaussian distribution, for better visualizing the difference of behaviours and speeds of convergence.
\vspace{-1ex}
\vspace{-2ex}
The three plots given in Figure~\ref{fig:ratecvdepth-asymp} highlight the difference of rates of decay of the halfspace depths according to the tail behaviour of the measure. The left and middle plots point out the fast convergence of halfspace depth for the Gaussian sample (decreasing from 2.6\% to less than 1\% (0.83\%) on the given range for $t_n(\text{Gaussian})$), and the impact of the heaviness for the Pareto samples, with a decrease from $26.7\%$ to $16\%$ on the given range for $t_n(\text{Pareto}(3.2))$ for the Pareto with 3rd moment, from 39\% to 30\% for the Pareto(2.5), while from 43\% to 35\% for the heaviest Pareto (with no 2nd moment), hence a very slow decrease compared with Pareto(3.2). The third plot allows for a direct comparison between light and heavy tails, considering the Gaussian scaling for $t_n$; the relation between the rate of decay of the halfspace depth to $0$ and the tail behaviour becomes even more obvious. 

Note that it would have been nice to look at the convergence towards $1$ of the normalized halfspace depth function $HD(t_n x,\pr_n)/N(t_n)$ (rather than $HD(t_n x,\pr_n)$), as a function of $t_n$ , where $N(t_n)$ corresponds to the speed of convergence, namely of order $N(t_n)=n^{-\eta}$ for the Gaussian case (see Example~\ref{exple:exp-normal}) and $N(t_n)=V(t_n)$ as defined in Theorem~\ref{thm:emp-hdepth-mrv-decay}  for the Pareto one. Nevertheless, to observe something informative in terms of convergence, it would require a large number of observations (more than $10^{20}$), which is computationally not feasible with the R-package we are using. We conjecture that a possible way to circumvent this computational hurdle would be to use the geometry of isoquantile (isodepth) contours, as they secrete immense amount of information about the underlying distribution.

\subsubsection{Which descriptive and inferential tool to use, geometric quantile or halfspace depth function?}
\label{sss:ToolComparison}

If one were to wonder which one of the two geometric measures to use as an inferential tool to identify the tail behaviour of the underlying measure, as usual, the answer cannot be binary: it would depend on the data set at hand.

As we know, geometric quantiles uniquely characterise the underlying distribution, whereas halfspace depths characterise the underlying distribution only in certain cases (\cite{Nagy2022}). It seems apparent then to believe that geometric quantiles should be a natural choice as an inferential tool. However, there are clear advantages of using halfspace depths in certain scenarios. For instance, halfspace depths provide an immediate visual estimate of the support of the underlying distribution, which geometric quantiles fail to provide in an easy way. Halfspace depth contours show close resemblance to the isodensity contours, and thus provide a visual tool in identifying the underlying measure, whereas the iso-geometric quantile contours of even an elliptically symmetric distribution do not appear convex near the extremes. Nevertheless, for a `nice' dataset, the geometric quantiles exhibit a much faster convergence than the halfspace depth, and can be computed for dimension higher than $2$, while it becomes quite a computational challenge for halfspace depth.


\section{Proofs}
\label{sec:proofs}

Note that all the (in)equalities involving random variables, appearing in this section, are to be understood as almost sure (in)equalities.


\subsection{Proof of Theorems \ref{thm:higher-order} \& \ref{thm:upperbd-q} -}

{\sf\large Proof of Theorem \ref{thm:higher-order} -}

We begin with the proof of \eqref{eqn:2nd-order} by expressing all the terms involved in terms of the orthonormal basis $\{u,w_1,...,w_{d-1}\}$ of $\real^d$. Let $b(\alpha)$ and $\{\beta_k(\alpha)\}_{k=1}^{d-1}$ be real numbers defined by 
\begin{equation}
\label{eqn:decomposition}
    \frac{q(\alpha u)}{\|q(\alpha u)\|}=b(\alpha )u+\sum _{k-1}^{d-1}\beta_k(\alpha)w_k.
\end{equation}
Therefore,
$$
q(\alpha u) - \|q(\alpha u)\| \, u = \|q(\alpha u)\|\left(b(\alpha ) - 1\right)u + \|q(\alpha u)\| \, \sum _{k-1}^{d-1}\beta_k(\alpha)w_k.
$$
Also, observe that
$$
\E[X - \langle X,u\rangle u] = \sum_{k=1}^{d-1} w_k\,\E\langle X,w_k\rangle.
$$
In view of the above two equations, can write
\begin{eqnarray*}
&& \|q(\alpha u)\|\, \Big[q(\alpha u) - \|q(\alpha u)\| \, u - \E(X - \langle X,u\rangle u) \Big] \\
&=&  \|q(\alpha u)\|^2 \left(b(\alpha ) - 1\right)u + \|q(\alpha u)\|^2 \, \sum _{k-1}^{d-1}\beta_k(\alpha)w_k - \|q(\alpha u)\|\, \sum_{k=1}^{d-1} \E\langle X,w_k\rangle \, w_k,
\end{eqnarray*}
which forms the expression of interest in Equation \eqref{eqn:2nd-order}.

Introducing the desired limit to this expression, we have
\begin{eqnarray}\label{eqn:part-I-II}
&& \|q(\alpha u)\|\, \Big[q(\alpha u) - \|q(\alpha u)\| \, u - \E(X - \langle X,u\rangle u) \Big] \nonumber \\
&& \mbox{} + \half\| \E(X-\langle X,u\rangle u)\|^2 u - \sum_{k=1}^{d-1} \text{cov}(\langle X,u\rangle , \langle X, w_k\rangle) w_k  \nonumber\\
& = & \underbrace{\left[\|q(\alpha u)\|^2 \left(b(\alpha ) - 1\right) + \half\| \E(X-\langle X,u\rangle u)\|^2 \right]}_{I} u \nonumber\\
&& \mbox{} + \underbrace{\sum _{k-1}^{d-1} \Big(\|q(\alpha u)\|^2 \,\beta_k(\alpha) - \|q(\alpha u)\|\,\E\langle X,w_k\rangle -  
\text{cov}(\langle X,u\rangle , \langle X, w_k\rangle)\Big)  \, w_k.}_{II}
\end{eqnarray}

Using Lemma 6.4 in \cite{Girard2017},
we can write
\begin{equation}\label{eqn:part-II}
\|q(\alpha u)\|^2 \beta_k(\alpha) - \|q(\alpha u)\| \E\langle X,w_k\rangle - \text{cov}(\langle X,u\rangle , \langle X,w_k\rangle) \to 0\,\,\,\,\,\,\,\,\,\,\text{ as } \alpha \to 1,
\end{equation}
which concludes that part--$II$ in the above equation converges to $0$ as $\alpha \to 1$.  Next, let us consider part--$I$. First, observe that we can write
\begin{equation}\label{eqn:simplify-part-I}
    \|q(\alpha u)\|^2 (1 - b(\alpha)) = \half \|q(\alpha u)\|^2\, (1-b^2(\alpha))\Bigg(1+\frac{1-b(\alpha)}{1+b(\alpha)}\Bigg).
\end{equation}
Now, notice that $b(\alpha) \to 1$ as $\alpha \to 1$ (by Property \ref{property:set-1}), and $b^2(\alpha) + \sum_{1}^{d-1} \beta^2_k(\alpha) =1$. Therefore, Equation \eqref{eqn:simplify-part-I} simplifies to
\begin{equation}\label{eqn:simplify-part-Ia}
\|q(\alpha u)\|^2 (1 - b(\alpha)) = 
\half \|q(\alpha u)\|^2 \, \sum_{k=1}^{d-1}\beta_k^2 (\alpha) \Big (1 + o(1) \Big ).\\
\end{equation}
However, by the definition of $\beta_k(\alpha)$ and orthogonality of the basis $\{u,w_1,\ldots,w_{d-1}\}$, we have
$$
\|q(\alpha u)\|\,\beta_k(\alpha) = \|q(\alpha u)\|\,\left\langle \frac{q(\alpha u)}{\|q(\alpha u)\|} - \alpha u,w_k\right\rangle.
$$
Now, using Lemma~6.2 in \cite{Girard2017},   
we have
$$ 
\|q(\alpha u)\|\,\beta_k(\alpha) =
\|q(\alpha u)\|\left\langle \frac{q(\alpha u)}{\|q(\alpha u)\|} - \alpha u,w_k\right\rangle \underset{\alpha\to 1}{\longrightarrow} \E\langle X-\langle X,u\rangle u,w_k\rangle.
$$
Therefore,
\begin{equation}
\|q(\alpha u)\|^2 (1-b(\alpha)) \,\,\, \underset{\alpha\to 1}{\longrightarrow} \,\,\, \frac12 \sum_{k=1}^{d-1} \left[\E\langle X-\langle X,u\rangle u,w_k\rangle\right]^2 = \frac12 \|\E\left(X-\langle X,u\rangle u\right) \|^2,
\end{equation}
which proves that part-$I$ of \eqref{eqn:part-I-II} converges to $0$ as $\alpha \to 1$. This, together with \eqref{eqn:part-II}, proves the first part of the theorem.

We now prove the second part of Theorem~\ref{thm:higher-order}.  
Using once again Lemmas 6.2 \& 6.3 in \cite{Girard2017}
we can write
\begin{eqnarray}\label{eqn:10}
& &\|q(\alpha u)\| \Bigg[\|q(\alpha u)\|^2\beta_k^2 (\alpha ) - |\E\langle X,w_k\rangle |^2\Bigg ]\nonumber\\
&=& \Bigg[\|q(\alpha u)\| \Big(\|q(\alpha u)\|\beta_k (\alpha ) - |\E\langle X,w_k\rangle |\Big)\Bigg ]\Bigg[\|q(\alpha u)\|\beta_k (\alpha ) + \E\langle X,w_k\rangle \Bigg] \nonumber\\
&\underset{\alpha\to 1}{\longrightarrow} & \text{cov} (\langle X,u\rangle , \langle X,w_k\rangle )\,\, 2\E\langle X,w_k\rangle.
\end{eqnarray}
Using the decomposition \eqref{eqn:decomposition} and Proposition 6.3 in \cite{Girard2017}
under the assumption $\E\|X\|^3<\infty$, we have
\begin{eqnarray}\label{eqn:11}
&& \|q(\alpha u)\| \Bigg (\|q(\alpha u)^2 (1-\alpha b(\alpha))-\half \E\|X-\langle X,u\rangle u\|^2\Bigg )\nonumber\\
& \underset{\alpha\to 1}{\longrightarrow}& \E\Bigg (\langle X,u\rangle \Big [ \|X-\langle X,u\rangle u\|^2-\langle X,\E(X-\langle X,u\rangle u)\rangle  \Big]\Bigg ) \,
 \stackrel{\Delta}{=}\, f_X(u). 
\end{eqnarray}
Therefore, from Equations \eqref{eqn:10} and \eqref{eqn:11}, we can conclude that
\begin{align}
&\|q(\alpha u)\|\Bigg (\|q(\alpha u)\|^2\Big (1-\alpha b(\alpha) \Big ) - \half \E\|X-\langle X,u\rangle u\|^2 - \half \sum \Big[\|q(\alpha u)\|^2\beta_k^2 (\alpha ) - |\E\langle X,w_k\rangle |^2\Big ]\Bigg ) \nonumber\\
& \underset{\alpha\to 1}{\longrightarrow} \, f_X(u) - \sum_{k=1}^{d-1} \E\langle X,w_k\rangle\,\,\text{cov}(\langle X,u\rangle , \langle X,w_k\rangle. \label{eqn:12}
\end{align}
But, notice that the LHS of Equation~\eqref{eqn:12} can be rewritten as 
$$\|q(\alpha u)\|\Bigg (\|q(\alpha u)\|^2\half\Big (1-\alpha^2+( b(\alpha)-\alpha)^2 \Big ) - \half \sum_{k=1}^{d-1} \text{var} \langle X,w_k \rangle\Bigg),$$
by observing that $(1-b^2(\alpha))=\sum_{k=1}^{d-1} \beta^2_k(\alpha)$, and after some further algebraic manipulation.

Additionally, note that $\|q(\alpha u)\|^3(\alpha - b(\alpha))^2 \to 0$ since $\|q(\alpha u)\|(\alpha -b(\alpha)) \to 0$ and $\|q(\alpha u)\|^2(\alpha -b(\alpha)) \to 0$, as a result of 
Proposition 6.3 in \cite{Girard2017}.
Therefore, we deduce that
\begin{equation*}
\|q(\alpha u)\|\Bigg (\|q(\alpha u)\|^2 \Big (1-\alpha\Big ) - \left( \text{tr} \Sigma - \langle\Sigma u,u\rangle \right)\Bigg)
\to f_X(u) - \sum_{k=1}^{d-1} \E\Bigg (\langle X,w_k\rangle)\,\,\text{cov}(\langle X,u\rangle , \langle X,w_k\rangle\Bigg ).
\end{equation*}

Simplifying the limit, we have
	\begin{eqnarray*}
        &&  \|q_X(\alpha u)\|\Big[\|q_X(\alpha u)\|^2 (1 - \alpha ) -  \half (\text{tr } \Sigma - \langle\Sigma u,u\rangle)\Big] \, \longrightarrow \nonumber \\
        && \E\Big(\langle X,u\rangle \big[\|X-\langle X,u\rangle u\|^2 - \langle X,\E\left(X-\langle X,u\rangle u\right)\rangle \big]\Big) -\sum_{k=1}^{d-1} \text{cov}(\langle X,u\rangle, \langle X,w_k\rangle)\, \E(\langle X,w_k\rangle ).\nonumber \\
        && = \sum_{k=1}^{d-1}\E\Big( \langle X,u\rangle \langle X,w_k\rangle^2 - \langle X,u\rangle\langle X,w_k\rangle \E\langle X,w_k\rangle \Big) -\sum_{k=1}^{d-1} \text{cov}(\langle X,u\rangle, \langle X,w_k\rangle)\, \E(\langle X,w_k\rangle ) \\
        && = \sum_{k=1}^{d-1}\text{cov}\big(\langle X,u\rangle\langle X,w_k\rangle, \langle X,w_k\rangle\big) -\sum_{k=1}^{d-1} \text{cov}(\langle X,u\rangle, \langle X,w_k\rangle)\, \E(\langle X,w_k\rangle )
        	\end{eqnarray*}
which proves \eqref{eqn:3rd-order}. 

\hfill \qed

{\sf\large Proof of Theorem \ref{thm:upperbd-q} -}\\
Recall that $q(\alpha u) =\argmin f_{\alpha u}(q)$, with $f_{\alpha u}(q)=\E\left[ \|X-q\|-\|X\|-\langle\alpha u,q \rangle\right]$. Observe that, since $f_{\alpha u} (0) = 0$, we must have 
\begin{equation}\label{eqn:fau-upper-bd}
f_{\alpha u} (q(\alpha u)) \le f(0) = 0.
\end{equation} 

For any $q$ and $k_{\alpha}>0$, we can obtain the following lower bound by applying triangle inequality,
\begin{eqnarray}\label{eqn:fau-lower-bd}
f_{\alpha u}(q) &\ge& \E\left[\{\|q\|-2\|X\|\}\one_{\|X\|\leq k_{\alpha}}\right]+\E\left[ -\|q\|\one_{\|X\|>k_{\alpha}}\right]-\alpha \|q\|\nonumber \\
& \ge & \|q\|\left\{ 1 -\alpha - 2\pr(\|X\| > k_{\alpha})\right\}  - 2\,k_{\alpha} \pr(\|X\| \le k_{\alpha})
\end{eqnarray}

Combining \eqref{eqn:fau-upper-bd} and \eqref{eqn:fau-lower-bd}, we conclude that
$$\|q(\alpha u) \| \left\{ 1 -\alpha - 2\pr(\|X\| > k_{\alpha})\right\}  - 2\,k_{\alpha} \pr(\|X\| \le k_{\alpha}) \le 0$$
which proves the assertion of the theorem. \hfill \qed


\subsection{Proof of Theorem~\ref{thm:asymp-depth-ell-symm}}

Let us consider $x,y\in\real^d$ such that $\|\Sigma x\| = \|\Sigma y\|$. Let $A$ be the orthogonal matrix such that $\Sigma x = A\Sigma y$. Setting $A^* = \Sigma^{-1} A\Sigma$, we have $x = A^*y$, and $\text{det}(A^*) = 1$. Recalling  the Definition~\ref{def:HD} of halfspace depth, 
and setting $$\displaystyle p_y = \text{arg}\min_{|p|\neq 0}\int_{\{w:\langle w-t y,p \rangle\geq 0\}}\!\!\! f(w) dw$$ we have, 
$$ \frac{HD(tx;\pr)}{HD(ty;\pr)} =  \frac{\inf_{\|p\|\neq 0} \int_{\{w: \langle w-tx,p \rangle\geq 0\}}f(w) dw}{\int_{\{w:\langle  w- ty, p_y \rangle \geq 0\}}f(w) dw} \le \frac{\int_{\{w: \langle w-tx,Bp_y \rangle\geq 0\}}f(w) dw}{\int_{\{w:\langle  w- ty, p_y \rangle \geq 0\}}f(w) dw},$$
for any linear transformation $B$.

Now, using the transformation $w=tz$, we have
$$ \frac{HD(tx;\pr)}{HD(ty;\pr)} \le  \frac{\int_{\{z: \langle z-x,Bp_y \rangle\geq 0\}}f(tz) dz}{\int_{\{z:\langle  z- y, p_y \rangle \geq 0\}}f(tz) dz}.$$
Recalling that $x=A^*y$, we introduce the transformation $z=A^*u$ in the numerator to observe that 
$$ \frac{HD(tx;\pr)}{HD(ty;\pr)} \le  \frac{\int_{\{u: \langle A^*(u-y),Bp_y \rangle\geq 0\}}f(tA^*u) du}{\int_{\{z:\langle  z- y, p_y \rangle \geq 0\}}f(tz) dz}.$$
Since the above inequality is satisfied for any nonsingular $B$, we set $B^T=(A^*)^{-1}$, to obtain
$$ \frac{HD(tx;\pr)}{HD(ty;\pr)} \le  \frac{\int_{\{u: \langle u-y,p_y \rangle\geq 0\}}f(tA^*u) du}{\int_{\{z:\langle  z- y, p_y \rangle \geq 0\}}f(tz) dz}.$$ Now while observing that $\|\Sigma A^*u\| = \|\Sigma u\|$, we invoke the assumption of asymptotic elliptical symmetry of $f$, to conclude that there exists $t_0$ large enough, such that
$$ \frac{HD(tx;\pr)}{HD(ty;\pr)} \le  \frac{\int_{\{u: \langle u-y,p_y \rangle\geq 0\}} \left[\frac{f(tA^*u)}{f(tu)}\right] f(tu) du}{\int_{\{z:\langle  z- y, p_y \rangle \geq 0\}}f(tz) dz} \le (1+\epsilon),\,\,\,\forall t\ge t_0$$

Using similar arguments, we can also conclude that
$$\frac{HD(tx;\pr)}{HD(ty;\pr)} \ge (1-\eps), \,\,\,\,\,\,\,\,\forall t\ge t_0$$
which concludes the result. \hfill \qed

\subsection{Proof of Theorem~\ref{thm:hd-decay-light-tail}}
\label{ss:theo2.12}

\begin{itemize}
    \item[(i)] The proof of the upper bound is based on a simple application of Markov inequality. Indeed, we can write for any $Y$ with distribution $\pr$,
$$
HD(tx,\pr) 
= \inf_{|p|=1}\pr[e^{\langle Y,p\rangle }\ge e^{t \langle x,p\rangle }]
\le \inf_{|p|=1} e^{-t\langle x,p\rangle }\;\E[e^{\langle Y,p\rangle}],
$$
from which the first result follows.
    \item[(ii)]
    Recall the Tukey definition of halfspace depth for any measure $\pr$ given in Definition~\ref{def:HD},  
    but also the alternative way to rewrite this definition (see \eqref{eqn:hdepth-project}), using the standard parameterisation of a halfspace in terms of its distance from origin, and its normal vector, namely:
    $$
    HD(x,\pr) = \inf_{\theta \in S^{d-1}} \pr\left[\langle \theta,X\rangle \ge \langle \theta,x\rangle\right].
    $$
    Let $c>0$ be such that, for all sufficiently large values of $R$,
    \begin{equation}
    \inf_{\theta \in S^{d-1}}\dfrac{\pr(\langle \theta, X\rangle \geq R)}{f(R)} >c>0.
    \end{equation}    
Invoking the definition of halfspace depth, we can write
$$
HD(t x, \pr)= \inf_{\theta \in S^{d-1}} \pr(\langle \theta, X\rangle \geq t \langle\theta, x\rangle) \geq \inf_{\theta \in S^{d-1}} \pr(\langle\theta, X\rangle \geq t \|x\|) \geq c\, f(t\|x\|),
$$
hence the lower bound.
\end{itemize}


\subsection{Proof of Theorem~\ref{thm:empirical-geom-quant-asymp1}}
\label{ss:theo3.2}

We shall begin with some regularity estimates for the function 
$$(x,q)\mapsto \frac{x-q}{\|x-q\|},$$
which is going to play a crucial role in the proof of Theorem \ref{thm:empirical-geom-quant-asymp1}.

\begin{lemma}
\label{lem:eqcontinuity}
Let $\{q_n\}_{n\ge 1} \subset\real^d$ be such that $\displaystyle q_n \underset{n\to\infty}{\longrightarrow} q_{\infty}$. Define $f: \real^d \times \real ^d \to \real ^d$ by,
\begin{equation}
    f(x,q) = 
     \begin{cases}
      \frac{x-q}{\|x-q\|} &x \neq q\\
       0 &x=q.\\
     \end{cases}
\end{equation}
Then, for $n$ large enough, $$\|f(x,q_n)-f(x,q_{\infty})\|<K_M\|q_n-q_{\infty}\|, \,\,\,\,\,\,\text{ whenever }\,\,\,\,\, \|x-q_{\infty}\|>M>0,$$ where $K_M$ is a constant depending on $M$.
\end{lemma}
\noindent{\bf Proof:} Observe that, for $x \neq q $,
\begin{equation}
    \nabla_q f(x,q)=-\frac{1}{\|x-q\|}\Bigg[I_{d \times d}-\frac{(x-q)^T(x-q)}{\|x-q\|^2}\Bigg].
\end{equation}
If $\|x-q_{n}\|>\frac{M}{2}$, then writing $\|\cdot\|_F$ for the Frobenius norm, we have
\begin{equation}
    \|\nabla_q f(x,q_n)\|_F \leq \frac{2}{M} \times d(d+1).
\end{equation} 
Since $\displaystyle q_n \underset{n\to\infty}{\to} q_{\infty}$, there exists $k_0$ large enough such that $\|q_n-q_{\infty}\|<\frac{M}{2}, \,\, \forall n \geq k_0$. Then, whenever $x$ satisfies $\|x-q_{\infty}\|>M$, it follows that, for all $n\ge k_0$, $\|x-q_n\|>\frac{M}{2}$. Thereby for such an $x$,
\begin{equation*}
    \|f(x,q_n)-f(x,q_{\infty})\|<\frac{2d(d+1)}{M}\|q_n-q_{\infty}\| ,\,\, \forall n \geq k_0,
\end{equation*}
which concludes the proof. \hfill \qed

In addition to the above observation, we also need primary estimates for the asymptotic behaviour of $\displaystyle \sum_{i=1}^n (X_i - \hat q_n(\alpha_n u))\|X_i - \hat q_n(\alpha_n u)\|^{-1}$, for $u\in S^d$, and $\alpha_n\in (0,1)$. Recall from Theorem 2.12 in \cite{Chaudhury1996} that
$$
    \left|\frac{1}{n}\sum_{\substack{{i=1}\\{X_i\neq \hat q_n(\alpha_n u)}}}^n \frac{X_i-\hat q_n(\alpha_n u)}{\|X_i-\hat q_n(\alpha_n u)\|} + \alpha_n u\right| \leq \frac{1}{n}\#\{i:X_i=\hat q_n(\alpha_n u)\}.
$$
Since $\pr \left(\#\{i:X_i=\hat q_n(\alpha_n u)\} > 1\right)=0$ for samples drawn from any continuous distribution, we have
$$
    \lim_{n \to \infty}\Bigg|\frac{1}{n}\sum_{\substack{{i=1}\\{X_i\neq \hat q_n(\alpha_n u)}}}^n \frac{X_i-\hat q_n(\alpha_n u)}{\|X_i-\hat q_n(\alpha_n u)\|} + \alpha_n u\Bigg|=0, \,\,\text{ a.s.}
$$
In other words,
\begin{equation}
\label{eqn:derivative}
\lim _{n \to \infty} \frac{1}{n}\sum_{\substack{{i=1}\\{X_i\neq \hat q_n(\alpha_n u)}}}^n \frac{X_i-\hat q_n(\alpha_n u)}{\|X_i-\hat q_n(\alpha_n u)\|} =-u,\,\,\,\, \text{a.s.}
\end{equation}

\noindent{\bf Proof of Theorem \ref{thm:empirical-geom-quant-asymp1}.}
We shall prove the result by contradiction. For part $(i)$, let us assume that $\{\hat q_n(\alpha_n u)\}$ is a bounded sequence; so, we can always extract a convergent subsequence. To avoid complexity of notation, let us consider 
$\displaystyle \lim_{n \to \infty} \|\hat q_{n}({\alpha_n}u)\|=q_{\infty}$ a.s. and $\|q_{\infty}\|< +\infty$. We can write
\begin{eqnarray*}
&& \Bigg|\frac{1}{n} \sum_{\substack{{i=1}\\{X_i\neq \hat q_n(\alpha_n u)}}}^n \frac{X_i-\hat q_n(\alpha_n u)}{\|X_i-\hat q_n(\alpha_n u)\|} - \frac{1}{n}\sum_{i=1}^n \frac{X_i-q_{\infty}}{\|X_i-q_{\infty}\|}\Bigg|\\
& \leq & \frac{1}{n}\sum_{\substack{{i=1}\\{X_i\neq \hat q_n(\alpha_n u)}}}^n \Bigg| \frac{X_i-\hat q_n(\alpha_n u)}{\|X_i-\hat q_n(\alpha_n u)\|}-\frac{X_i-q_{\infty}}{\|X_i-q_{\infty}\|}\Bigg| +\frac{1}{n}\,\,\,\,\,\,(\text{since }\pr(\#\{i: X_i=\hat q_n(\alpha_n u)\}>1)=0)\\
& = & \frac{1}{n}\sum_{\substack{{i=1}\\{X_i\neq \hat q_n(\alpha_n u)}}}^n \Bigg| \frac{X_i-\hat q_n(\alpha_n u)}{\|X_i-\hat q_n(\alpha_n u)\|}-\frac{X_i-q_{\infty}}{\|X_i-q_{\infty}\|}\Bigg| \mathbbm{1} _{\{|X_i-q_{\infty}|>M>0\}}\\
&& \mbox{} +\frac{1}{n}\sum_{\substack{{i=1}\\{X_i\neq \hat q_n(\alpha_n u)}}}^n \Bigg| \frac{X_i-\hat q_n(\alpha_n u)}{\|X_i-\hat q_n(\alpha_n u)\|}-\frac{X_i-q_{\infty}}{\|X_i-q_{\infty}\|}\Bigg| \mathbbm{1}_{\{|X_i-q_{\infty}|\leq M\}}+\frac{1}{n}\\
& \leq & K_M \|\hat q_n(\alpha_n u)-q_{\infty}\| + \frac{1}{n}\times 2\sum_{\substack{{i=1}\\{X_i\neq \hat q_n(\alpha_n u)}}}^n \mathbbm{1}_{\{|X_i-q_{\infty}|\leq M\}} + \frac1n,
\end{eqnarray*}
where the first term is coming from Lemma \ref{lem:eqcontinuity}, and the second term from the triangle inequality. Let us set $\epsilon >0$ small. Then, by Glivenko--Cantelli  theorem, there exists a small enough $M$ such that,
$$
    \frac{1}{n}\sum_{\substack{{i=1}\\{X_i\neq \hat q_n(\alpha_n u)}}}^n \mathbbm{1}_{\{|X_i-q_{\infty}|\leq M\}}<\frac{\eps}{6}
$$
almost surely for sufficiently large $n$. Now choose $n$ large so that $\|\hat q_n(\alpha_n u)-q_{\infty}\| < \frac{\eps}{3K_M}$, and $1/n < \eps/3$.
Therefore, we obtain
$$
    \Bigg|\frac{1}{n}\sum_{\substack{{i=1}\\{X_i\neq \hat q_n(\alpha_n u)}}}^n \frac{X_i-\hat q_n(\alpha_n u)}{\|X_i-\hat q_n(\alpha_n u)\|} - \frac{1}{n}\sum_{i=1}^n \frac{X_i-q_{\infty}}{\|X_i-q_{\infty}\|}\Bigg|< \eps.
$$
Since $\epsilon$ is arbitrary, we conclude that
\begin{equation}
\label{eqn:samp conv}
    \lim_{n\to \infty} \Bigg|\frac{1}{n}\sum_{\substack{{i=1}\\{X_i\neq \hat q_n(\alpha_n u)}}}^n \frac{X_i-\hat q_n(\alpha_n u)}{\|X_i-\hat q_n(\alpha_n u)\|} - \frac{1}{n}\sum_{i=1}^n \frac{X_i-q_{\infty}}{\|X_i-q_{\infty}\|}\Bigg|=0,\,\,\,\,\,\,\text{a.s.}
\end{equation}
Now, using the strong law of large numbers, it comes
\begin{equation}
\label{eqn:slln}
\lim_{n \to \infty} \frac{1}{n}\sum_{i=1}^n \frac{X_i-q_{\infty}}{\|X_i-q_{\infty}\|} = \E\Bigg[ \frac{X-q_{\infty}}{\|X-q_{\infty}\|} \Bigg],\,\,\,\,\,\,\text{a.s.}
\end{equation}
Combining Equations \eqref{eqn:derivative}, \eqref{eqn:samp conv} and \eqref{eqn:slln}, we obtain, almost surely,
$$
    -u=\lim _{n \to \infty}\frac{1}{n}\sum_{\substack{{i=1}\\{X_i\neq \hat q_n(\alpha_n u)}}}^n \frac{X_i-\hat q_n(\alpha_n u)}{\|X_i-\hat q_n(\alpha_n u)\|}= \E\Bigg[ \frac{X-q_{\infty}}{\|X-q_{\infty}\|} \Bigg],
$$
implying that
\begin{equation}\label{eqn:to-contradict}
    \|-u\|=1=\left\|\E\left[ \frac{X-q_{\infty}}{\|X-q_{\infty}\|} \right]\right\|.
\end{equation}
Since the distribution of $X$ does not lie on a single straight line, therefore,
\begin{equation*}
\left\|\E\Bigg[ \frac{X-q_{\infty}}{\|X-q_{\infty}\|} \Bigg]\right\|
< \E\left\| \frac{X-q_{\infty}}{\|X-q_{\infty}\|} \right\| 
= 1,
\end{equation*}
which contradicts \eqref{eqn:to-contradict}. Therefore, 
$\displaystyle \|\hat q_n(\alpha_n u)\| \overset{a.s.}{\underset{n\to\infty}{\rightarrow}} \infty$, thus proving part $(i)$ of the theorem.

Let us now turn to $(ii)$ of the theorem, which we again prove via contradiction. Assume that there exists a subsequence of $\left\{\frac{\hat q_{n}({\alpha_n}u)}{\|\hat q_{n}({\alpha_n}u)\|}\right\}$ that converges, as $n\to\infty$, to some $v\in S^{d-1}$ such that $v \neq u$. For simplicity, let us keep the same notation for the subsequence, {\it i.e.}, $\left\{\frac{\hat q_{n}({\alpha_n}u)}{\|\hat q_{n}({\alpha_n}u)\|}\right\}_{n\ge 1}$. 
           
We have already established in Equation~\eqref{eqn:derivative} that
\begin{equation}\label{eqn:derivative-1}
\lim _{n \to \infty} \frac{1}{n}\sum_{i=1}^n \frac{X_i-\hat q_n(\alpha_n u)}{\|X_i-\hat q_n(\alpha_n u)\|}\mathbbm{1}_{\{X_i\neq \hat q_n(\alpha_n u)\}}=-u, \,\,\,\,\,\,\,\,\,\text{a.s.}
\end{equation}
Let us prove that the left hand side converges to $-v$ almost surely, which will imply $u=v$, thus contradicting the assumption. We can write, using similar  arguments as previously,
\begin{eqnarray}\label{eqn:proof-3.6(b)}
&& \left|\frac{1}{n}\sum_{i=1}^n \frac{X_i-\hat q_n(\alpha_n u)}{\|X_i-\hat q_n(\alpha_n u)\|}\mathbbm{1}_{\{X_i\neq \hat q_n(\alpha_n u)\}} + \frac{1}{n}\sum_{i=1}^n v \right| 
\le  \frac1n + \frac1n \sum_{i=1}^n \left|\frac{X_i-\hat q_n(\alpha_n u)}{\|X_i-\hat q_n(\alpha_n u)\|} + v \right|\,\mathbbm{1}_{\{X_i\neq \hat q_n(\alpha_n u)\}}\nonumber \\
&& \leq \, \frac1n + \frac{1}{n}\sum_{i=1}^n \left|\frac{\frac{X_i}{\|\hat q_n(\alpha_n u)\|}-\frac{\hat q_n(\alpha_n u)}{\|\hat q_n(\alpha_n u)\|}}{\left|\frac{X_i}{\|\hat q_n(\alpha_n u)\|}-\frac{\hat q_n(\alpha_n u)}{\|\hat q_n(\alpha_n u)\|}\right|}+v \right|\mathbbm{1}_{\{X_i\neq \hat q_n(\alpha_n u),\, |X_i|\leq M\}} + 2 \, \frac{1}{n}\sum_{i=1}^n  \mathbbm{1}_{\{|X_i|> M\}}.
\end{eqnarray}
Fix $\epsilon>0$ and choose $N_0$ such that
\begin{equation}\label{eqn:proof-3.6(b)-choose-n}
n\epsilon >4,\,\,\,\,\forall n\ge N_0.
\end{equation}
Next, observe that
$$
\frac{1}{n}\sum_{i=1}^n  \mathbbm{1}_{\{|X_i|> M\}} \le \frac{1}{n}\sum_{i=1}^n  \left|\mathbbm{1}_{\{|X_i|> M\}} - (1-F_{|X|}(M))\right| + (1-F_{|X|}(M)).
$$
Now, invoking the Glivenko--Cantelli theorem again, there exist $M_0>0$ and integer $N_1$ such that 
\begin{equation}\label{eqn:proof-3.6(b)-GC}
\frac{1}{n}\sum_{i=1}^n  \mathbbm{1}_{\{|X_i|> M\}} \le \frac{\eps}8\,\,\,\,\,\,\forall M> M_0\,\,\,\text{ and } n>N_1.
\end{equation}
Thus, combining Equations \eqref{eqn:proof-3.6(b)-choose-n} and \eqref{eqn:proof-3.6(b)-GC} with \eqref{eqn:proof-3.6(b)}, we have
\begin{eqnarray}\label{eqn:proof-3.6(b)-last}
&& \left|\frac{1}{n}\sum_{i=1}^n \frac{X_i-\hat q_n(\alpha_n u)}{\|X_i-\hat q_n(\alpha_n u)\|}\mathbbm{1}_{\{X_i\neq \hat q_n(\alpha_n u)\}} + \frac{1}{n}\sum_{i=1}^n v \right|\nonumber \\
& \leq & \frac{\eps}2 + \frac{1}{n}\sum_{i=1}^n \left|\frac{\frac{X_i}{\|\hat q_n(\alpha_n u)\|}-\frac{\hat q_n(\alpha_n u)}{\|\hat q_n(\alpha_n u)\|}}{\left|\frac{X_i}{\|\hat q_n(\alpha_n u)\|}-\frac{\hat q_n(\alpha_n u)}{\|\hat q_n(\alpha_n u)\|}\right|}+v \right|\mathbbm{1}_{\{X_i\neq \hat q_n(\alpha_n u),\, |X_i|\leq M\}}.
\end{eqnarray}

Now recall that we have proved that $\displaystyle\|\hat q_n(\alpha_n u)\| \overset{a.s.}{\underset{n\to\infty}{\to}} \infty$. Therefore, under our assumption that $\displaystyle \frac{\hat q_{n}({\alpha_n}u)}{\|\hat q_{n}({\alpha_n}u)\|}\overset{a.s.}{\underset{n\to\infty}{\longrightarrow}} v$, we have that, whenever $|X_i|\leq M$ and $X_i\neq \hat q_n(\alpha_n u)$, 
$$
\left(\frac{X_i}{\|\hat q_n(\alpha_n u)\|}-\frac{\hat q_n(\alpha_n u)}{\|\hat q_n(\alpha_n u)\|}\right)\left|\frac{X_i}{\|\hat q_n(\alpha_n u)\|}-\frac{\hat q_n(\alpha_n u)}{\|\hat q_n(\alpha_n u)\|}\right|^{-1} \overset{a.s.}{\underset{n\to\infty}{\longrightarrow}} -v 
$$ 
and this convergence is uniform in $i\leq n$. Therefore, the middle term in \eqref{eqn:proof-3.6(b)} can be made smaller than $\eps/3$ by choosing a sufficiently large $n$. Thus, concluding that 
\begin{equation}\label{eqn:derivative-2}
\lim _{n \to \infty} \frac{1}{n}\sum_{i=1}^n \frac{X_i-\hat q_n(\alpha_n u)}{\|X_i-\hat q_n(\alpha_n u)\|}\mathbbm{1}_{\{X_i\neq \hat q_n(\alpha_n u)\}}=-v, \,\,\,\,\,\,\,\,\text{a.s.}
\end{equation}
together with \eqref{eqn:derivative-1}, contradicts the assumption about the existence of $v\neq u$, thereby completing the proof of part $(ii)$ of Theorem \ref{thm:empirical-geom-quant-asymp1}. \hfill \qed


\subsection{Proof of Theorem~\ref{thm:empirical-geom-quant-asymp2}}\label{ss:theo3.3}

We now state (and prove) the following auxilliary result that provides an upperbound on the rate of growth of the sample geometric quantiles, and this will form a necessary part in the proof of Theorem~\ref{thm:empirical-geom-quant-asymp2}.
\begin{proposition}\label{thm:rate of quantile} 
Let $\{X_n\}_{n\ge 1}$ be an i.i.d. sample drawn from distribution $\pr$ on $\real^d$, whose support is not contained in any one dimensional affine subspace of $\real^d$. Let $\{\delta_n\}_{n\ge 1}$, $\{k_n\}_{n\ge 1}$ and $\{\alpha_n\}_{n\ge 1}$ be sequences of real numbers satisfying the following conditions:
\begin{itemize}
\item $\delta_n$ be such that $\sum_{n=1}^{\infty}\exp{(-n \delta_n^2)}< \infty$
\item $k_n$ be such that $\pr(\|X_i\|>k_n)\leq \delta_n$
\item $\alpha_n \in (0,1)$ be such that $(1-4\delta_n-\alpha_n)>0$ and $\displaystyle \alpha_n \underset{n\to\infty}{\to} 1$.
\end{itemize}
Then, we have
$$\|\hat q_n (\alpha_n u)\|\leq (M_n +2)k_n \,\,\,\,\,\text{a.s.}$$ 
for any unit vector $u$ and any $M_n$ satisfying $\displaystyle M_n>\frac{\alpha_n+2 \delta_n}{1-4\delta_n-\alpha_n}$.
\end{proposition}
\noindent{\bf Proof of Proposition~\ref{thm:rate of quantile}:} Let $\{k_n\}_{n\ge 1}$ and $\{\delta_n\}_{n\ge 1}$ be sequences of non--negative real numbers such that $k_n \uparrow \infty$ and $\delta_n\downarrow 0$ as $n\to\infty$, such that
\begin{equation} \label{eqn:decay rate}
\pr(\|X\|> k_n)\leq\delta_n,
\end{equation} 
for all $n\geq 1$ (exact form of $\delta_n$ will be chosen later). Now invoking Theorem A as stated on p.201 of \cite{Serfling1980} (used in \cite{Chaudhury1992}), Fact 5.1), 
we introduce
$$
p_n(X_i)=\textbf{1}_{\{\|X_i\|>k_n\}}-\pr({\|X_i\|>k_n}).
$$
Clearly, $p_n$ is symmetric and $|p_n(.)|\leq 2$, for all $n$. Also, $\E[p_n(X_i)]=0$ and 
\begin{equation*}
\sigma^2 = \text{var}(p_n(X_i)) = \pr\left(\|X_i\|>k_n\right)\, \pr\left(\|X_i\|\leq k_n\right) <1
\end{equation*}
Therefore, applying \cite[Theorem A,p.201]{Serfling1980} for large $n$ with $m=1$ and $b=2$ gives:
\begin{eqnarray*}
\pr\left(\left|\frac{1}{n}\sum_{i=1}^n\textbf{1}_{\{\|X_i\|>k_n\}}-P(\|X_i\|>k_n)\right|\geq \delta_n\right)
&\leq & 2\exp\left(-\frac{n\delta_n^2}{2\sigma^2+\frac{4}{3}\,\delta_n}\right) \,\,
\leq  2\exp\bigg(-n\frac{\delta_n^2}{4}\bigg).
\end{eqnarray*}
Note that the first assumption on $\delta_n$ in the statement of Proposition~\ref{thm:rate of quantile} ensures that \\
$\displaystyle
\sum_{i=1}^n \exp \left(-n\frac{\delta_n^2}{4}\right)<\infty$. Hence, by Borel--Cantelli, we conclude that
$$ \Big|\frac{1}{n}\sum_{i=1}^n\textbf{1}_{\{\|X_i\|>k_n\}}- \pr(\|X_i\|>k_n)\Big| < \delta_n, \,\,\,\text{ a.s.} $$ 
for all but finitely many $n$. Therefore, using the assumption on $\{k_n\}$ (in the same proposition), we have
\begin{equation}
\label{eqn:emp prob est1}
\frac{1}{n}\sum_{i=1}^n\mathbbm{1}_{\{\|X_i\|>k_n\}} \leq \pr(\|X_i\|>k_n)+\delta_n = 2 \delta_n \,\,\,\,\text{a.s.}
\end{equation}
for all but finite many $n$, meaning,
\begin{equation}
\label{eqn:emp prob est2}
\frac{1}{n}\sum_{i=1}^n\mathbbm{1}_{\{\|X_i\|\leq k_n\}}\geq (1-2\delta_n) \,\,\,\,\text{a.s.}
\end{equation}
for all but finitely many $n$. 

Having obtained the preliminary estimates, we are now set to estimate $\|\hat q_n(\alpha_n u)\|$. Recall that,
$$\hat q_n(\alpha_n u) = \argmin_q \frac{1}{n}\sum_{i=1}^n \left(\|X_i-q\|-\alpha_n \langle u,q\rangle -\|X_i\| \right)$$
So, let us study the objective functional, which we can write as
\begin{eqnarray*}
&& \frac{1}{n}\sum_{i=1}^n \left( \|X_i-q\|-\alpha_n \langle u,q\rangle-\|X_i\|\right)\\
& = & \frac{1}{n}\sum_{i=1}^n \left(\| X_i-q\| - \|X_i\|\right)\mathbbm{1}_{\{\|X_i\|\leq k_n\}}+\frac{1}{n}\sum_{i=1}^n \left(\|X_i-q\| -\|X_i\|\right)\mathbbm{1}_{\left\{\|X_i\|> k_n\right\}}-\alpha_n \langle u,q\rangle\\
& \geq & \frac{1}{n}\sum_{i=1}^n \left( \|q\|-\|2X_i\|\right)\mathbbm{1}_{\{\|X_i\|\leq k_n\}}+\frac{1}{n}\sum_{i=1}^n -\|q\|\mathbbm{1}_{\{\|X_i\|> k_n\}}-\alpha_n \|q\|,
\end{eqnarray*}
where we used Cauchy--Schwartz inequality for the last term, and the triangle inequality in two different ways in the last inequality.

Now, setting $\|q\|>(M_n+2)k_n$ for some $M_n>1$, whose exact form will be chosen later, we have, whenever $\|X_i\|\leq k_n$, $\left( \|q\|-2\|X_i\|\right)>\|q\|-2k_n> M_n k_n>\frac{M_n}{M_n+1}|q|$. Therefore, the above objective functional can be further reduced to
$$
\frac{1}{n}\sum_{i=1}^n \left(\|X_i-q\|-\alpha_n  \langle u,q\rangle-\|X_i\|\right)
\ge \frac{M_n}{M_n+2} \|q\| \frac{1}{n}\sum_{i=1}^n \mathbbm{1}_{\{\|X_i\|\leq k_n\}}- \|q\| \frac{1}{n}\sum_{i=1}^n \mathbbm{1}_{\{\|X_i\|> k_n\}}-\alpha_n \|q\|.
$$
Subsequently, using Equations~\eqref{eqn:emp prob est1} and \eqref{eqn:emp prob est2} in the previous inequality provides
\begin{eqnarray}\label{eqn:qn-minimiser}
\frac{1}{n}\sum_{i=1}^n \left(\|X_i-q\|-\alpha_n  \langle u,q\rangle-\|X_i|\right)
& \geq & \frac{M_n}{M_n+2} \|q\| (1-2\delta_n)-\|q\| 2 \delta_n - \|q\| \alpha_n \nonumber\\
& \geq & \|q\|\, \left[\frac{M_n}{M_n+2} (1-2\delta_n)-(\alpha_n+2 \delta_n)\right]\nonumber\\
& > & 0,  \,\,\,\,\,\,\,\,   \text{ if } M_n>\frac{\alpha_n+2 \delta_n}{1-4\delta_n-\alpha_n}.
\end{eqnarray}
As observed in \eqref{eqn:qn-minimiser}, the objective functional stays positive whenever $\|q\| > (M_n+2)k_n$. On the other hand, the objective functional equals $0$ for $q = 0$. Therefore, $\|\hat q_n(\alpha_n u)\|\leq (M_n+2)k_n$, a.s. \hfill \qed

\paragraph{Proof of Theorem~\ref{thm:empirical-geom-quant-asymp2} $\mathbf{(i)}$ -}
Let us express
$$\hat{q}_n=\argmin_{q}\frac{1}{n}\sum_{i=1}^n \phi (u,X_i-q)$$ 
with $\phi (u,x)= \|x\|+\langle u,x\rangle$. Clearly, for $x\neq 0$,
\begin{equation}\label{eqn:phi-derivative-I}
\lim_{t\rightarrow 0^+} \frac{\phi(u,x+th) - \phi(u,x)}{t} = \langle u,h\rangle+\left\langle \frac{x}{\|x\|},h\right\rangle.
\end{equation}
Similarly, for $x=0$, 
\begin{equation}\label{eqn:phi-derivative-II}
\lim_{t \to 0^+} \frac{\phi (u,th)-\phi (u,0)}t
= \lim_{t \to 0^+}\frac{\|th\|+t\langle u,h\rangle}t
=  \|h\|+\langle u,h\rangle.
\end{equation}
Observe that the function $\frac{1}{n}\sum_{i=1}^n \phi(u,X_i-q)$ is strictly convex in $q$, therefore if $\hat{q}_n(u)$ is a solution of this optimisation problem, then, for any $h\in\real^d$,
$$\lim_{t\rightarrow{0^+}} t^{-1}\left[\frac{1}{n}\sum_{i=1}^n \phi (u,X_i-\hat q _n(u)+ ht )-\frac{1}{n}\sum_{i=1}^n\phi (u,X_i-\hat q_n(u))\right]\geq 0,$$
which implies, using \eqref{eqn:phi-derivative-I} when $x\neq 0$, and \eqref{eqn:phi-derivative-II} otherwise,
\begin{eqnarray}
\label{eqn:inequality}
&&\frac{1}{n}\sum_{\substack{{i=1}\\X_i\neq \hat q_n}}^n\left\langle\frac{X_i-\hat q_n(u)}{\|X_i-\hat q_n(u)\|}+u,h\right\rangle+\frac{1}{n}\sum_{\substack{{i=1}\\X_i= \hat q_n}}^n \left(\|h\|+\langle u,h\rangle \right)\geq 0 \nonumber\\
& \implies & \langle u,h \rangle + \frac{1}{n}\sum_{\substack{{i=1}\\X_i\neq \hat q_n}}^n \left\langle\frac{X_i-\hat q_n(u)}{\|X_i-\hat q_n(u)\|},h\right\rangle \geq - \frac{1}{n}\sum_{\substack{{i=1}\\X_i= \hat q_n}}^n \|h\|.
\end{eqnarray}

By replacing $h$ with $-h$ in  \eqref{eqn:inequality}, we can conclude that
\begin{equation}
\label{eqn1:inequality}
\left| \langle u,h \rangle + \frac{1}{n}\sum_{\substack{{i=1}\\X_i\neq \hat q_n}}^n \left\langle\frac{X_i-\hat q_n(u)}{\|X_i-\hat q_n(u)\|},h\right\rangle\right| \leq  \frac{1}{n}\sum_{\substack{{i=1}\\X_i= \hat q_n}}^n \|h\|.
\end{equation}
Let us consider the expression $\left\langle\alpha_n u-\frac{\hat q_n(\alpha_n u)}{\|\hat q_n(\alpha_n u)\|},v\right\rangle$. Using \eqref{eqn1:inequality} to estimate $\langle \alpha_nu,v\rangle$, we obtain
$$\left\langle\alpha_n u-\frac{\hat q_n(\alpha_n u)}{\|\hat q_n(\alpha_n u)\|},v\right\rangle\\
\leq  -\frac{1}{n}\sum_{\substack{{i=1}\\X_i\neq \hat q_n}}^n\left\langle\frac{X_i-\hat q_n(\alpha_n u)}{\|X_i-\hat q_n(\alpha_n u)\|},v\right\rangle +\frac{1}{n}\sum_{\substack{{i=1}\\X_i= \hat q_n}}^n|v|-\left\langle\frac{\hat q_n( \alpha_n u)}{\|\hat q_n( \alpha_n u)\|},v\right\rangle .$$
Adding and subtracting the sample average of $\left\langle\frac{X_i}{\|\hat q_n(\alpha_n u)\|},v\right\rangle$, rearranging the terms, and denoting 
$$W(x,q)=\left[\left|\frac{x}{\|q\|}-\frac{q}{\|q\|}
\right|^{-1}-1\right]\left(\frac{x}{\|q\|}-\frac{q}{\|q\|}\right),$$ 
we rewrite the above inequality as
\begin{eqnarray*}
\left\langle\alpha_n u-\frac{\hat q_n(\alpha_n u)}{\|\hat q_n(\alpha_n u)\|},v\right\rangle
&=&  -\frac{1}{n}\sum_{\substack{{i=1}\\X_i\neq \hat q_n}}^n\left\langle W(X_i,\hat q_n(\alpha_n u)),v\right\rangle+ \frac{1}{n}\sum_{\substack{{i=1}\\X_i= \hat q_n}}^n \|v\| \\ 
&& \mbox{} +\frac{1}{n}\sum_{\substack{{i=1}\\X_i= \hat q_n}}^n\left\langle\frac{X_i}{\|\hat q_n(\alpha_n u)\|}-\frac{\hat q_n( \alpha_n u)}{\|\hat q_n( \alpha_n u)\|},v\right\rangle - \frac{1}{n}\sum_{i=1}^n \left\langle\frac{X_i}{\|\hat q_n(\alpha_n u)\|},v\right\rangle.
\end{eqnarray*}
Multiplying the above with $\|\hat q_n(\alpha_n u)\|$, we further have
\begin{eqnarray*}
&& \|\hat q_n(\alpha_n u)\|\left\langle\alpha_n u-\frac{\hat q_n(\alpha_n u)}{\|\hat q_n(\alpha_n u)\|},v\right\rangle + \frac{1}{n}\sum_{\substack{{i=1}\\x_i\neq \hat q_n}}^n \|\hat q_n(\alpha_n u)\|\left\langle  W(X_i,\hat q_n(\alpha_n u)),v\right\rangle + \frac{1}{n}\sum_{i=1}^n\left\langle X_i,v\right\rangle \\
&&\leq \, \|\hat q_n(\alpha_n u)\|\frac{\|v\|}{n} \# \{i: X_i=\hat q_n(\alpha_n u)\}.
\end{eqnarray*}
Similarly, by replacing $v$ with $-v$ in the above inequality,  we have
\begin{eqnarray*}
&& \left|\|\hat q_n(\alpha_n u)\|\left\langle\alpha_n u-\frac{\hat q_n(\alpha_n u)}{\|\hat q_n(\alpha_n u)\|},v\right\rangle + \frac{1}{n}\sum_{\substack{{i=1}\\x_i\neq \hat q_n}}^n \|\hat q_n(\alpha_n u)\|\left\langle  W(X_i,\hat q_n(\alpha_n u)),v\right\rangle + \frac{1}{n}\sum_{i=1}^n\left\langle X_i,v\right\rangle \right|\\
&&\leq \,\|\hat q_n(\alpha_n u)\|\frac{\|v\|}{n} \# \{i: X_i=\hat q_n(\alpha_n u)\}.
\end{eqnarray*}
The proof of part $(i)$ follows trivially given the following claims:

\textbf{Claim 1:} 
\begin{equation}\label{eqn:claim-1-empirical-geom-quant-asymp2}
\lim_{n \to \infty} \frac{1}{n}\sum_{\substack{{i=1}\\ X_i\neq \hat q_n}}^n \|\hat q_n(\alpha_n u)\|\left\langle  W(X_i,\hat q_n(\alpha_n u)),v\right\rangle = -\langle u,v\rangle E\langle X, u \rangle,\,\,\,\,\text{a.s.}
\end{equation}
\textbf{Claim 2:}
\begin{equation}\label{eqn:claim-2-empirical-geom-quant-asymp2}
    \lim_{n \to \infty}\frac{\|\hat q_n(\alpha_n u)\|}{n} = 0,\,\,\,\text{a.s.}
\end{equation}
	
{\it Proof of Claim 1.} Set  
$$F_{n,i}=\|\hat q_n(\alpha_n u)\|\left\langle W(X_i, \hat q_n(\alpha_n u)),v\right\rangle \mathbbm{1}_{\{X_i \neq \hat q_n(\alpha_n u)\}}+\left \langle u,v\right \rangle \langle X_i, u \rangle
$$ and 
$$F_n=\frac{1}{n}\sum_{i=1}^n F_{n,i}.$$ 
Clearly, it suffices to show that $|F_n| \to 0$ almost surely, as $n \to \infty$, to obtain \eqref{eqn:claim-1-empirical-geom-quant-asymp2}.

We shall break the sum in two parts 
$$
F_n=\frac{1}{n}\sum_{i=1}^n F_{n,i}\one_{\{\|X_i\|\leq M\}}+\frac{1}{n}\sum_{i=1}^nF_{n,i}\one_{\{\|X_i\|\geq M\}},
$$ 
where $M>0$ will be chosen later accordingly, and show that each of the two terms converges to $0$ almost surely. 
For the second term, we recall the following estimate (see the proof of Lemma 6.2 in \cite{Girard2017}, p.134),
\begin{equation*}
\big|\|\hat q_n(\alpha_n u)\|\left\langle W(X_i, \hat q_n(\alpha_n u)),v\right\rangle\big|\mathbbm{1}_{\{X_i \neq \hat q_n(\alpha_n u)\}} \leq 3\|v\|\|X_i\|,
\end{equation*}
which implies $$|F_{n,i}| \leq 4 \|v\|\|X_i\|.$$ 
Thus, we have
\begin{equation}\label{eqn:claim-1-empirical-geom-quant-asymp2-truncate}
\frac{1}{n}\sum_{i=1}^nF_{n,i}\mathbbm{1}_{\{\|X_i\|\geq M\}}\leq 4\|v\|\frac{1}{n}\sum_{i=1}^n \|X_i\| \mathbbm{1}_{\{\|X_i\|\geq M\}} \overset{a.s.}{\underset{n\to\infty}{\longrightarrow}} 4\,\|v\|\E\left[\|X\|\mathbbm{1}_{\{\|X\|\geq M\}}\right].
\end{equation}
Let $M$ be large enough enough so that \, $\displaystyle 4\|v\|\E\left[\|X\|\mathbbm{1}_{\|X\|\geq M}\right]<\frac{\eps}{2}$.

Next, considering the first term, $F_{n,i}\mathbbm{1}_{\{\|X_i\|\leq M\}}$, observe that, whenever $X_i \neq \hat q_n(\alpha_n u)$, $W(X_i, \hat q_n(\alpha_n u))$ can be rewritten as
\begin{eqnarray*}
&& \|\hat q_n(\alpha_n u)\|\, W(X_i, \hat q_n(\alpha_n u))\\
&=& \|\hat q_n(\alpha_n u)\|\,\left\| \frac{X_i}{\|\hat q_n(\alpha_n u)\|}-\frac{\hat q_n(\alpha_n u)}{\|\hat q_n(\alpha_n u)\|}\right\|^{-1}\left[1- \left\| \frac{X_i}{\|\hat q_n(\alpha_n u)\|}-\frac{\hat q_n(\alpha_n u)}{\|\hat q_n(\alpha_n u)\|}\right\|^2\right]\\
&& \mbox{} \times  \left[1+\left\| \frac{X_i}{\|\hat q_n(\alpha_n u)\|}-\frac{\hat q_n(\alpha_n u)}{\|\hat q_n(\alpha_n u)\|}\right\|\right]^{-1} \left(\frac{X_i}{\|\hat q_n(\alpha_n u)\|}-\frac{\hat q_n(\alpha_n u)}{\|\hat q_n(\alpha_n u)\|}\right)\\
&=& \left\| \frac{X_i}{\|\hat q_n(\alpha_n u)\|}-\frac{\hat q_n(\alpha_n u)}{\|\hat q_n(\alpha_n u)\|}\right\|^{-1} \left[-\frac{\|X_i\|^2}{\|\hat q_n(\alpha_n u)\|}+2\left\langle X_i,\frac{\hat q_n(\alpha_n u)}{\|\hat q_n(\alpha_n u)\|}  \right\rangle\right] \\
&& \mbox{} \times \left[1+\left\| \frac{X_i}{\|\hat q_n(\alpha_n u)\|}-\frac{\hat q_n(\alpha_n u)}{\|\hat q_n(\alpha_n u)\|}\right\|\right]^{-1} \left(\frac{X_i}{\|\hat q_n(\alpha_n u)\|}-\frac{\hat q_n(\alpha_n u)}{\|\hat q_n(\alpha_n u)\|}\right).
\end{eqnarray*}
Since $\displaystyle \|\hat q_n(\alpha_n u)\| \overset{a.s.}{\underset{n\to\infty}{\longrightarrow}} \infty$ and $\frac{\hat q_n(\alpha_n u)}{\|\hat q_n(\alpha_n u)\|} \overset{a.s.}{\underset{n\to\infty}{\longrightarrow}} u$ almost surely, we thus have
\begin{equation}
\Big|\|\hat q_n(\alpha_n u)\|\left\langle W(X_i, \hat q_n(\alpha_n u)),v\right\rangle 
+ \left \langle u,v\right \rangle \langle X_i, u \rangle \Big| \mathbbm{1}_{\{X_i \neq \hat q_n(\alpha_n u),\,\|X_i\|\le M\}}
\overset{a.s.}{\underset{n\to\infty}{\longrightarrow}} 0.
\end{equation}
Observe that the above convergence is uniform in $i\leq n$. Therefore, for any $\eps>0$ and $\omega\in\{\|X_i\|\leq M\}$, there exists a random number $N_{\eps}(\omega)<\infty$ such that 
$$
\frac{1}{n}\sum_{i=1}^n F_{n,i}(\omega)\,\one_{\{\|X_i\|\leq M\}} < \frac{\eps}{2} \,\,\,\,\,\forall \,\,\, n \geq N_{\eps}(\omega).
$$ 
The claim follows by collating the estimates obtained above and in \eqref{eqn:claim-1-empirical-geom-quant-asymp2-truncate}.

{\it Proof of Claim 2.} We shall invoke Proposition~\ref{thm:rate of quantile} to prove this claim. Let us set $\displaystyle\delta_n=\frac{(1-\alpha_n)}{5}$, implying $1-4\delta_n-\alpha_n>0$. Since $\E\|X_i\|<\infty$, a simple application of Markov's inequality gives us 
$$\pr\left(\|X_i\|>\frac{5\E\|X_i\|}{(1-\alpha_n)}\right)\leq\frac{(1-\alpha_n)}{5}=\delta_n.$$ 
Now the assumption $\displaystyle\sum_{n=1}^{\infty} \exp(-n(1-\alpha_n)^2)<\infty$ guarantees that we can apply Proposition~\ref{thm:rate of quantile}. Therefore, $\displaystyle\|\hat q_n (\alpha_n u)\|\leq (M_n +2)\frac{5\E\|X_i\|}{(1-\alpha_n)} $ where $M_n$ is any number greater than $\displaystyle\frac{\alpha_n+2\delta_n}{1-4\delta_n-\alpha_n}$. Let $\displaystyle M_n=\frac{5}{1-\alpha_n}$, then $\|\hat q_n (\alpha_n u)\|=\textrm{O}_{\pr}(\frac{1}{(1-\alpha_n)^2})$. Since $\displaystyle\sum_{n=1}^{\infty} \exp(-n(1-\alpha_n)^2)<\infty$, we have $\displaystyle\frac{\|\hat q_n(\alpha_n u)\|}{n}\to 0$ almost surely, as $n\to \infty$, which completes the proof of Claim 2, and hence concludes the proof of part $(i)$ of Theorem \ref{thm:empirical-geom-quant-asymp2}. \hfill \qed

\paragraph{Proof of Theorem~\ref{thm:empirical-geom-quant-asymp2} $\mathbf{(ii)}$ -}
Substituting $u$ by $\alpha_n u$ and $h$ by $\frac{\hat q_n(\alpha_n u)}{\|\hat q_n(\alpha_n u)\|}$ in \eqref{eqn1:inequality}, we obtain
\begin{equation}
\label{eqn:sample quantile characterisation2}
\left| \left\langle \alpha_n u, \frac{\hat q_n(\alpha_n u)}{\|\hat q_n(\alpha_n u)\|} \right\rangle + \frac{1}{n}\sum_{\substack{{i=1}\\X_i\neq \hat q_n}}^n \left\langle\frac{X_i-\hat q_n(\alpha_n u)}{\|X_i-\hat q_n(\alpha_n u)\|}, \frac{\hat q_n(\alpha_n u)}{\|\hat q_n(\alpha_n u)\|}\right\rangle\right| \leq  \frac{1}{n}\# \{i: X_i= \hat q_n\}.
\end{equation}
Then, we have
\begin{eqnarray*}
\left\langle\alpha_n u-\frac{\hat q_n(\alpha_n u)}{\|\hat q_n(\alpha_n u)\|},\frac{\hat q_n(\alpha_n u)}{\|\hat q_n(\alpha_n u)\|}\right\rangle 
& = & \left\langle\alpha_n u,\frac{\hat q_n(\alpha_n u)}{\|\hat q_n(\alpha_n u)\|}\right\rangle-1\\
& \leq & -\frac{1}{n}\sum_{\substack{{i=1}\\X_i\neq \hat q_n}}^n\left\langle\frac{X_i-\hat q_n(\alpha_n u)}{\|X_i-\hat q_n(\alpha_n u)\|}, \frac{\hat q_n(\alpha_n u)}{\|\hat q_n(\alpha_n u)\|} \right\rangle -\frac{1}{n}\sum_{\substack{{i=1}\\X_i\neq  \hat q_n}}^n 1\\
& = & -\frac{1}{n}\sum_{\substack{{i=1}\\X_i\neq \hat q_n}}^n \left[ \left\langle\frac{X_i-\hat q_n(\alpha_n u)}{\|X_i-\hat q_n(\alpha_n u)\|},\frac{\hat q_n(\alpha_n u)}{\|\hat q_n(\alpha_n u)\|}\right\rangle+1\right].
\end{eqnarray*}
Similarly,
\begin{eqnarray*}
\left\langle\alpha_n u-\frac{\hat q_n(\alpha_n u)}{\|\hat q_n(\alpha_n u)\|},\frac{\hat q_n(\alpha_n u)}{\|\hat q_n(\alpha_n u)\|}\right\rangle
& \ge & -\frac{1}{n}\sum_{\substack{{i=1}\\X_i\neq \hat q_n}}^n\left[\left\langle\frac{X_i-\hat q_n(\alpha_n u)}{\|X_i-\hat q_n(\alpha_n u)\|},\frac{\hat q_n(\alpha_n u)}{\|\hat q_n(\alpha_n u)\|}\right\rangle+1\right]-\frac{2}{n}\sum_{\substack{{i=1}\\X_i=  \hat q_n}}^n 1.
\end{eqnarray*}
Combining the two estimates, and setting 
\begin{equation}\label{eqn:z-n-i}
Z_{n,i}\stackrel{\Delta}{=}\left\langle\frac{X_i-\hat q_n(\alpha_n u)}{\|X_i-\hat q_n(\alpha_n u)\|},\frac{\hat q_n(\alpha_n u)}{\|\hat q_n(\alpha_n u)\|}\right\rangle+1, 
\end{equation}
we conclude
\begin{equation}
0\le \left|\left\langle\alpha_n u-\frac{\hat q_n(\alpha_n u)}{\|\hat q_n(\alpha_n u)\|},\frac{\hat q_n(\alpha_n u)}{\|\hat q_n(\alpha_n u)\|}\right\rangle+\frac{1}{n}\sum_{\substack{{i=1}\\X_i\neq \hat q_n}}^n Z_{n,i} \right| \leq \frac{2}{n} \,\,\,\text{ a.s.}
\end{equation}
since $\#\{i:X_i=\hat q_n(\alpha_n u)\}\leq 1$ almost surely for continuous random variables $X_i$.

Equivalently, almost surely, we have
$$\| \hat q_n(\alpha_n u)\|^2 \left|\left\langle\alpha_n u-\frac{\hat q_n(\alpha_n u)}{\|\hat q_n(\alpha_n u)\|},\frac{\hat q_n(\alpha_n u)}{\|\hat q_n(\alpha_n u)\|}\right\rangle+\frac{1}{n}\sum_{\substack{{i=1}\\X_i\neq \hat q_n}}^n Z_{n,i}\right| \leq 2 \frac{\| \hat q_n(\alpha_n u)\|^2}{n}.$$ 
To finalize the proof of the theorem, we prove the following two claims:

\textbf{Claim 3:} $$\frac{\|\hat q_n(\alpha_n u)\|^2}{n} \overset{a.s.}{\underset{n\to\infty}{\longrightarrow}} \,0.$$

\textbf{Claim 4:} 
\begin{equation*}
    \frac{\|\hat q_n(\alpha_n u)\|^2}{n}\sum_{\substack{{i=1}\\X_i\neq \hat q_n}}^n Z_{n,i} \;\overset{a.s.}{\underset{n\to\infty}{\longrightarrow}} \;\half\E\|X-\langle X,u\rangle u\|^2.
\end{equation*}
{\it Proof of Claim 3:} We will use Proposition~\ref{thm:rate of quantile} to prove this claim. Setting $\displaystyle \delta_n=\frac{1-\alpha_n}{8}$, we have $\displaystyle 1-4\delta_n-\alpha_n=\frac{1-\alpha_n}{2}>0$. Also assume $\displaystyle k_n=\sqrt{\frac{8\E(\|X_i\|^2)}{1-\alpha_n}}$. Then, by Chebyshev's inequality, we have
$$
\pr(\|X_i\|>k_n)\leq \frac{\E(\|X_i\|^2)}{k_n^2}=\frac{1-\alpha_n}{8}=\delta_n.
$$
Therefore, using Proposition~\ref{thm:rate of quantile}, we can write
$$
\|\hat q_n(\alpha_n u)\| \leq (M_n+2)k_n,
$$ 
where $M_n$ is any number greater than $\displaystyle \frac{\alpha_n+2 \delta_n}{1-4\delta_n-\alpha_n}$. Assuming $\displaystyle M_n=\frac{2}{1-\alpha_n}$, the previous inequality becomes 
$$\|\hat q_n(\alpha_n u)\|\leq \left(\frac{2}{1-\alpha_n}+2\right)\sqrt{\frac{8\E(\|X_i\|^2)}{1-\alpha_n}}.
$$ 
Recall that, by assumption, $\displaystyle n (1-\alpha_n)^3 \underset{n\to\infty}{\to} \infty$. Therefore, 
$$
\frac{\|\hat q_n(\alpha_n u)\|^2}{n}=\textrm{O}_{\pr} \left(\frac{1}{n(1-\alpha_n)^3}\right) \to 0 \,\,\,\,\,\text{as}\,\,\,\, n \to \infty,
$$ which establishes the identity stated in Claim 3.

{\it Proof of Claim 4.} We shall begin by breaking the primary expression in two parts:
\begin{eqnarray}\label{eqn:thm3.8-claim2-break}
&& \frac{1}{n}\sum_{\substack{i=1\\X_i \neq \hat q_n}}^n \left[\|\hat q_n(\alpha_n u)\|^2 Z_{n,i}-\half  \|X_i-\langle X_i,u\rangle u\|^2\right] = \nonumber\\
&&\frac{1}{n}\sum_{\substack{i=1\\X_i \neq \hat q_n}}^n \left[\|\hat q_n(\alpha_n u)\|^2 Z_{n,i}-\half  \|X_i-\langle X_i,u\rangle u\|^2\right]\mathbbm{1}_{\{\|X_i\|<M\}}\nonumber\\
&& \quad + \frac{1}{n}\sum_{\substack{i=1\\X_i \neq \hat q_n}}^n \left[\|\hat q_n(\alpha_n u)\|^2 Z_{n,i}-\half  \|X_i-\langle X_i,u\rangle u\|^2\right]\mathbbm{1}_{\{\|X_i\|\geq M\}},
\end{eqnarray}
where $M>0$ will be chosen appropriately to make both the terms arbitrarily small. 

Let us analyse the second term. First, notice that, from the proof of \cite[Lemma 6.3]{Girard2017}, $\|\hat q_n(\alpha_n u)\|^2 Z_{n i} \leq 2 \|X_i\|^2$. Therefore,
$$
\frac{1}{n}\sum_{i=1}^n \left|\|\hat q_n(\alpha_n u)\|^2 Z_{n,i}-\half  \|X_i-\langle X_i,u\rangle u\|^2\right|\mathbbm{1}_{\{|X_i|\geq M\}} \leq \frac{1}{n}\sum_{i=1}^n 3 \|X_i\|^2 \mathbbm{1}_{\{|X_i|\geq M\}}.
$$
Now, for a fixed $\epsilon>0$, there exists $N_M(\epsilon)$ such that, for all $n\ge N_M(\epsilon)$,
$$
\frac{1}{n}\sum_{i=1}^n \|X_i\|^2 \mathbbm{1}_{|X_i|\geq M} \le \frac{\eps}{4} +\E[|X_1|^2 \mathbbm{1}_{|X_1|>M}].
$$
Next, by integrability of $\|X_1\|^2$, there exists $M$ large enough such that
$$
\E[\|X_1\|^2 \mathbbm{1}_{\{|X_1|>M\}}] \le \frac{\epsilon}{4}.
$$
Collating the above results provides
\begin{equation}
\label{eqn:second part}
\frac{1}{n}\sum_{i=1}^n \left|\|\hat q_n(\alpha_n u)\|^2 Z_{n,i}-\half  \|X_i-\langle X_i,u\rangle u\|^2\right|\mathbbm{1}_{|X_i|\geq M} \le \frac{3\,\eps}4.
\end{equation}

We now analyse the first part of \eqref{eqn:thm3.8-claim2-break}, by rewriting $Z_{n,i}$, defined in \eqref{eqn:z-n-i}, as
\begin{eqnarray*}
Z_{n,i}
&=& \left\|\frac{X_i}{\|\hat q_n(\alpha_n u)\|}-\frac{\hat q_n(\alpha_n u)}{\|\hat q_n(\alpha_n u)\|}\right\|^{-1} \\
&& \mbox{} \times \left\{\left\|\frac{X_i}{\|\hat q_n(\alpha_n u)\|}-\frac{\hat q_n(\alpha_n u)}{\|\hat q_n(\alpha_n u)\|}\right\|-\left[1-\frac{1}{\|\hat q_n(\alpha_n u)\|}\left\langle X_i,\frac{\hat q_n(\alpha_n u)}{\|\hat q_n(\alpha_n u)\|} \right\rangle\right]\right\}\\
&=& A_{i,n} \times \frac{B_{i,n}}{\|\hat q_n(\alpha_n u)\|^2},
\end{eqnarray*}
where 
\begin{eqnarray*}
A_{i,n}&:=& \left\|\frac{X_i}{\|\hat q_n(\alpha_n u)\|}-\frac{\hat q_n(\alpha_n u)}{\|\hat q_n(\alpha_n u)\|}\right\|^{-1} \\
&& \mbox{} \times \left\{\left\|\frac{X_i}{\|\hat q_n(\alpha_n u)\|}-\frac{\hat q_n(\alpha_n u)}{\|\hat q_n(\alpha_n u)\|}\right\|+\left[1-\frac{1}{\|\hat q_n(\alpha_n u)\|}\left\langle X_i,\frac{\hat q_n(\alpha_n u)}{\|\hat q_n(\alpha_n u)\|} \right\rangle\right]\right\}^{-1}\\
& \overset{a.s.}{\underset{n\to\infty}{\longrightarrow}} & \frac12 
\end{eqnarray*}
and 
\begin{eqnarray*}
B_{i,n}&:=&\left\langle X_i-\left\langle X_i,\frac{\hat q_n(\alpha_n u)}{\|\hat q_n(\alpha_n u)\|}\right\rangle \frac{\hat q_n(\alpha_n u)}{\|\hat q_n(\alpha_n u)\|},X_i\right\rangle \\
&\overset{a.s.}{\underset{n\to\infty}{\longrightarrow}}& \langle X_i-\langle X_i,u \rangle u, X_i\rangle,
\end{eqnarray*} 
using Theorem \ref{thm:empirical-geom-quant-asymp1}.

Moreover, observe that when the $X_i$'s are bounded, these convergences are uniform for $1\le i\le n$, namely, if $\|X_i\|< M$, then, for $\eps>0$, there exists a random $N_M(\eps)$ such that, almost surely,
$$
\left|\left(\|\hat q_n(\alpha_n u)\|^2 Z_{n,i}-\half  \|X_i-\langle X_i,u\rangle u\|^2\right)\mathbbm{1}_{|X_i|< M}\right| < \frac{\eps}{4} , \quad\forall n \ge N_M(\eps), \,\forall i \leq n.
$$
Hence, for $n \geq N_M(\eps)$, we have
\begin{equation}
\label{eqn:first part}
\frac{1}{n}\sum_{\substack{i=1\\X_i \neq \hat q_n}}^n \left[\|\hat q_n(\alpha_n u)\|^2 Z_{n,i}-\half  \|X_i-\langle X_i,u\rangle u\|^2\right]\mathbbm{1}_{\{|X_i|<M\}}<\frac{\eps}{4}.
\end{equation}
Combining \eqref{eqn:second part} and \eqref{eqn:first part}, we get, for $n \geq \max\{ N_M(\eps),\bar N_M(\eps)\}$,
$$
\frac{1}{n}\sum_{\substack{i=1\\X_i \neq \hat q_n}}^n \left[\|\hat q_n(\alpha_n u)\|^2 Z_{n,i}-\half  \|X_i-\langle X_i,u\rangle u\|^2\right]< \frac{3\eps}{4}+\frac{\eps}{4}=\eps.
$$
Since $\eps$ is arbitrary, the claim is proved, which in turn concludes the proof of Theorem \ref{thm:empirical-geom-quant-asymp2}. \hfill \qed

\subsection{Proof of Theorem~\ref{thm:numerics}.}\label{ss:theo3.4} 

We shall begin with the proof of the first part, by expressing all the terms involved in terms of the orthonormal basis $\{u,w_1,...,w_{d-1}\}$ of $\real^d$. Let $b^{(n)}(\alpha_n)$ and $\{\beta_k^{(n)}(\alpha_n)\}_{k=1}^{d-1}$ be real numbers defined by 
$$
\frac{\hat q_n(\alpha_n u)}{\|\hat q_n(\alpha_n u)\|}=b^{(n)}(\alpha_n )u+\sum _{k-1}^{d-1}\beta_k^{(n)}(\alpha_n)w_k.
$$
Note that
\begin{equation}\label{eqn:orthonormality}
b^{(n)}(\alpha_n)^2 + \sum_{k=1}^{d-1}\beta_k^{(n)}(\alpha_n)^2=1.
\end{equation}
Using the orthonormal basis $\{u,w_1,...,w_{d-1}\}$, we rewrite the expression of interest as
\begin{eqnarray*}
&& \hat q_n(\alpha_n u) - \|\hat q_n(\alpha_n u)\|\,u - \frac1n \sum_{i=1}^n \left( X_i-\langle X_i,u\rangle\right) \\
&=& \|\hat q_n(\alpha_n u)\| \left(b^{(n)}(\alpha_n)-1\right) u + \sum_{k=1}^{d-1} \left(\|\hat q_n(\alpha_n u)\| \beta_k^{(n)}(\alpha_n)-\frac{1}{n}\sum_{i=1}^n \langle X_i,w_k\rangle \right) w_k.
\end{eqnarray*}
Recall that, as a consequence of Theorem \ref{thm:empirical-geom-quant-asymp1}, $$b^{(n)}(\alpha_n)  \overset{a.s.}{\underset{n\to\infty}{\longrightarrow}} 1, $$ 
which also means, using \eqref{eqn:orthonormality}, that 
\begin{equation}\label{eqn:quantile-other-projection}
    \beta^{(n)}_k(\alpha_n) \overset{a.s.}{\underset{n\to\infty}{\longrightarrow}} 0,\,\,\,\forall k=1,\ldots, (d-1).
\end{equation}
Additionally, we know from Theorem~\ref{thm:empirical-geom-quant-asymp2} $(i)$ that, for all $k \in 1,...,d-1$, 
\begin{eqnarray}
\label{eqn:projection of quantile}
\| \hat q_n(\alpha_n u)\|\left\langle\alpha_n u- \frac{\hat q_n(\alpha_n u)}{\|\hat q_n(\alpha_n u)\|}, w_k\right\rangle
&=& \mbox{} - \| \hat q_n(\alpha_n u)\| \beta_k^{(n)}(\alpha_n) \nonumber \\
&\overset{a.s.}{\underset{n\to\infty}{\longrightarrow}} & -\E\langle X-\langle X, u \rangle u,w_k \rangle.
\end{eqnarray}

Therefore, combining \eqref{eqn:quantile-other-projection} and \eqref{eqn:projection of quantile} together with \eqref{eqn:orthonormality}, we obtain
\begin{eqnarray}\label{eqn:quantile-projection-1}
\| \hat q_n(\alpha_n u)\|\left(1-b^{(n)}(\alpha_n)\right) &=& \| \hat q_n(\alpha_n u)\|\left(1-b^{(n)}(\alpha_n)^2\right)(1+b^{(n)}(\alpha_n)) \nonumber \nonumber \\
&=& \sum_{k=1}^{d-1}\| \hat q_n(\alpha_n u)\|\beta_k^{(n)}(\alpha_n)^2 \, \left(1+b^{(n)}(\alpha_n)\right) 
\overset{a.s.}{\underset{n\to\infty}{\longrightarrow}}  0.
\end{eqnarray}
Now, combining \eqref{eqn:projection of quantile}, \eqref{eqn:quantile-projection-1} and the strong law of large numbers, we can conclude that
$$
\hat q_n(\alpha_n u) - \|\hat q_n(\alpha_n u)\|\,u - \frac1n \sum_{i=1}^n \left( X_i-\langle X_i,u\rangle\right) \overset{a.s.}{\underset{n\to\infty}{\longrightarrow}} 0,
$$
which proves the first part of the theorem.

Let us move to the second part of the theorem for which we recall that $\Sigma$ denotes the covariance matrix corresponding to $\pr$.

Observe that, using a combination of \eqref{eqn:projection of quantile} and the continuous mapping theorem, we can write
\begin{equation}\label{eqn:beta-k-n-cgs}
\Big(\| \hat q_n(\alpha_n u)\| \beta_k^{(n)}(\alpha_n u)\Big)^2 \overset{a.s.}{\underset{n\to\infty}{\longrightarrow}} \Big(\E \big\langle X-\langle X, u \rangle u,w_k \big\rangle\Big)^2 = \left(\E\langle X,w_k\rangle\right)^2 ,
\end{equation}
for all $k=1,\ldots,d-1$. 
From Theorem \ref{thm:empirical-geom-quant-asymp2} $(ii)$ and Equation \ref{eqn:orthonormality}, we have,
\begin{equation}\label{eqn:b-n-cgs}
\| \hat q_n(\alpha_n u)\|^2[\alpha_n b^{(n)}(\alpha_n)-1] \overset{a.s.}{\underset{n\to\infty}{\longrightarrow}}  -\half \E \|X-\langle X,u\rangle u\|^2=-\half \sum_{k=1}^{d-1}\E \langle X,w_k\rangle^2.
\end{equation}

Combining \eqref{eqn:beta-k-n-cgs} and \eqref{eqn:b-n-cgs} provides
\begin{eqnarray}
\label{eqn:rate1}
\| \hat q_n(\alpha_n u)\|^2\left[1-\alpha_n b^{(n)}(\alpha_n)-\half\sum_{k=1}^{d-1}\beta^{(n)}_k(\alpha_n)^2\right] & \overset{a.s.}{\underset{n\to\infty}{\longrightarrow}} &\half \sum_{k=1}^{d-1}\Big(\E \langle X,w_k\rangle^2-\left[\E\langle X,w_k\rangle\right]^2\Big) \nonumber\\
&& = \half \sum_{k=1}^{d-1} \text{Var}\langle X,w_k\rangle.
\end{eqnarray}
Note that, from Theorem \ref{thm:empirical-geom-quant-asymp2}, 
$\displaystyle
\| \hat q_n(\alpha_n u)\|\left\langle\alpha_n u-\frac{\hat q_n(\alpha_n u)}{\|\hat q_n(\alpha_n u)\|},u\right\rangle \overset{a.s.}{\underset{n\to\infty}{\longrightarrow}} 0$,
which then implies
\begin{equation}\label{eqn:rate2}
\| \hat q_n(\alpha_n u)\|^2(\alpha_n -b^{(n)}(\alpha_n))^2 \overset{a.s.}{\underset{n\to\infty}{\longrightarrow}} 0.
\end{equation}
Finally, summing \eqref{eqn:rate1} and \eqref{eqn:rate2} provides
$$ 
\| \hat q_n(\alpha_n u)\|^2\left[1-\alpha_n b^{(n)}(\alpha_n)-\half\sum_{k=1}^{d-1}\beta^{(n)}_k(\alpha_n)^2+(\alpha_n -b^{(n)}(\alpha_n))^2\right] \overset{a.s.}{\underset{n\to\infty}{\longrightarrow}} \half \sum_{k=1}^{d-1} \text{Var}\langle X,w_k\rangle.
$$
Finally, observe that
$$
\| \hat q_n(\alpha_n u)\|^2\left(1-\alpha_n\right) - \| \hat q_n(\alpha_n u)\|^2\left[1-\alpha_n b^{(n)}(\alpha_n)-\half\sum_{k=1}^{d-1}\beta^{(n)}_k(\alpha_n)^2+(\alpha_n -b^{(n)}(\alpha_n))^2\right] = o_{\pr}(1),
$$
which proves the result since $\text{Var}\langle X,u\rangle+\sum_{k=1}^{d-1} \text{Var}\langle X,w_k\rangle=\text{tr}\, (\Sigma)$. \hfill \qed

\subsection{Proof of Theorem~\ref{thm:empHDdecay-vs-popHD}}
\label{ss:proofTheoEmpHDdecay}

The proof is articulated in three steps, the first one based on the observation made in Lemma~\ref{prop:hdepth-sup-emp-ratio-gl}, the second one adapting the latter result to the framework given in Theorem~\ref{thm:empHDdecay-vs-popHD}, the third and last step using Theorem 5.1 in \cite{Alexander1987}.

{\it Step 1.} First notice that Lemma~\ref{prop:hdepth-sup-emp-ratio-gl} follows directly from the observations that: 
\begin{itemize}
    \item[] If $\displaystyle HD(t_nx,\pr_n) \ge HD(t_nx,\pr)$, then for $\displaystyle\widetilde{H}_{t_nx}$ such that  $\displaystyle\pr(\widetilde{H}_{t_nx}) = HD(t_nx,\pr)$, we can write:
\begin{equation}\label{eq:lemPn>P}
\left|\dfrac{HD(t_n x, \pr_n)}{HD(t_n x, \pr)}-1\right| = \dfrac{HD(t_n x, \pr_n)}{HD(t_n x, \pr)}-1  \le \left(\frac{\pr_n(\widetilde{H}_{t_nx})}{\pr(\widetilde{H}_{t_nx})} - 1 \right) \le \sup_{H \in \calH_{t_nx}} \left|\frac{\pr_n(H)}{\pr(H)} - 1 \right|.
\end{equation}
\item[] Similarly, if $HD(t_nx,\pr_n) \le HD(t_nx,\pr)$, then for $H^*_{t_nx}$ s.t. $ \pr_n(H^*_{t_nx})= HD(t_nx,\pr_n)$, we have
$$\left|\dfrac{HD(t_n x, \pr_n)}{HD(t_n x, \pr)}-1\right| =  1- \dfrac{HD(t_n x, \pr_n)}{HD(t_n x, \pr)} \le \left(1 - \frac{\pr_n(H^*_{t_nx})}{\pr(H^*_{t_nx})}  \right) \le \sup_{H \in \calH_{nx}} \left|\frac{\pr_n(H)}{\pr(H)} - 1 \right|.$$
\end{itemize}

{\it Step 2.} In view of Condition $(C2)$ of Theorem~\ref{thm:empHDdecay-vs-popHD}, we need to consider a subset of $\calH$, hence to adapt Lemma~\ref{prop:hdepth-sup-emp-ratio-gl} to this smaller class, as follows:
\begin{proposition}\label{prop:hdepth-sup-emp-ratio}
    Under Condition $(C2)$ of Theorem~\ref{thm:empHDdecay-vs-popHD}, we have
    \begin{equation}
    \left| \frac{HD(t_nx,\pr_n)}{HD(t_nx,\pr)}-1\right| \le \sup\left\{ \left| \frac{\pr_n(H)}{\pr(H)}-1\right|:\,\,\, H\in \calH_{t_nx},\,\,\pr(H) \ge \gamma_n\right\}.
    \end{equation}
\end{proposition} 
{\bf Proof of Proposition~\ref{prop:hdepth-sup-emp-ratio}.} The arguments for the proof of Proposition~\ref{prop:hdepth-sup-emp-ratio} are identical to those used to prove Lemma~\ref{prop:hdepth-sup-emp-ratio-gl}, but taking into account Condition $(C2)$ of Theorem~\ref{thm:empHDdecay-vs-popHD}. Let us consider the first case, when $HD(t_nx,\pr_n) \ge HD(t_nx,\pr)$, for which we have \eqref{eq:lemPn>P}. Then, Condition $(C2)$ implies $HD(t_nx,\pr) = \pr(\widetilde{H}_{t_n x}) \ge \gamma_n$, leading to
$$
\left| \frac{HD(t_nx,\pr_n)}{HD(t_nx,\pr)}-1\right| \le \sup\left\{ \left| \frac{\pr_n(H)}{\pr(H)}-1\right|:\,\,\, H\in \calH_{t_nx},\,\,\pr(H) \ge \gamma_n\right\}.
$$
Next, let us consider the case $HD(t_nx,\pr_n) \le HD(t_nx,\pr)$, where we have
$$\left|\dfrac{HD(t_n x, \pr_n)}{HD(t_n x, \pr)}-1\right| =  1- \dfrac{HD(t_n x, \pr_n)}{HD(t_n x, \pr)} = 1 - \frac{\pr_n(H^*_{t_nx})}{\pr(\widetilde{H}_{t_nx})},$$
where $\widetilde{H}_{t_nx}$ and $H^*_{t_nx}$ are optimal halfspaces for $\pr$ and $\pr_n$, respectively, as defined in the proof of Lemma~\ref{prop:hdepth-sup-emp-ratio-gl} above.
Note that, by the definition of halfspace depth,
$$ \pr(H^*_{t_nx}) \ge HD(t_nx,\pr) = \pr(\widetilde{H}_{t_nx})$$
Now invoking condition $HD(t_nx,\pr) \ge \gamma_n$, we conclude that $\pr(H^*_{t_nx}) \ge \gamma_n$, which leads to the following upper bound,
$$ \left|\dfrac{HD(t_n x, \pr_n)}{HD(t_n x, \pr)}-1\right| \le 1 - \frac{\pr_n(H^*_{t_nx})}{\pr(H^*_{t_nx})} 
 \le  \sup\left\{ \left| \frac{\pr_n(H)}{\pr(H)}-1\right|:\,\,\, H\in \calH_{t_nx},\,\,\pr(H) \ge \gamma_n\right\}
$$
concluding the proposition. \hfill \qed

{\it Step 3.} This step is based on Theorem 5.1 from \cite{Alexander1987},
which we recall for self containedness in the appendix, so that we have direct access to the conditions under which it holds.

First, note that the collection $\calH$ of all halfspaces is a VC class, so that Theorem 5.1 from \cite{Alexander1987} holds for the halfspaces. 

Next, since the domain of supremum in Proposition~\ref{prop:hdepth-sup-emp-ratio} is a subset of $\calH$, as a consequence of Alexander's theorem, we can conclude that
\begin{eqnarray}\label{eqn:hdepth-fullsup-reducedsup}
\limsup_{n \to \infty} \,\sup \left\{\left|\frac{\pr_n(H)}{\pr(H)}-1\right|: H \in \calH_{t_nx}, \pr(H)\geq \gamma_n\right\} 
&\le & \limsup_{n \to \infty} \,\sup \left\{\left|\frac{\pr_n(H)}{\pr(H)}-1\right|: H \in \calH, \pr(H)\geq \gamma_n\right\} \nonumber\\
& = & 0\,\,\,\,a.s. 
\end{eqnarray}
Combining this last result \eqref{eqn:hdepth-fullsup-reducedsup} with Proposition~\ref{prop:hdepth-sup-emp-ratio} concludes to Theorem \ref{thm:empHDdecay-vs-popHD}.

\subsection{Proof of Theorems~\ref{thm:emp-hdepth-mrv-decay} and \ref{thm:emp-hdepth-light-decay}.}
\label{ss:proofTheosTaildecayHD}

Note that the statement of Theorem~\ref{thm:empHDdecay-vs-popHD} also has certain growth conditions on the sequence $t_n$, which in turn are related to $\gamma_n$. 

\paragraph{Proof of Theorem~\ref{thm:emp-hdepth-mrv-decay}.}~

Observe that
\begin{equation}
     \frac{HD(t_n x,\pr_n)}{V(t_n)} = \frac{HD(t_n x, \pr)}{V(t_n)}\left[1+\left(\dfrac{HD(t_n x, \pr_n)}{HD(t_n x, \pr)}-1\right) \right]
\end{equation}
Now, from Theorem~\ref{thm:hd-decay-heavy-tail}, for large $n$, there exists $c>0$ such that,
$$
HD(t_n x, \pr) > c\, V(t_n).
$$
Combining this last inequality with the given condition $V(t_n)>\gamma_n$ gives, $HD(t_n x, \pr)>c \, \gamma_n$. 

Finally, the result follows from Theorem~\ref{thm:empHDdecay-vs-popHD}. \hfill \qed

\paragraph{Proof of Theorem~\ref{thm:emp-hdepth-light-decay}.}~

\begin{itemize}
    \item[(i)] The proof for the upper bound follows the same line of arguments as for the proof of Theorem~\ref{thm:emp-hdepth-mrv-decay}, together with the result of Theorem~\ref{thm:hd-decay-light-tail}.
    \item[(ii)] The result is straightforward, combining the lower bound given in Theorem~\ref{thm:hd-decay-light-tail}, (ii), 
    with \eqref{eqn:emp-pop-hdepth}, considering a unit vector $x$.
\end{itemize}

\section{Conclusion}
\label{sec:concl}

Much literature has been developed so far on the two types of geometric measures we consider in this paper, geometric quantiles and halfspace depth functions, mainly looking at their properties, such as continuity, convexity, affine equivariance, invariance through (orthogonal) transformations, etc. 
Our focus is on the asymptotic properties of these geometric measures, in particular questioning their relation with the tail behaviour of the underlying distribution. 

First, we have considered the population side, completing the asymptotics literature to get a full picture and laying the basis for the sample side. 
Then we ask the same questions on the asymptotics of the two geometric measures when considering the empirical distribution. This is a very important problem in view of applications, questioning the relevance of those tools when working on samples.   
It is worth recalling that geometric quantiles uniquely identify the underlying probability measure, but this property is not always true for halfspace depth, as recently proven in \cite{Nagy2021a}. Nevertheless, the characterisation is unique when having measures with finite support, as is the case of empirical measures. This motivated us to study further the halfspace depth for samples. 

We were able to provide rates of growth for geometric quantiles and rates of decay for halfspace depth functions when considering the empirical distribution, but also specify these rates depending on the type of tail behaviour of the measure, light or heavy (one of the main questions in risk analysis). These results are important, theoretically, but also in view of providing adequate tools to tackle extremes when analysing data sets. 

The two geometric measures studied in this paper already satisfy some `must have' properties, each one depending on the context, or on the use of it, as there exists no ideal or perfect tool. Nevertheless, we believe that our results contribute to obtain a better idea on these tools, as well as to validate their empirical use depending on the framework. Our next interest is to study another type of geometric measure based on the alternative approach of optimal transport maps. 

\section*{Acknowledgments}

    This paper benefited a lot from the mutual visits at ESSEC Business School, Paris, France and at TIFR-CAM, Bangalore, India, respectively. The authors are grateful to the hosting institutions and for funding to the {\it Fondation des Sciences de la Mod\'elisation} (ANR-11-LABX-0023-01), and SERB--MATRICS grant MTR/2020/000629.
    
\bibliographystyle{chicago}
\bibliography{LitSibsankar.bib}

\newpage
\appendix

\section{Appendix}

\subsection{Higher order expansions for geometric quantiles: sample counterpart}
\label{secSupp:higherOrder-sample}

The proof of Theorem~\ref{thm:numerics} demonstrates that, with additional effort and under suitable conditions on $\alpha_n$ as a function of the sample size $n$, it is possible to derive sample versions of the higher order expansions elicited in  Theorem~\ref{thm:higher-order} $(ii)$.

In this section, we present empirical data illustrating the establishment of the sample version of Equation~\eqref{eqn:3rd-order} from Theorem~\ref{thm:higher-order}. This demonstration will reinforce the idea that the growth pattern of geometric quantiles is intricately tied to the tail behavior of the underlying measure. Consequently, it underscores the utility of geometric quantiles as a tool to distinguish between distributions with light and heavy tails.

Writing $\hat{\Sigma}$ for the sample covariance matrix, set
\begin{align}\label{eqn:3rd-order-appendix}		
        z_n(\alpha_n) &= L_n(\alpha_n) - R_n,  \\
        \text{with} \quad L_n(\alpha_n)& :=
        \|\hat{q}_n(\alpha_n u)\|\Big[\|\hat{q}_n(\alpha_n u)\|^2 (1 - \alpha_n) - \half (\text{tr } \hat{\Sigma}_n - \langle\hat{\Sigma}_n u,u\rangle)\Big] \nonumber\\
        \text{and}\quad R_n &:= - \sum_{k=1}^{d-1}\widehat{\text{cov}}(u, w_k; w_k) +\sum_{k=1}^{d-1} \widehat{\text{cov}}(u; w_k)\, \hat{\mu}_{w_k},\nonumber
\end{align}

where $\hat{\mu}_{w_k}$ denotes the sample mean projected in $w_k$ direction, $\widehat{\text{cov}}(u; w_k)$ the sample covariance between $\langle X,u\rangle$ and $\langle X,w_k\rangle$, and $\widehat{\text{cov}}(u, w_k; w_k)$ the sample covariance between the product $\langle X,u\rangle\langle X,w_k\rangle$ and $\langle X,w_k\rangle$.
Observe that $z_n(\alpha_n)$, as defined in \eqref{eqn:3rd-order-appendix},  contains two distinct sets of terms, one reliant solely on the $n$-sample and another contingent on both the $n$-sample and an additional dependency on $\alpha_n$. 

As we concentrate on the asymptotic behavior of $\hat{q}_n(\alpha_n u)$ when $\alpha_n$ approaches unity (as $n\to\infty$), we plot $z_n(\alpha_n)$ while incrementing both $n$ and $\alpha_n$ simultaneously. It is important to note that the computation regarding the proportional growth of $\alpha_n$ in relation to the sample size $n$ is intricately detailed in the proof of Theorem~\ref{thm:numerics}. Throughout the subsequent plots in this section, we have adopted the following values for $\alpha_n$:
\begin{equation}\label{eqn:alpha-n-plots-appendix}
\alpha_n = 1 - \sqrt{2\frac{\log n}{n}}.
\end{equation}
%
Clearly, the asymptotics of $R_n$ as $n\to\infty$ is a simple consequence of the law of large numbers, whereas the asymptotics of $L_n(\alpha_n)$ (for large $n$ and $\alpha_n$) involves intricate analysis. 

In addition to plotting $z_n(\alpha_n)$, we also plot $L_n(\alpha_n)$ and $R_n$ separately to assess the source of asymptotic fluctuations / aberrations (if any) in the plots corresponding to $z_n(\alpha_n)$.


\noindent\underline{\bf Standard Gaussian}

\begin{figure}[H]
\centering
\begin{minipage}{0.49\textwidth}
\includegraphics[width=.7\linewidth]{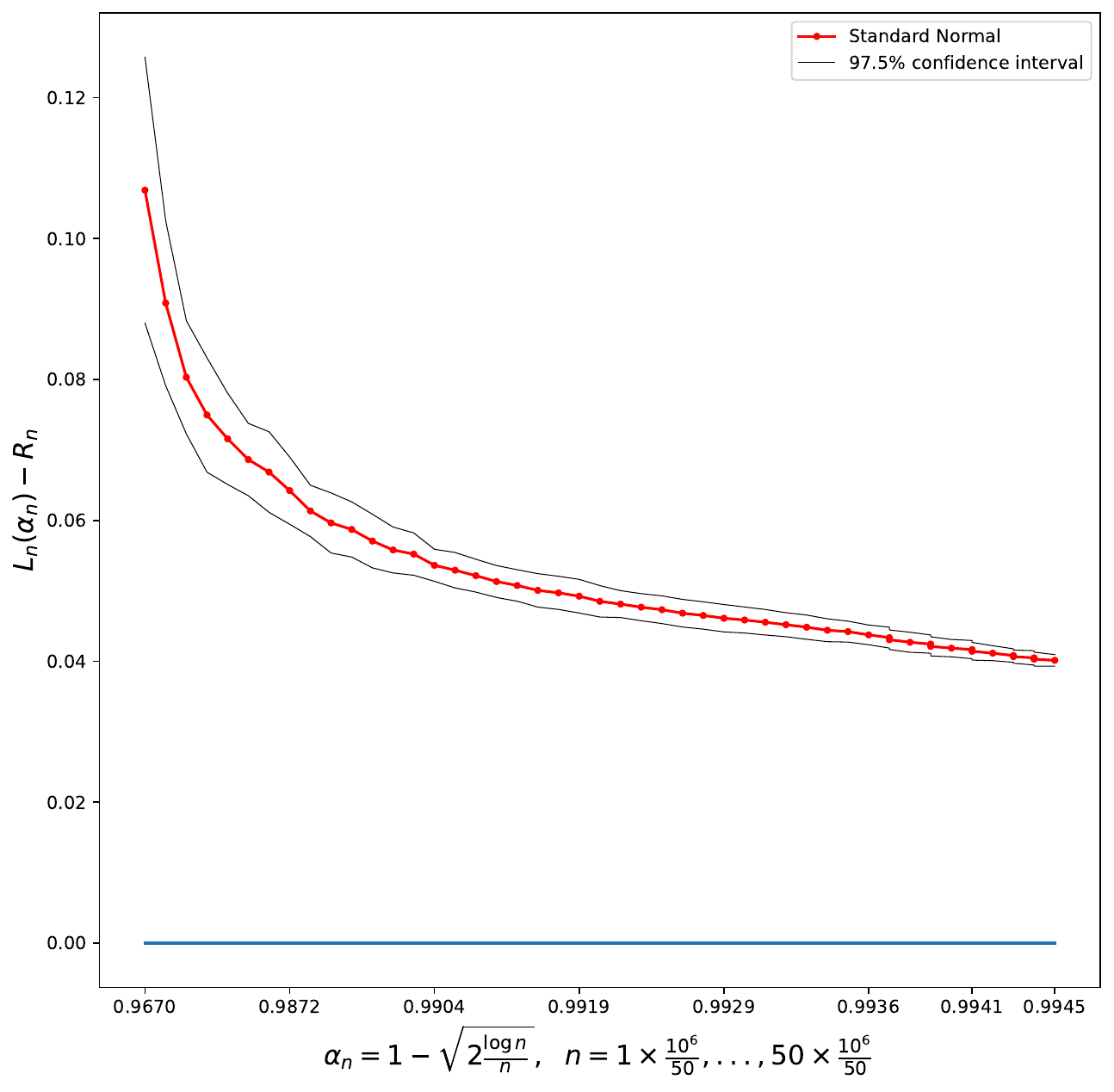}
\end{minipage}
\hfill 
\begin{minipage}{0.49\textwidth}
\includegraphics[width=.7\linewidth]{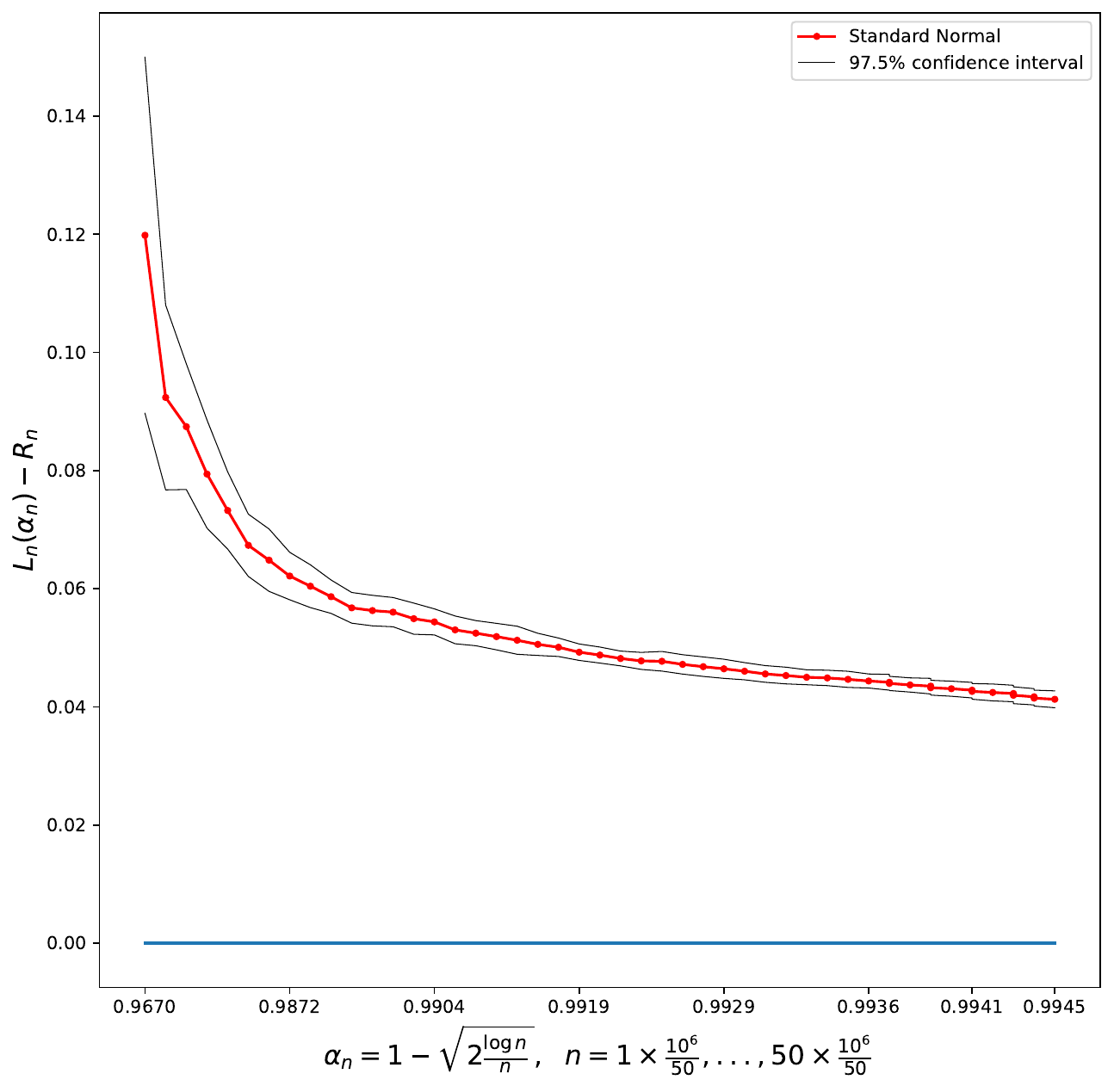}
\end{minipage}
\parbox{420pt}{\caption{\sf \footnotesize Gaussian case. Decay of $z_n(\alpha_n)$ defined in \eqref{eqn:3rd-order-appendix}, along two chosen directions: $u=(1,0)$ (left plot) and $u=(1,1)$ (right plot).}
\label{fig:gauss-2nd-order-zn}}
\end{figure}

The following plots show the decay of $z_n(\alpha_n)$, along two chosen directions $(1,0)$ and $(1,1)$, as $n$ increases linearly and $\alpha_n$ increases according to \eqref{eqn:alpha-n-plots-appendix}. Specifically, $10^6$ samples were generated from bivariate standard Gaussian distribution. Subsequently,  $z_n(\alpha_n)$ was iteratively computed for $n= \frac{k}{50}\,10^6$, for $k=1,\ldots,50$. To mitigate fluctuations induced by the samples, these computations were redone using jackknife resampling. The averages of all computed $z_n(\alpha_n)$ values are then plotted in Figure~\ref{fig:gauss-2nd-order-zn}.
Observe that not only the confidence intervals are diminishing as $n$ and $\alpha_n$ increase, but also the value of $z_n(\alpha_n)$ distinctly approaches $0$. 

\begin{figure}[H]
\centering
\begin{minipage}{0.49\textwidth}
\includegraphics[width=.7\linewidth]{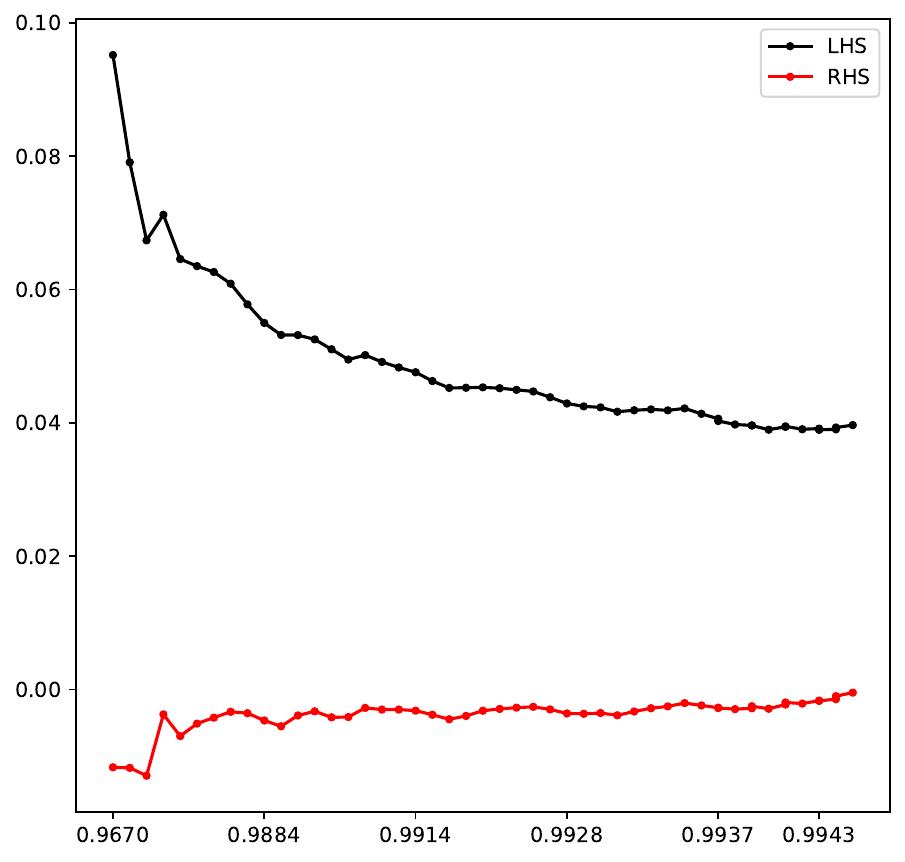}
\end{minipage}
\hfill 
\begin{minipage}{0.49\textwidth}
\includegraphics[width=.7\linewidth]{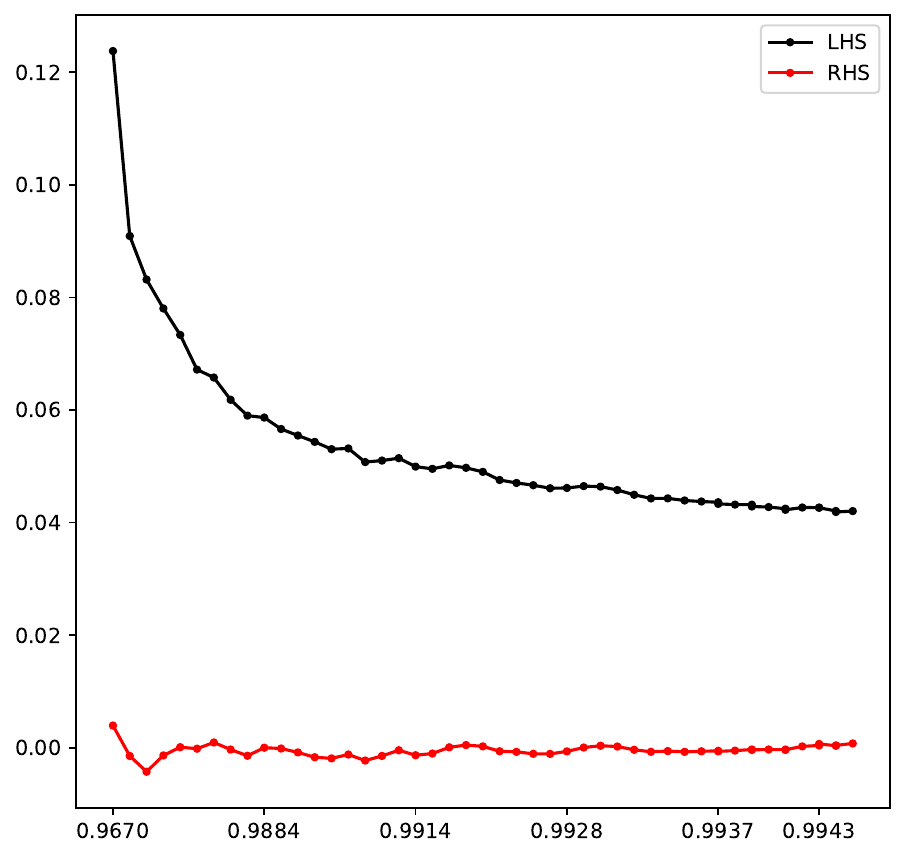}
\end{minipage}
\parbox{460pt}{\caption{\sf \small Gaussian case. $L_n(\alpha_n)$ and $R_n$, denoted as LHS and RHS, respectively, are represented on a same plot. Two directions are considered: $u=(1,0)$ (left plot), and $u=(1,1)$ (right plot).}
\label{fig:gauss-2nd-order-Ln-Rn}}
\end{figure}

We shall now plot $L_n(\alpha_n)$ and $R_n$ separately to observe whether one or the other presents a distinctive pattern. Considering that the samples are drawn from bivariate standard Gaussian distribution, a straightforward consequence of the law of large numbers suggests convergence of $R_n$ towards $0$, which clearly appears in Figure~\ref{fig:gauss-2nd-order-Ln-Rn}. Consequently, if $z_n(\alpha_n) \to 0$ as $n\to\infty$, it logically follows that $L_n(\alpha_n)$ must also converge to $0$. In Figure~\ref{fig:gauss-2nd-order-Ln-Rn}, we observe that the gradual approach of $z_n(\alpha_n)$ to zero predominantly arises from the slow convergence of $L_n(\alpha_n)$ towards zero.
%


\noindent\underline{\bf Multivariate Student-$t$ with $\nu=4.2$}

We now redo the above plots for samples generated from a bivariate $t$--distribution whose density is given by
\begin{equation}\label{eqn:mul-t-density}
\frac{\Gamma(1 +\nu/2)}{\Gamma(\nu/2)\,\nu\,\pi} \left[ 1 + \frac{\|x\|^2}{\nu}\right]^{-\left( 1+ \frac{\nu}2\right)}
\end{equation}
starting with $\nu=4.2$.

\begin{figure}[H]
\centering
\begin{minipage}{0.49\textwidth}
\includegraphics[width=.7\linewidth]{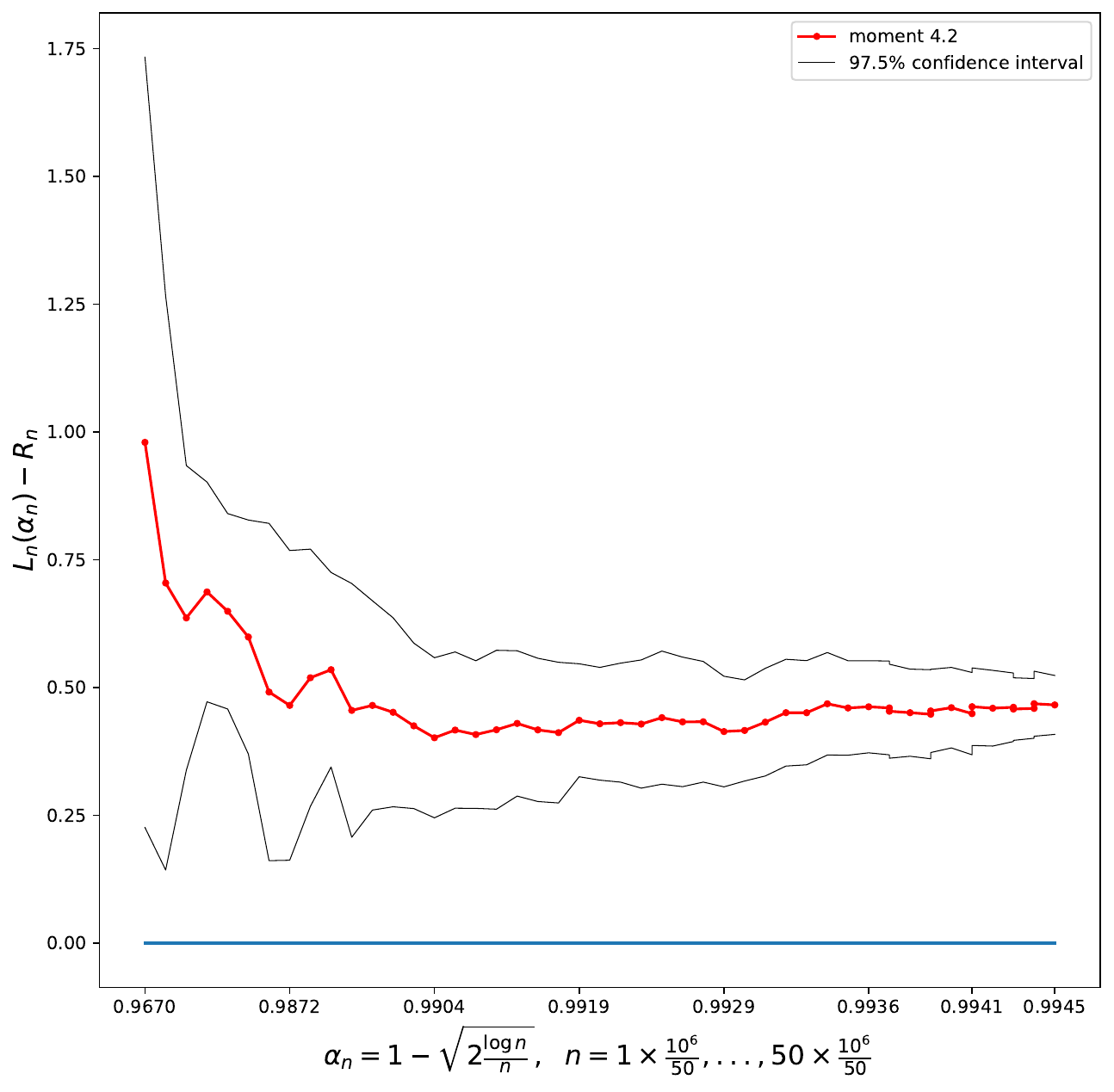}
\end{minipage}
\hfill 
\begin{minipage}{0.49\textwidth}
\includegraphics[width=.7\linewidth]{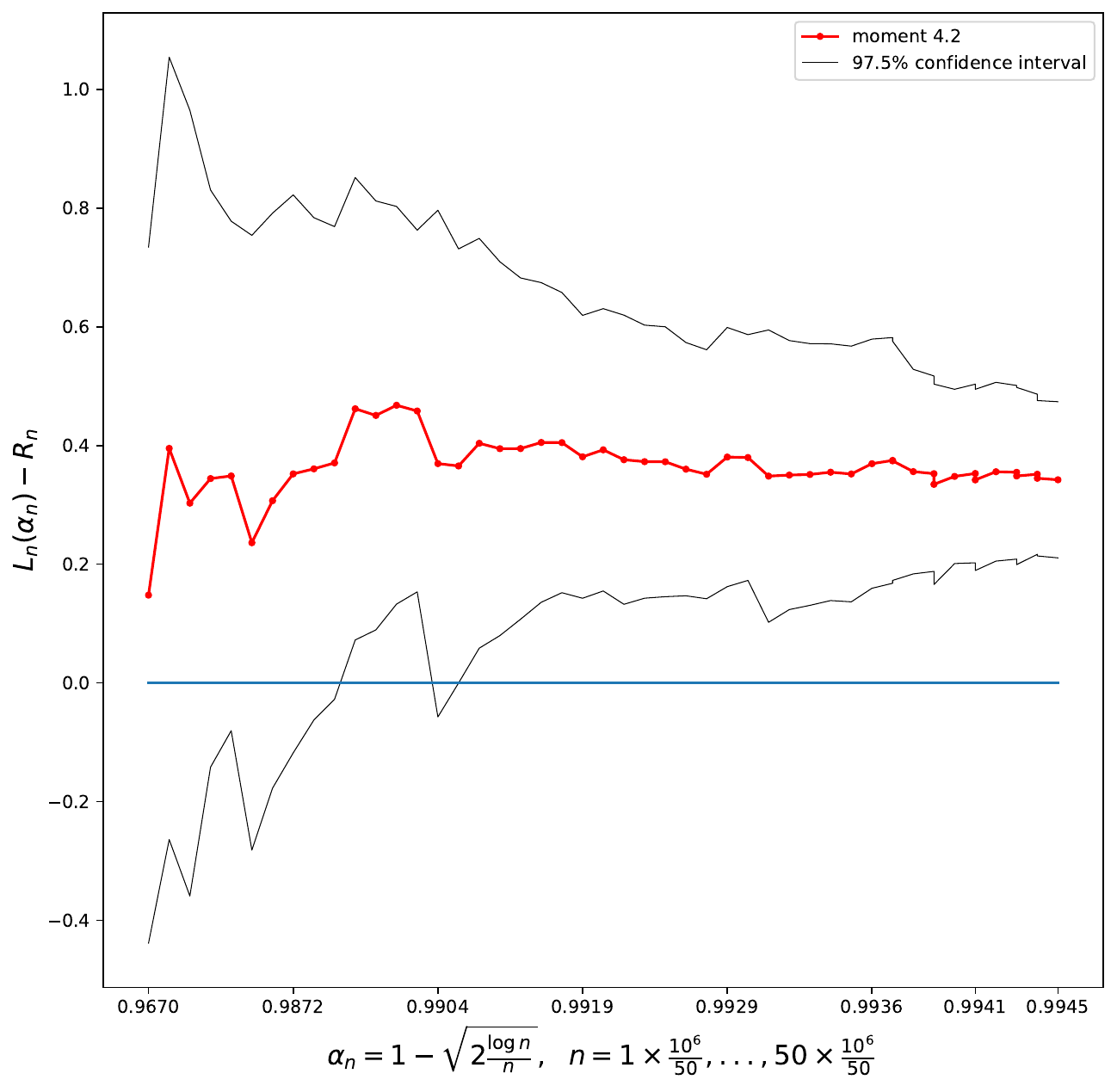}
\end{minipage}
\parbox{460pt}{\caption{\sf \small Student-$t$ case with $\nu=4.2$. Decay of $z_n(\alpha_n)$ defined in \eqref{eqn:3rd-order-appendix}, along two chosen directions: $u=(1,0)$ (left plot) and $u=(1,1)$ (right plot).}
\label{fig:t-4-2nd-order-zn}}
\end{figure}

After generating $10^6$ samples from the above distribution, we computed $z_n(\alpha_n)$ iteratively for $n= \frac{k}{50}\,10^6$, for $k=1,\ldots,50$. Employing jackknife resampling, as previously done, we smoothed out fluctuations induced by the samples. Subsequently, we computed the average of all derived $z_n(\alpha_n)$ values, which are depicted in Figure~\ref{fig:t-4-2nd-order-zn}. These plots showcase the decay of $z_n(\alpha_n)$ along two specified directions, $(1,0)$ and $(1,1)$, as $n$ progresses linearly while $\alpha_n$ follows \eqref{eqn:alpha-n-plots-appendix}.

Observe that $z_n(\alpha_n)$ is getting closer to $0$ as $n$ and $\alpha_n$ increase, however, the order of the magnitude of the difference is higher than in the case of standard Gaussian, with a multiplicative factor of order 10. Also note that the confidence bands are broader in comparison to those observed in the standard Gaussian scenario. 

\begin{figure}[H]
\centering
\begin{minipage}{0.49\textwidth}
\includegraphics[width=.7\linewidth]{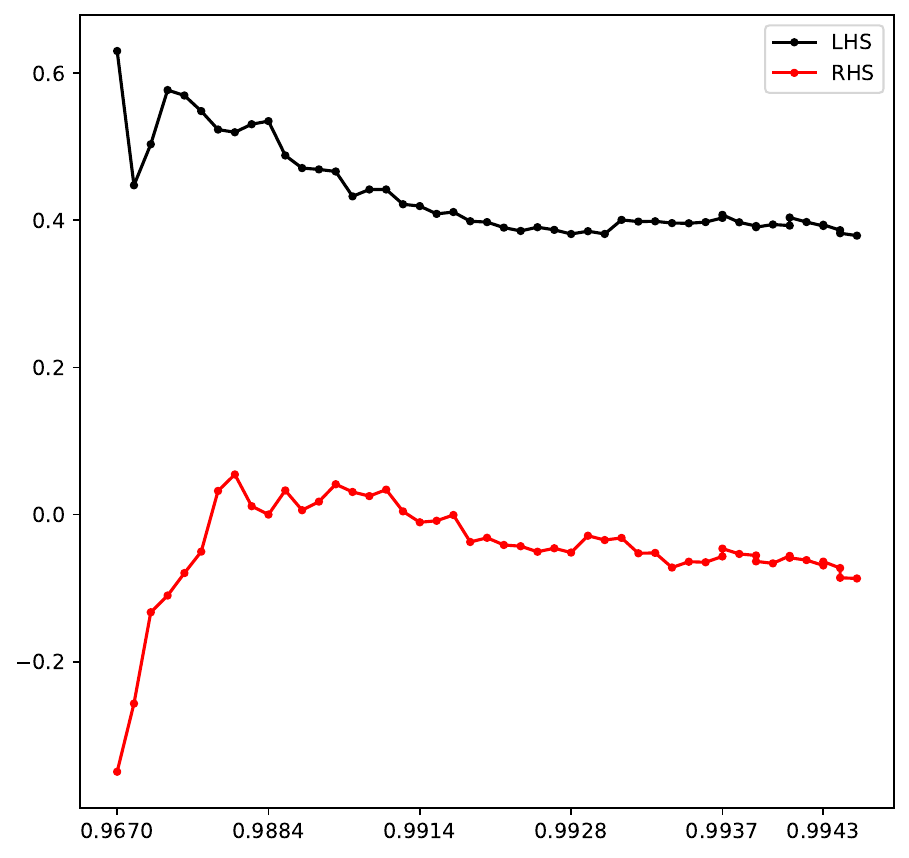}
\end{minipage}
\hfill 
\begin{minipage}{0.49\textwidth}
\includegraphics[width=.7\linewidth]{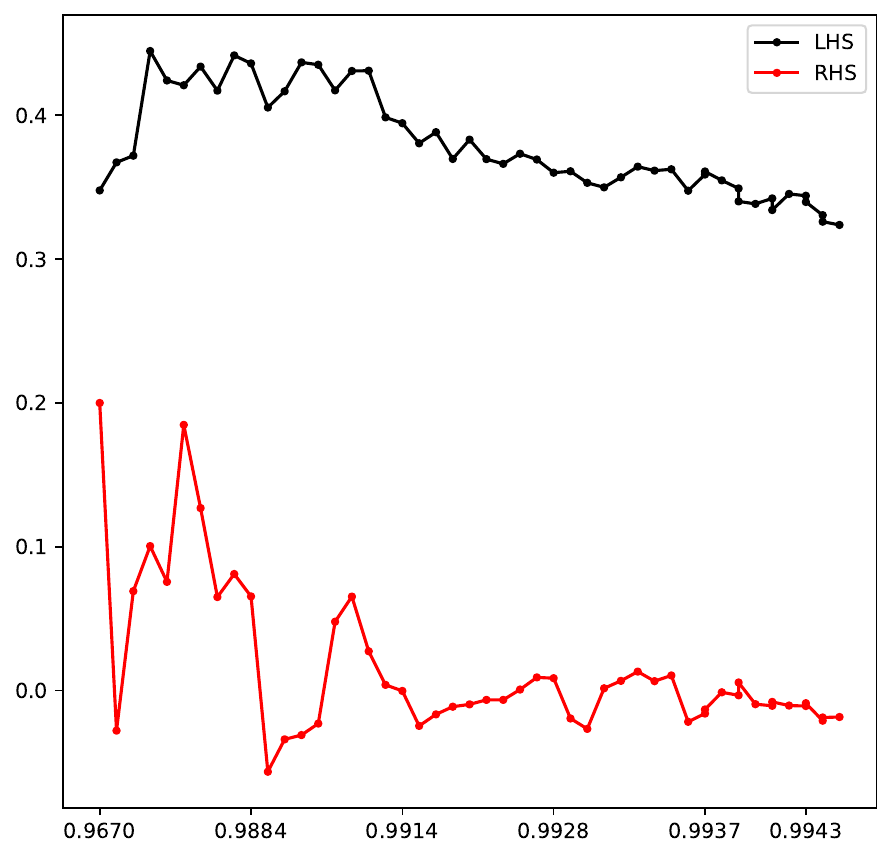}
\end{minipage}
\parbox{460pt}{\caption{\sf \small Student-$t$ case with $\nu=4.2$. $L_n(\alpha_n)$ and $R_n$, denoted as LHS and RHS, respectively, are represented on a same plot. Two directions are considered: $u=(1,0)$ (left plot), and $u=(1,1)$ (right plot).}
\label{fig:t-4-2nd-order-Ln-Rn}}
\end{figure}

In order to further investigate the reason for slower rate of convergence of $z_n(\alpha_n)$ to $0$, we plot $L_n(\alpha_n)$ and $R_n$ separately in Figure~\ref{fig:t-4-2nd-order-Ln-Rn}.
Note that, given the choice of $\nu$, the fourth order moments (hence the third ones) of the density given by \eqref{eqn:mul-t-density} are finite. It implies that $R_n$, by means of the law of large numbers, will converge to a finite limit, which, due to the symmetry of the underlying distribution, equals $0$. As a consequence, $L_n(\alpha_n)$ must also converge to $0$. We observe in Figures~\ref{fig:t-4-2nd-order-zn} and \ref{fig:t-4-2nd-order-Ln-Rn} that $R_n$ appears to stabilise around its eventual limit. On the other hand, $L_n(\alpha_n)$ appears to approach $0$ rather slowly, in comparison to the standard Gaussian scenario. Such relatively slow decay of $z_n(\alpha_n)$ (and of $L_n(\alpha_n)$) as compared to the case of standard Gaussian, might be attributed to the relatively heavier tail of the $t$--distribution.


\noindent\underline{\bf Multivariate $t$ with $\nu=3.2$}

Further exploring the impact of the tail behaviour of the underlying distribution on the asymptotics of $z_n(\alpha_n)$, we now consider samples from the density given by \eqref{eqn:mul-t-density} with $\nu=3.2$. Note that with this choice of $\nu$, the third order moments of this distribution remain finite (but not the fourth ones).

\begin{figure}[H]
\centering
\begin{minipage}{0.49\textwidth}
\includegraphics[width=.7\linewidth]{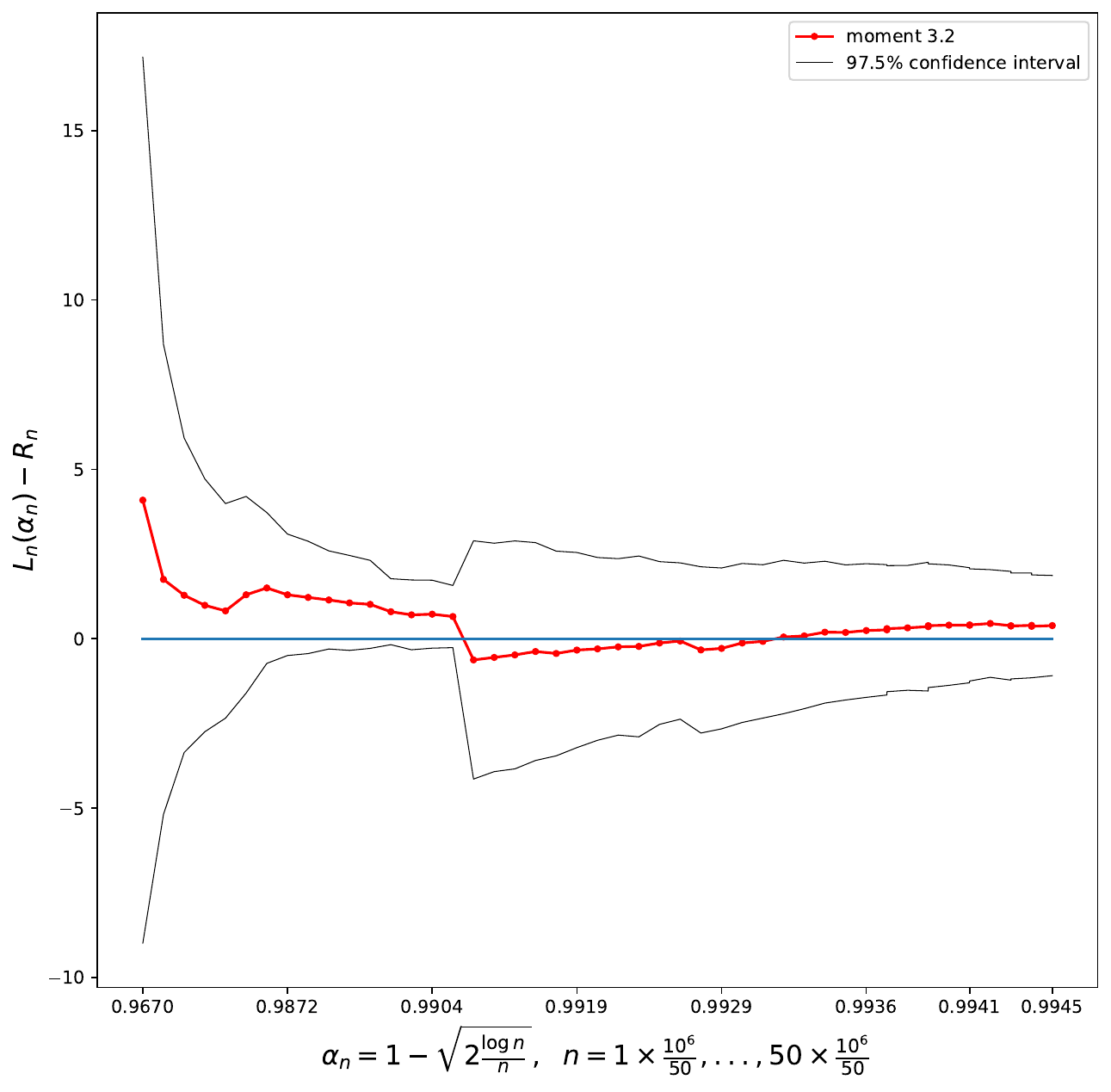}
\end{minipage}
\hfill 
\begin{minipage}{0.49\textwidth}
\includegraphics[width=.7\linewidth]{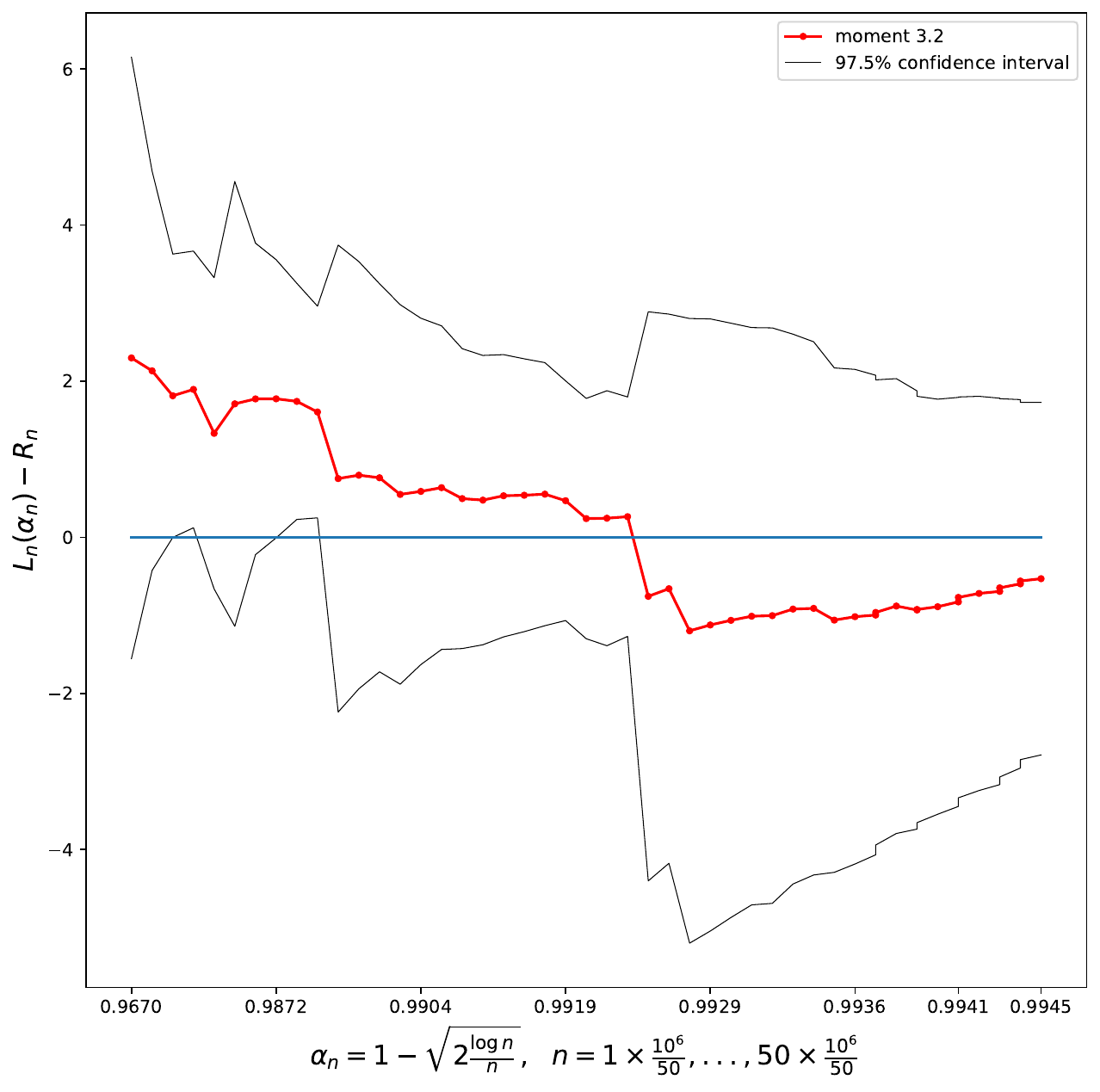}
\end{minipage}
\parbox{460pt}{\caption{\sf \small Student-$t$ case with $\nu=3.2$. Decay of $z_n(\alpha_n)$ defined in \eqref{eqn:3rd-order-appendix}, along two chosen directions: $u=(1,0)$ (left plot) and $u=(1,1)$ (right plot).}
\label{fig:t-3-2nd-order-zn}}
\end{figure}

\begin{figure}[H]
\centering
\begin{minipage}{0.49\textwidth}
\includegraphics[width=.7\linewidth]{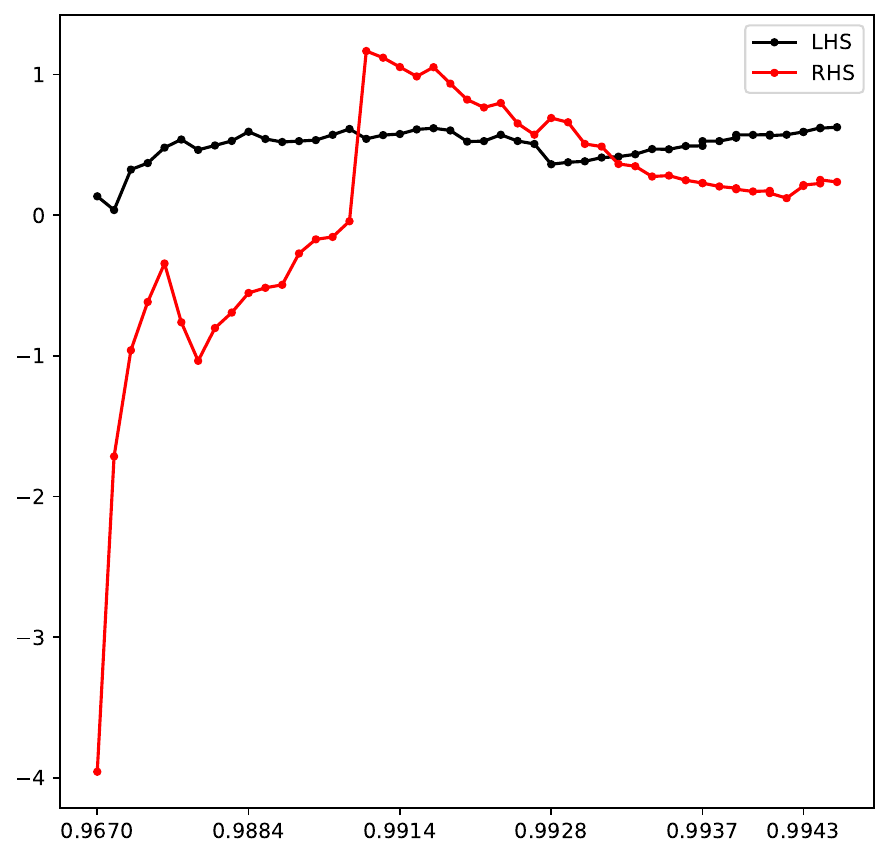}
\end{minipage}
\hfill 
\begin{minipage}{0.49\textwidth}
\includegraphics[width=.7\linewidth]{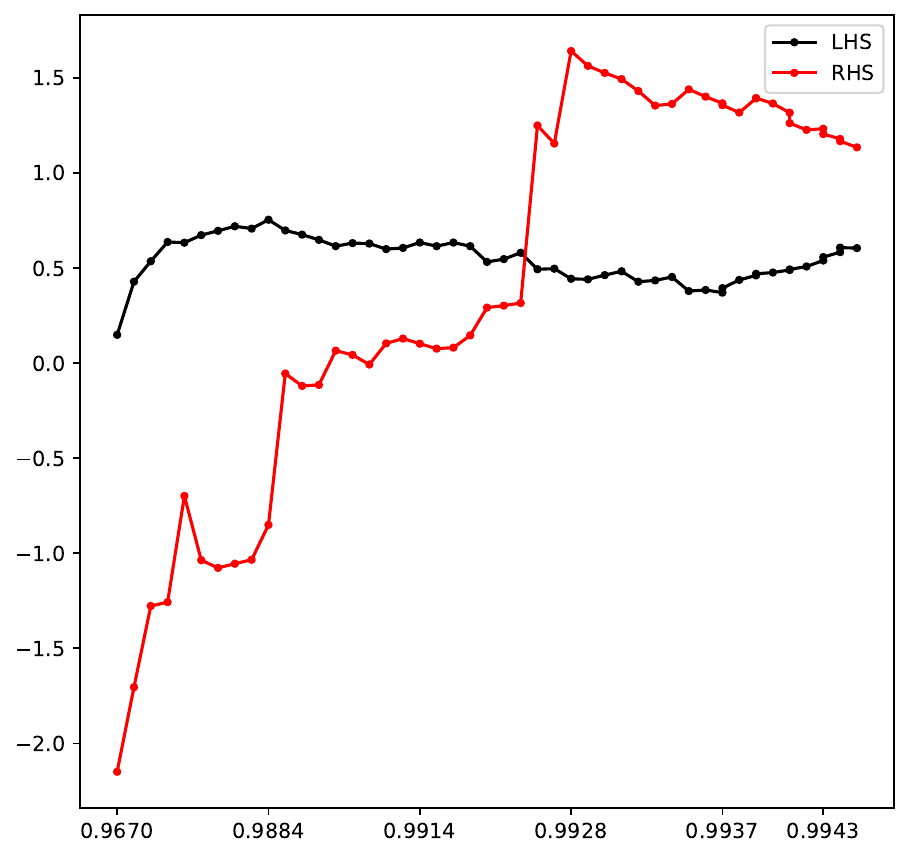}
\end{minipage}
\parbox{460pt}{\caption{\sf \small Student-$t$ case with $\nu=3.2$. $L_n(\alpha_n)$ and $R_n$, denoted as LHS and RHS, respectively, are represented on a same plot. Two directions are considered: $u=(1,0)$ (left plot), and $u=(1,1)$ (right plot).}}
\label{fig:t-3-2nd-order-Ln-Rn}
\end{figure}

After generating $10^6$ samples from this distribution, we repeated the entire process as explained above in the cases of the standard Gaussian and multivariate $t$-distribution with $\nu=4.2$. Subsequently, we plot the resulting average $z_n(\alpha_n)$ in Figure~\ref{fig:t-3-2nd-order-zn}.
Evidently, the convergence of $z_n(\alpha_n)$ to $0$ is much slower in comparison to the previous two scenarios. The order of the magnitude of $z_n(\alpha_n)$ is significantly different from what was computed in earlier examples of standard Gaussian and multivariate $t$ with $\nu=4.2$ (with a factor of order 10 (or 6) between the two Student-$t$'s). Indeed, it is noticeable that the confidence bands tend to widen as the tail of the distribution becomes heavier. To identify the source of the slower rate of convergence, $L_n(\alpha_n)$ and $R_n$ are plotted in Figure~\ref{fig:t-3-2nd-order-Ln-Rn}.
%
Notice first that $L_n(\alpha_n)$ and $R_n$ tend to concentrate around each other as $n$ and $\alpha_n$ increase, but we cannot conclusively determine the convergence. Interestingly, $L_n(\alpha_n)$, which involves computing $\hat{q}_n(\alpha_n u)$, shows a steady profile, but does not conclusively exhibit {\it convergence} to $0$. In contrast, despite exhibiting large jumps, $R_n$ appears to show tendency to be closer to $0$. Observe that the size of the jumps is almost of the order of the value of $R_n$ itself. Such behaviour is characteristic of heavy tail distributions, thus conclusively indicating the impact of the tail distribution on the convergence, or lack thereof.


\noindent\underline{\bf Multivariate $t$ with $\nu=2.2$}

Continuing on the same theme to explore the convergence of $z_n(\alpha_n)$ with restrictive moment conditions, we now repeat the above computations for multivariate -- $t$ as given by \eqref{eqn:mul-t-density} with $\nu=2.2$. The output is shown below in Figures~\ref{fig:t-2-2nd-order-zn} and \ref{fig:t-2-2nd-order-Ln-Rn}.
\begin{figure}[H]
\centering
\begin{minipage}{0.49\textwidth}
\includegraphics[width=.7\linewidth]{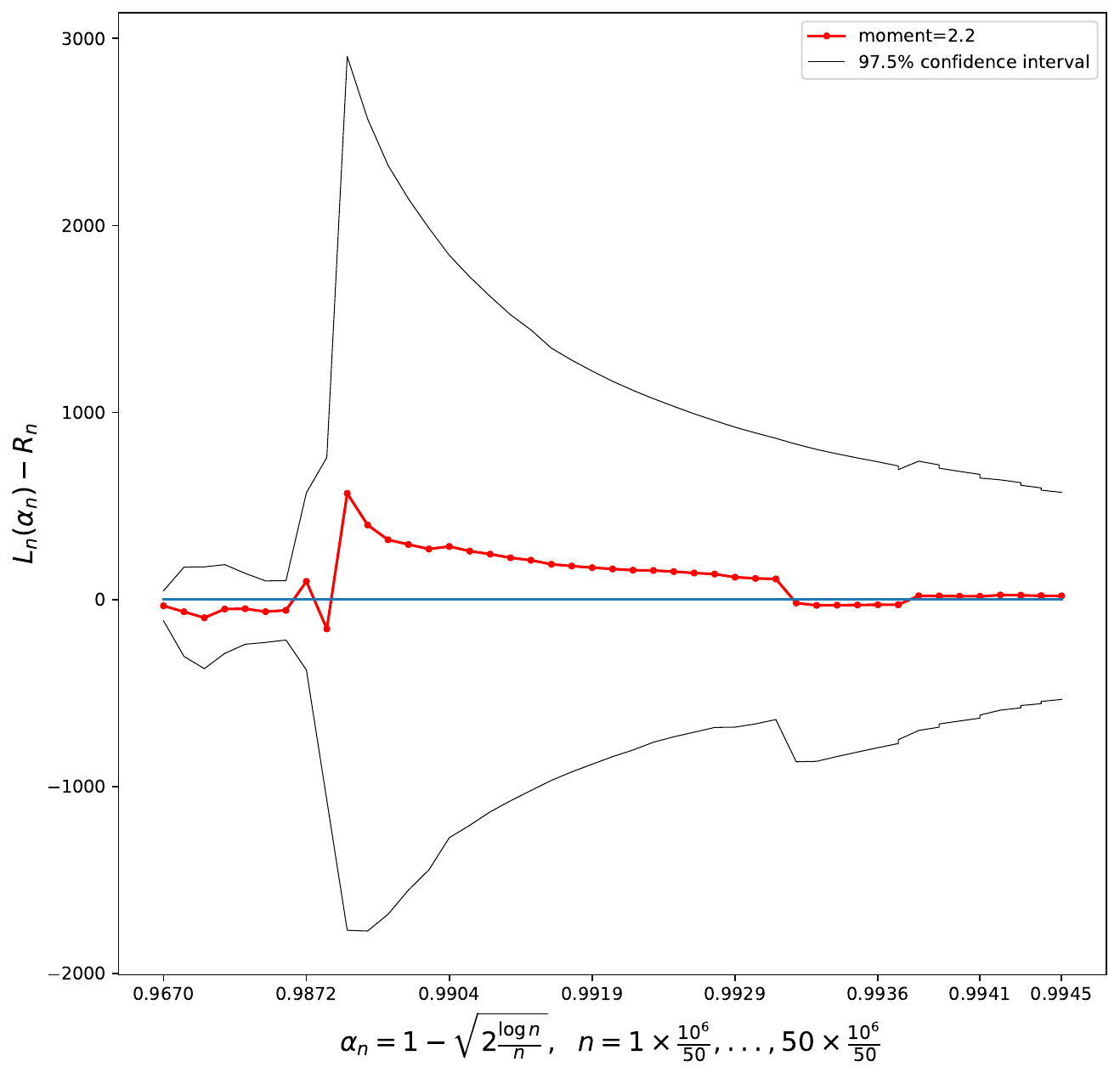}
\end{minipage}
\hfill 
\begin{minipage}{0.49\textwidth}
\includegraphics[width=.7\linewidth]{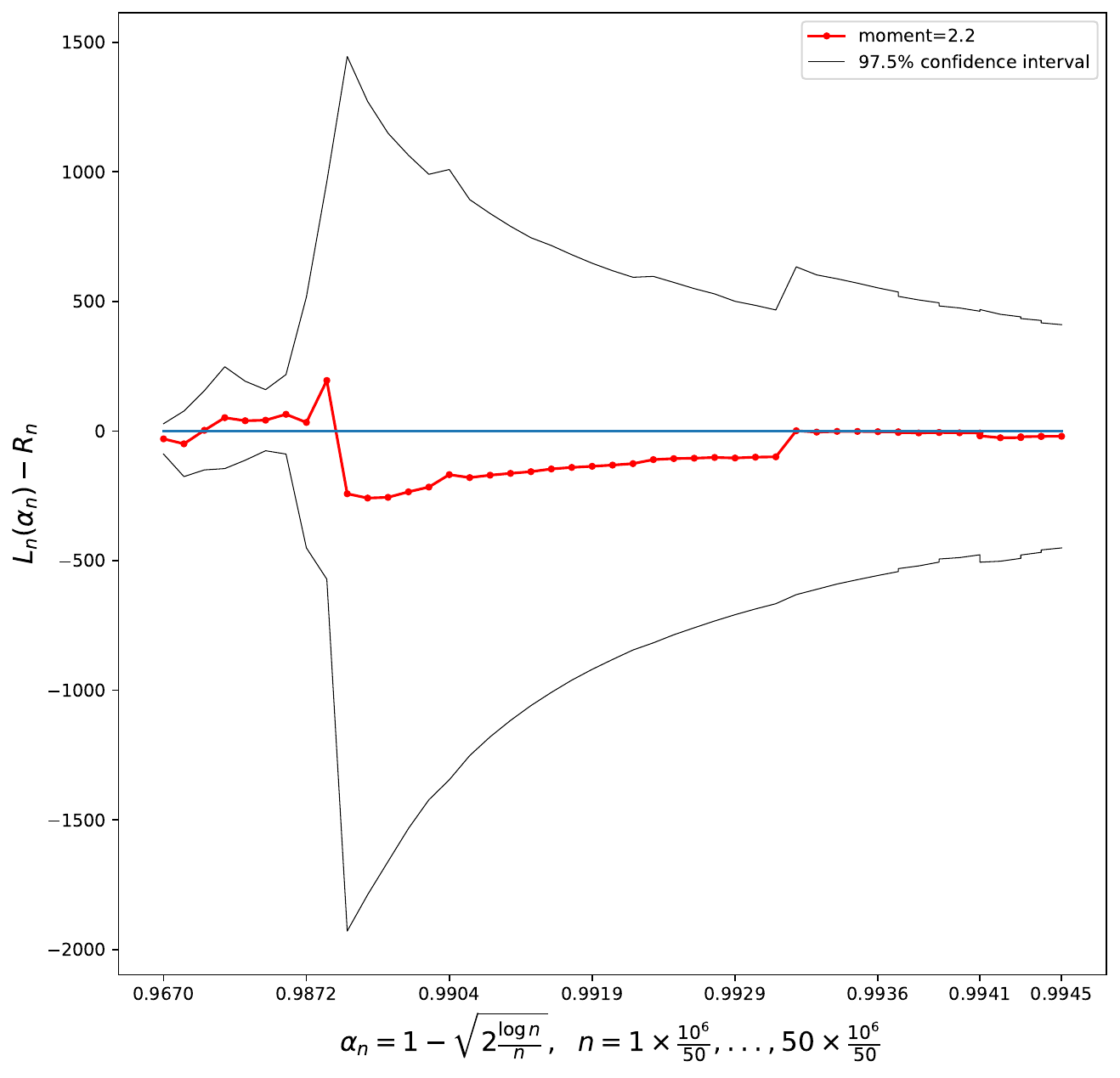}
\end{minipage}
\parbox{460pt}{\caption{\sf \small Student-$t$ case with $\nu=2.2$. Decay of $z_n(\alpha_n)$ defined in \eqref{eqn:3rd-order-appendix}, along two chosen directions: $u=(1,0)$ (left plot) and $u=(1,1)$ (right plot).}
\label{fig:t-2-2nd-order-zn}}
\end{figure}
\begin{figure}[H]
\centering
\begin{minipage}{0.49\textwidth}
\includegraphics[width=.7\linewidth]{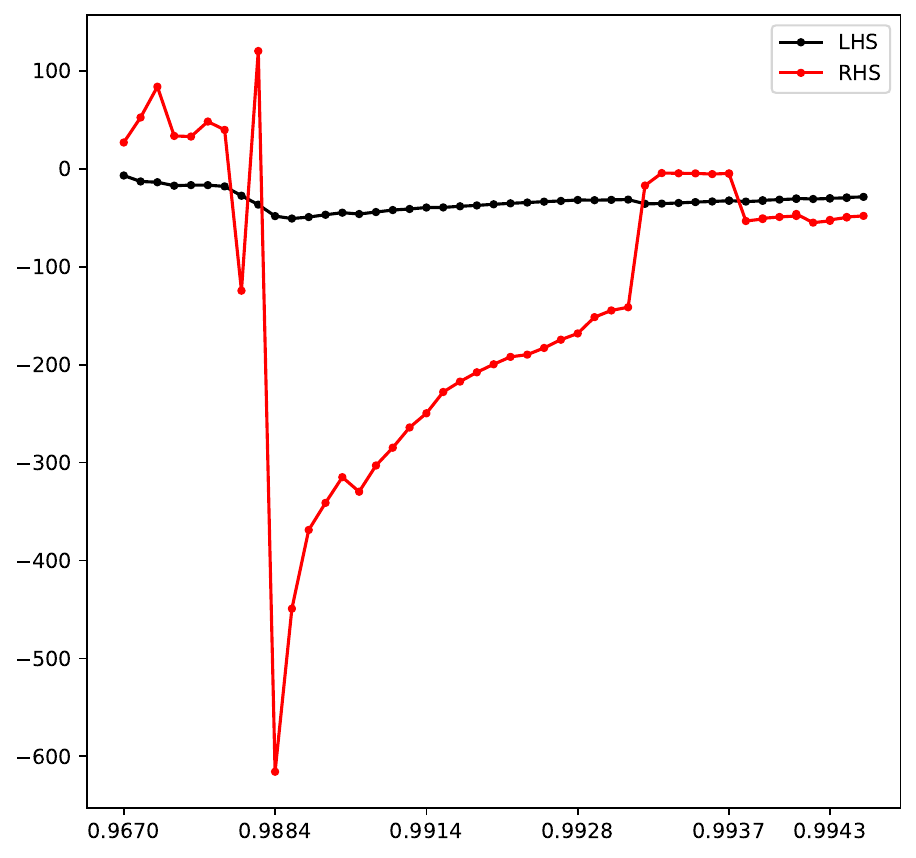}
\end{minipage}
\hfill 
\begin{minipage}{0.49\textwidth}
\includegraphics[width=.7\linewidth]{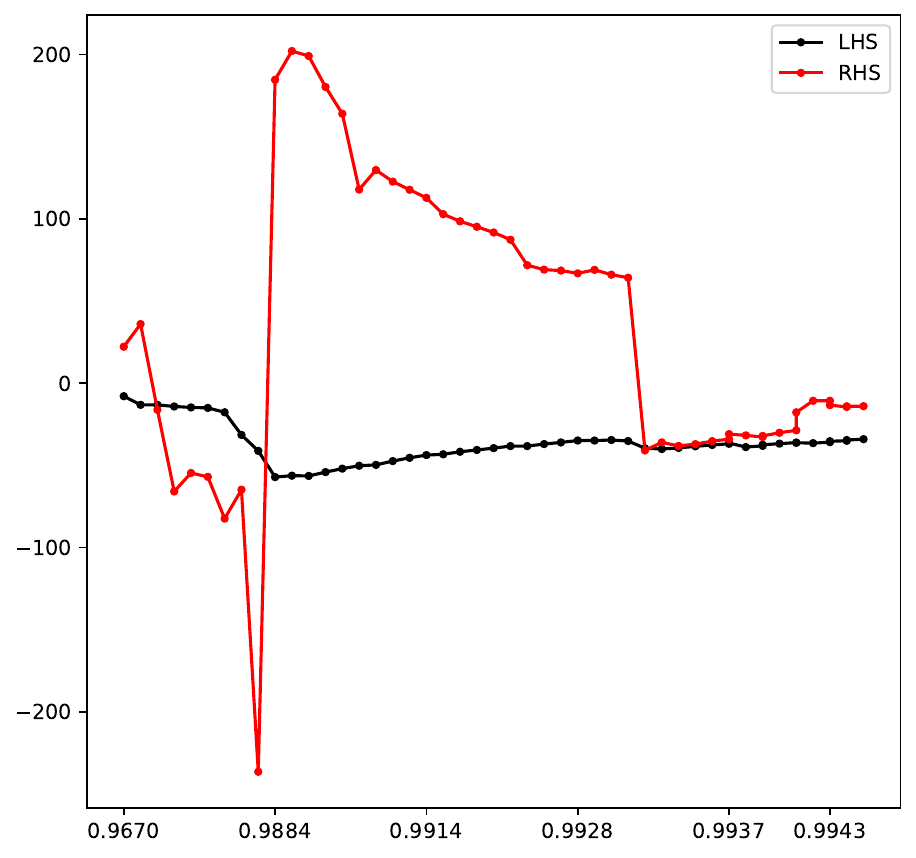}
\end{minipage}
\parbox{460pt}{\caption{\sf \small Student-$t$ case with $\nu=2.2$. $L_n(\alpha_n)$ and $R_n$, denoted as LHS and RHS, respectively, are represented on a same plot. Two directions are considered: $u=(1,0)$ (left plot), and $u=(1,1)$ (right plot).}}
\label{fig:t-2-2nd-order-Ln-Rn}
\end{figure}
The expression in $z_n(\alpha_n)$ involves third order empirical moments, which are always finite and can be computed for any given $n$. Nevertheless, as $n\to\infty$, the convergence of $z_n(\alpha_n)$ is challenged due to the absence of moments of order larger or equal to 3 in the case of $\nu=2.2$.
The behaviour of $R_n$ with several jumps is consistent with our previous observation in the case of $\nu=3.2$, as is the behaviour of $L_n(\alpha_n)$, considering again a multiplicative factor of order about 10. Nevertheless, notice the difference between $L_n(\alpha_n)$ and $R_n$, significantly of larger order.


\noindent\underline{\bf Summary}

We illustrated the sample version of Equation~\eqref{eqn:3rd-order} from Theorem~\ref{thm:higher-order} by plotting the decay of $z_n(\alpha_n)$ defined in \eqref{eqn:3rd-order-appendix} and of its additive components $L_n(\alpha_n)$ and $R_n$, using simulation experiments from light to heavy-tailed distributions, varying also the direction.

We can observe that the behavior of $z_n(\alpha_n)$ is intricately tied to the tail behavior of the underlying measure: The lighter the tail, the faster the convergence of $z_n(\alpha_n)$ towards 0. 

This property suggests that geometric quantiles may be a useful tool to distinguish between distributions with light and heavy tails, contrarily to earlier results based on second moments only.

\subsection{Alexander Theorem}
For the paper to be self-contained, let us recall the following theorem:

{\sf Theorem 5.1 in \cite{Alexander1987}}

    Let $\pr$ be a probability measure defined on $(\real^d,\calB(\real^d))$, and let $g$ be capacity function corresponding to $\pr$ as defined in \cite{Alexander1987}. Consider a sequence $\gamma_n$ satisfying the following conditions:
    $$    
    n^{-1}\log (g(\gamma_n)) = {\textrm o}(\gamma_n) \qquad \text{and}\qquad n^{-1}\log \log n = {\textrm o}(\gamma_n).
    $$
    Then,
    \begin{equation}
    \limsup_{n \to \infty} \,\,\,\sup \left\{\left|\frac{\pr_n(C)}{\pr(C)}-1\right|: C \in \calS, \pr(C)\geq \gamma_n\right\}=0\,\,\,\,a.s.
    \end{equation}

\end{document}